\documentclass[onecolumn,final]{elsarticle}

\usepackage{type1cm}          % activate if the above 3 fonts are
\usepackage{makeidx}          % allows index generation
\usepackage{graphicx}         % standard LaTeX graphics tool
\usepackage{multicol}         % used for the two-column index
\usepackage[bottom]{footmisc} % places footnotes at page bottom

\usepackage{newtxtext}        % 
\usepackage{newtxmath,bm,bbm} % selects Times Roman as basic font
\usepackage{setspace}
\usepackage{placeins,subfigure,overpic,rotating,subfigure}
\usepackage{braket}
\usepackage{chngpage}
\usepackage{cases}
\usepackage{tabularx}
\usepackage{multirow}

\usepackage{lettrine}
\usepackage{booktabs}

\usepackage{setspace,comment}

\makeindex                    % used for the subject index
\newtheorem{theorem}{Theorem}[section]

\newtheorem{remark}[theorem]{Remark}

\newcommand{\Eh} {\mathcal{E}_h}

\newcommand{\vrtx}{\mathsf{v}}
\newcommand{\normV} [1]{\left|\!\left|#1\right|\!\right|_{\V}}

\newcommand{\PiPr}[1]{\Pi^{\nabla,\P}_{#1}} %% elliptic projector, local
\newcommand{\Pizr}[1]{\Pi^{0,\P}_{#1}}  %% orthogonal projector, local
\newcommand{\PiDr}[1]{\Pi^{\Delta,\P}_{#1}} %% elliptic projector for CH, local

\def\trait #1 #2 #3 {\vrule width #1pt height #2pt depth #3pt}
\def\fin{\hfill
        \trait .3 5 0
        \trait 5 .3 0
        \kern-5pt
        \trait 5 5 -4.7
        \trait 0.3 5 0
\medskip}

\newcommand{\symmSpace}{\mathbb{R}^{2\times2}_{\textrm{sym}}} 
 \def\Vhcahn{{W}_h}
\def\VE{W_{h|\P}}
\def\VEt{\widetilde{V}_{h|\P}}
\def\P{{\mathbb P}}
\def\Vhbc{{W}_h^0}
 \newcommand{\nx}{n_x}
\newcommand{\ny}{n_y}
\newcommand{\tx}{\tau_x}
\newcommand{\ty}{\tau_y}
\newcommand{\btau}{\bm\tau}

\newcommand{\Vhrp} [1]{V^{p}_{h,#1}}
\newcommand{\VhPrp}[1]{V^{p}_{h,#1}(\P)}

\newcommand{\tVhPrp}[1]{\tilde{V}^{p}_{h,#1}(\P)}

\newcommand{\INTP}{\footnotesize{I}}

\newcommand{\REAL}{\mathbbm{R}}
\newcommand{\INTG}{\mathbbm{N}}

\newcommand{\TERM}[1]{\textbf{(#1)}}

\newcommand{\fv}{\mathbf{f}}
\newcommand{\gv}{\mathbf{g}}

\newcommand{\nv}{\mathbf{n}}

\newcommand{\uv}{\mathbf{u}}
\renewcommand{\vv}{\mathbf{v}}
\newcommand{\wv}{\mathbf{w}}
\newcommand{\xv}{\mathbf{x}}
\newcommand{\yv}{\mathbf{y}}

\newcommand{\Vv}{\mathbf{V}}

\newcommand{\as}{a}

\newcommand{\cs}{c}

\newcommand{\fs}{f}

\newcommand{\ms}{m}

\newcommand{\ps}{p}
\newcommand{\qs}{q}
\newcommand{\rs}{r}
\renewcommand{\ss}{s}
\newcommand{\ts}{t}
\newcommand{\us}{u}
\newcommand{\vs}{v}
\newcommand{\ws}{w}
\newcommand{\xs}{x}
\newcommand{\ys}{y}
\newcommand{\zs}{z}

\newcommand{\Cs}{C}
\newcommand{\Ds}{D}

\newcommand{\Fs}{F}

\newcommand{\Ms}{M}
\newcommand{\Ns}{N}

\newcommand{\Ts}{T}

\newcommand{\Vs}{V}

\newcommand{\Vsp}{V^{\prime}}

\newcommand{\matD}{\mathsf{D}}

\newcommand{\matI}{\mathsf{I}}

\newcommand{\calD}{\mathcal{D}}
\newcommand{\calE}{\mathcal{E}}

\newcommand{\calH}{\mathcal{H}}

\newcommand{\calM}{\mathcal{M}}

\newcommand{\calT}{\mathcal{T}}

\newcommand{\calV}{\mathcal{V}}

\newcommand{\HONE}  {H^1}
\newcommand{\HONEzr}{H^1_0}

\newcommand{\HTWO}  {H^2}
\newcommand{\HONEgm}{H^1_{\Gamma_D}}

\newcommand{\LTWO}  {L^2}
\newcommand{\LINF}  {L^{\infty}}

\newcommand{\HS}[1] {H^{#1}}

\newcommand{\CS}[1] {C^{#1}}
\newcommand{\VS}[1] {V^{#1}}

\newcommand{\PS}[1] {\mathbbm{P}_{#1}}

\renewcommand{\P} {\textsf{P}}            % polyhedral element
\newcommand  {\E} {e}% {\textsf{e}}            % edge
\newcommand  {\V} {\mathsf{v}}% {\textsf{v}}            % vertex
\newcommand{\hh}{h}
\newcommand{\Th}{\Omega_{\hh}}

\newcommand{\DIM} {d}              % space dimension
\newcommand{\hP}{\hh_{\P}}

\newcommand{\hE}{\hh_{\E}}
\newcommand{\hV}{\hh_{\V}}

\newcommand{\mP}{\ABS{\P}}

\newcommand{\mE}{\ABS{\E}}

\newcommand{\Pset}{\mathcal{P}}    % mesh polyhedrons
\newcommand{\NMB}{N}
\newcommand{\NP}{\NMB^{\Pset}}   % cardinality of the set of polyhedrons
\newcommand{\dV}{\,dV}
\newcommand{\dS}{\,ds}
\newcommand{\dx}{\,d\xv}

\newcommand{\fsh}{\fs_{\hh}}

\newcommand{\qsh}{\qs_{\hh}}

\newcommand{\ush}{\us_{\hh}}
\newcommand{\vsh}{\vs_{\hh}}

\newcommand{\wsh}{\ws_{\hh}}

\newcommand{\Fsh}{\Fs_{\hh}}

\newcommand{\uvh}{\uv_{\hh}}

\newcommand{\vvh}{\vv_{\hh}}

\newcommand{\wvh}{\wv_{\hh}}

\newcommand{\asP}{\as_{\P}}

\newcommand{\msP}{\ms_{\P}}
\newcommand{\rsP}{\rs_{\P}}

\newcommand{\ash}{\as_{\hh}}

\newcommand{\msh}{\ms_{\hh}}
\newcommand{\zsh}{\zs_{\hh}}
\newcommand{\rsh}{\rs_{\hh}}

\newcommand{\ashP}{\as_{\hh,\P}}

\newcommand{\mshP}{\ms_{\hh,\P}}
\newcommand{\rshP}{\rs_{\hh,\P}}

\newcommand{\bil} [2]{\big<#1,#2\big>}

\newcommand{\SP} {S^{\P}}

\newcommand{\nlen}{\hspace{-0.2mm}}
\newcommand{\hskp}{\hspace{0.2mm}}

\newcommand{\snorm}  [2]{|#1|_{#2}}

\newcommand{\norm}   [2]{|\nlen|#1|\nlen|_{#2}}
\newcommand{\Norm}   [2]{\Big|\nlen\nlen\Big|\hskp #1\hskp\Big|\nlen\nlen\Big|_{#2}}
\newcommand{\NORM}   [2]{\left|\nlen\left|#1\right|\nlen\right|_{#2}}

\newcommand{\TNORM}  [2]{\left|\nlen\left|\nlen\left|#1\right|\nlen\right|\nlen\right|_{{}_{#2}}}

\newcommand{\abs}    [1]{|#1|}

\newcommand{\ABS}    [1]{\left|#1\right|}

\newcommand{\Vhk} {\VS{\hh}_{k}}

\newcommand{\Vvhk} {\textbf{V}^{\hh}_{k}}

\newcommand{\Pin}[1]{\Pi^{\nabla}_{#1}}
\newcommand{\Piz}[1]{\Pi^{0}_{#1}}

\newcommand{\restrict}[2]{{#1}_{|{#2}}}

\newcommand{\EOD}{\end{document}}

\newcommand{\beps}   {{\bm\varepsilon}}
\newcommand{\bsig}   {{\bm\sigma}}

\newcommand{\DIRI}{D}
\newcommand{\NEUM}{N}
\newcommand{\GamD}{\Gamma_{\DIRI}}
\newcommand{\GamN}{\Gamma_{\NEUM}}

\newcommand{\TRACE}{\textrm{tr}}

\newcommand{\gvN}{\gv_{\NEUM}}

\newcommand{\zero}{\mathbf{0}}

\newcommand{\dxv}{d\xv}

\newcommand{\ssP}{\ss_{\P}}

\newcommand{\VRTX}[1]{\calV_{#1}}

\usepackage[a4paper,total={6in,8in}]{geometry}

\begin{document}

\begin{frontmatter} 

  \title{The conforming virtual element method for polyharmonic and
    elastodynamics problems: a review}

  \author[MOX]  {P.~F.~Antonietti}
  \author[IMATI]{G.~Manzini}
  \author[MOX]  {I.~Mazzieri}
  \author[MOX]  {S.~Scacchi}
  \author[MOX]  {M.~Verani}
  
  \address[MOX]{MOX, Dipartimento di Matematica, Politecnico di
    Milano, Italy} \address[IMATI]{IMATI, Consiglio Nazionale delle
    Ricerche, Pavia, Italy }

  %% ----------------------------
  %% Abstract
  %% ----------------------------
  \begin{abstract}
    In this paper we review recent results on the conforming virtual
    element approximation of polyharmonic and elastodynamics
    problems.
    The structure and the content of this review is motivated by three
    paradigmatic examples of applications: classical and anisotropic
    Cahn-Hilliard equation and phase field models for brittle
    fracture, that are briefly discussed in the first part of the
    paper.
    We present and discuss the mathematical details of the conforming
    virtual element approximation of linear polyharmonic problems, the
    classical Cahn-Hilliard equation and linear elastodynamics
    problems.
  \end{abstract}
  
  %\begin{keyword}
  %\end{keyword}
  
\end{frontmatter}

%\tableofcontents

%% section-1
%% \input{Sec1_Introduction.tex}
\section{Introduction}
\label{sec:introduction}

In the recent years, there has been a tremendous interest to numerical
methods that approximate partial differential equations (PDEs) on
computational meshes with arbitrarily-shaped polytopal elements.
One of the most successful method is the virtual element method (VEM),
originally proposed
in~\cite{BeiraodaVeiga-Brezzi-Cangiani-Manzini-Marini-Russo:2013} for
second-order elliptic problems and then extended to a wide range of
applications.
The VEM was originally developed as a variational reformulation of the
\emph{nodal} mimetic finite difference (MFD) method~\cite{%
%BeiraodaVeiga-Manzini-Putti:2015,%
BeiraodaVeiga-Lipnikov-Manzini:2011,% 
Brezzi-Buffa-Lipnikov:2009,%
Manzini-Lipnikov-Moulton-Shashkov:2017% 
} for solving diffusion
problems on unstructured polygonal meshes.
A survey on the MFD method can be found in the review
paper~\cite{Lipnikov-Manzini-Shashkov:2014} and the research
monograph~\cite{BeiraodaVeiga-Lipnikov-Manzini:2014}.
The VEM inherits the flexibility of the MFD method with respect to the
admissible meshes and this feature is well reflected in the many
significant applications using polytopal meshes that have been
developed so far, see, for
example,~\cite{%
Dassi-Mascotto:2018,%
Antonietti-BeiraodaVeiga-Scacchi-Verani:2016,%
Cangiani-Georgoulis-Pryer-Sutton:2016,%
%BeiraodaVeiga-Lovadina-Vacca:2016,%
BeiraodaVeiga-Chernov-Mascotto-Russo:2016,%
%BeiraodaVeiga-Brezzi-Marini-Russo:2016a,%
BeiraodaVeiga-Brezzi-Marini-Russo:2016b,%
%BeiraodaVeiga-Brezzi-Marini-Russo:2016c,%
%BeiraodaVeiga-Brezzi-Marini-Russo:2016d,%
%Benedetto-Berrone-Borio-Pieraccini-Scialo:2016b,%
%Berrone-Borio-Scialo:2016,%
%Benedetto-Berrone-Scialo:2016,%
%Berrone-Pieraccini-Scialo:2016,%
Perugia-Pietra-Russo:2016,%
Wriggers-Rust-Reddy:2016,
%BeiraodaVeiga-Lovadina-Mora:2015,%
BeiraodaVeiga-Manzini:2015,
Mora-Rivera-Rodriguez:2015,%
%Natarajan-Bordas-Ooi:2015,
Berrone-Pieraccini-Scialo-Vicini:2015,%
%Vacca-BeiraodaVeiga:2015,
Paulino-Gain:2015,%
%Antonietti-BeiraodaVeiga-Mora-Verani:2014,%
%BeiraodaVeiga-Brezzi-Marini-Russo:2014,%
%BeiraodaVeiga-Brezzi-Marini-Russo:2014b,
BeiraodaVeiga-Manzini:2014,%
%Benedetto-Berrone-Pieraccini-Scialo:2014,%
%BeiraodaVeiga-Brezzi-Marini:2013,
%Brezzi-Marini:2013,%
%Berrone-Benedetto-Borio:2016chapter,%
%Berrone-Borio:2017,%
%Benedetto-Berrone-Borio-Pieraccini-Scialo:2016eccomas,%
Benvenuti-Chiozzi-Manzini-Sukumar:2019:CMAME:journal,%
%Antonietti-Manzini-Verani:2019:CAMWA:journal,%%
Certik-Gardini-Manzini-Mascotto-Vacca:2020,%
Certik-Gardini-Manzini-Vacca:2018:ApplMath:journal%%
%BeiraodaVeiga:2019:NumerMath:journal%
}.
Meanwhile, the mixed VEM for elliptic problems were introduced
in setting \emph{a la`} Raviart-Thomas
in~\cite{BeiraodaVeiga-Brezzi-Marini-Russo:2016c}
and in a BDM-like setting
in~\cite{Brezzi-Falk-Marini:2014}. 
The nonconforming formulation for diffusion problems was proposed
in~\cite{AyusodeDios-Lipnikov-Manzini:2016} as the finite element
reformulation of~\cite{Lipnikov-Manzini:2014} and later extended to
general elliptic
problems~\cite{Cangiani-Manzini-Sutton:2017,Berrone-Borio-Manzini:2018:CMAME:journal},
Stokes problem~\cite{Cangiani-Gyrya-Manzini:2016}, eigenvalue
problems~\cite{Gardini-Manzini-Vacca:2019:M2AN:journal}, and the
biharmonic
equation~\cite{Antonietti-Manzini-Verani:2018,Zhao-Chen-Zhang:2016}.
equation~\cite{Antonietti-Manzini-Verani:2018}.
Moreover, the connection between the VEM and the finite elements on
polygonal/polyhedral meshes is thoroughly investigated
in~\cite{Manzini-Russo-Sukumar:2014,Cangiani-Manzini-Russo-Sukumar:2015},
between VEM and discontinuous skeletal gradient discretizations
in~\cite{DiPietro-Droniou-Manzini:2018},
and between the VEM and the BEM-based FEM method
in~\cite{Cangiani-Gyrya-Manzini-Sutton:2017:GBC:chbook}.
The VEM was originally formulated
in~\cite{BeiraodaVeiga-Brezzi-Cangiani-Manzini-Marini-Russo:2013} as a
conforming FEM for the Poisson problem.
Then, it was later extended to convection-reaction-diffusion problems
with variable coefficients
in~\cite{BeiraodaVeiga-Brezzi-Marini-Russo:2016b}.
%%
%The connection with De Rham diagrams and Nedelec elements with
%application to electromagnetism has been explored
%in~\cite{BeiraodaVeiga-Brezzi-Marini-Russo:2016a}.

The virtual element method combines a great flexibility in using
polytopal meshes with a great versatility and easiness in designing
approximation spaces with high-order continuity properties on general
polytopal meshes.
These two features turn out to be essential in the numerical treatment
of the classical plate bending problem, for which a $\CS{1}$-regular
conforming virtual element approximation has been introduced
in~\cite{Brezzi-Marini:2013,Chinosi-Marini:2016}.
Virtual elements with $C^1$- regularity have been proposed to solve
elliptic problems on polygonal
meshes~\cite{BeiraodaVeiga-Manzini:2014} and polyedral meshes
in~\cite{BeiraodaVeiga-Dassi-Russo:2020}, the transmission eigenvalue
problem in~\cite{Mora-Velasquez:2018}, the vibration problem of
Kirchhoff plates in~\cite{Mora-Rivera-Velasquez:2018}, the buckling
problem of {K}irchhoff-{L}ove plates in~\cite{Mora-Velasquez:2020}.
The use of $C^1$-virtual elements has also been employed in the
conforming approximation of the Cahn-Hilliard problem
~\cite{Antonietti-BeiraodaVeiga-Scacchi-Verani:2016} and the von
K\'arm\'an equations~\cite{Lovadina-Mora-Velasquez:2019}, and in the
context of residual based a posteriori error estimators for
second-order elliptic problems~\cite{BeiraodaVeiga-Manzini:2015}.

Higher-order of regularity of the numerical approximation is also
required when addressing PDEs with differential operators of order
higher than two as the already mentioned biharmonic problem and the
more general case of the polyharmonic equations.
An example { of the latter} is found in the work of
Reference~\cite{Antonietti-Manzini-Verani:2019}.

In this paper we consider three paradigmatic examples of
applications where the conforming discretization requires highly
regular approximation spaces.
The first two examples are the the classical and the anisotropic
Cahn-Hilliard equations, that are used in modeling a wide range of
problems such as the tumor growth, the origin of the Saturn rings, the
separation of di-block copolymers, population dynamics, crystal growth, image
processing and even the clustering of mussels, see
\cite{Antonietti-BeiraodaVeiga-Scacchi-Verani:2016} and the
references therein.
The third example highlights the importance of coupling phase field
equations with the elastodynamic equation in the context of modeling
fracture propagation (see also
\cite{Aldakheel-Hudobivnik-Hussein-Wriggers:2018} for a phase-field
based VEM and the references therein).
These three examples motivate the structure of this review, where we
consider the conforming virtual approximation of the polyharmonic
equation, the classical Cahn-Hilliard equation and the time-dependent
elastodynamics equation.

%% POLY-HARMONIC
Historically, the numerical approximation of polyharmonic problems
dates back to the eighties~\cite{Bramble-Falk:1985}, and more
recently, this problem has been addressed in the context of the finite
element method
by~\cite{Barrett-Langdon-Nurnberg:2004,Gudi-Neilan:2011,Wang-Xu:2013,Schedensack:2016,Gallistl:2017}.
The conforming virtual element approximation of the biharmonic problem
has been addressed in \cite{Brezzi-Marini:2013,Chinosi-Marini:2016}.
while a non-conforming approximation has been proposed
in~\cite{Zhao-Chen-Zhang:2016,Antonietti-Manzini-Verani:2018,Zhao-Zhang-Chen-Mao:2018}.
In Section~\ref{sec:polyharmonic}, we review the conforming virtual
element approximation of polyharmonic problems
following~\cite{Antonietti-Manzini-Verani:2019}.
A nonconforming approximation is studied in~\cite{Chen-Huang:2020}.

%% CAHN-HILLIARD
%%
The Cahn-Hilliard equation involves fourth-order spatial derivatives
and the conforming finite element method is not really popular
approach because primal variational formulations of fourth-order
operators { requires the use of} finite element basis
functions that are piecewise-smooth and globally $\CS{1}$- continuous.
Only a few finite element formulations exists with the
$\CS{1}$-continuity property, see for
example~\cite{Elliott-Zheng:1986,Elliott-French:1987}, but in general,
these methods are not simple and easy to implement.
This high-regularity issue has successfully been addressed in the
framework of isogeometric analysis~\cite{Gomez-Calo-Bazilevs-Hughes:2008}.
The virtual element method provides a very effective framework for the
design and development of highly regular conforming approximation, and
in Section~\ref{sec:Cahn-Hilliard} we review the method proposed
in~\cite{Antonietti-BeiraodaVeiga-Scacchi-Verani:2016}.

Alternative approaches are offered by nonconforming
methods~\cite{Elliott-French:1989} or discontinuous
methods~\cite{Wells-Kuhl-Garikipati:2006}), but these methods do not
provide $\CS{1}$-regular approximations.
Another common strategy employed \emph{in practice} to solve the
Cahn-Hilliard equation by finite elements resorts to mixed methods;
see, e.g., \cite{Elliott-French-Milner:1989, Elliott-Larsson:1992} and
\cite{Kay-Styles-Suli:2009} for the continuous and discontinuous
setting, respectively.
Recently, mixed based discretization schemes on polytopal meshes {have been addressed in} \cite{Chave-DiPietro-Marche-Pigeonneau:2016} in the
context of the Hybrid High Order Method, and in  \cite{Xin-Zhangxin:2019} in
the context of the mixed Virtual Element Method.
However, mixed finite element methods requires a bigger number of
degrees of freedom, which implies, as a drawback, an increased
computational cost.

%% ELASTODYNAMICS
%%
Very popular strategies for numerically solving the time-dependent
elastodynamics equations in the \emph{displacement formulation} are
based on spectral
elements~\cite{Komatitsch-Tromp:1999,Faccioli-Maggio-Quarteroni-Tagliani:1996},
discontinuous Galerkin and discontinuous Galerkin spectral
elements~\cite{Riviere-Wheeler:2003,Antonietti-AyusodeDios-Mazzieri-Quarteroni:2016,Antonietti-Mazzieri-Quarteroni-Rapetti:2012}.
High-order DG methods for elastic and elasto-acoustic wave propagation
problems have been extended to arbitrarily-shaped polygonal/polyhedral
grids \cite{Antonietti-Mazzieri:2018,Antonietti-Bonaldi-Mazzieri:2018}
to further enhance the geometrical flexibility of the discontinuous
Galerkin approach while guaranteeing low dissipation and dispersion
errors.
Recently, the lowest-order Virtual Element Method has been applied for
the solution of the elastodynamics equation on nonconvex polygonal
meshes \cite{Park-Chi-Paulino:2019-CMAME,Park-Chi-Paulino:2019-IJNME}.
See also \cite{BeiraodaVeiga-Brezzi-Marini:2013} for the approximation
of the linear elastic problem, \cite{BeiraodaVeiga-Lovadina-Mora:2015}
for elastic and inelastic problems on polytopal meshes,
\cite{Vacca:2016} for virtual element approximation of hyperbolic
problems.
In Section~\ref{sec:elastodynamics}, we review
the conforming virtual element method of arbitrary order of accuracy
proposed
in~\cite{Antonietti-Manzini-Mazzieri-Mourad-Verani:2020}.

\subsection{Notation and technicalities}  
Throughout the paper, we consider the usual multi-index notation.
In particular, if $v$ is a sufficiently regular bivariate function and
$\alpha=(\alpha_1,\alpha_2)$ a multi-index with $\alpha_1$, $\alpha_2$ nonnegative
integer numbers, the function
$D^{\alpha}v=\partial^{|\alpha|}v\slash{\partial x_1^{\alpha_1}\partial
  x_2^{\alpha_2}}$ is the partial derivative of $v$ of order
$|\alpha|=\alpha_1+\alpha_2>0$.
For $\alpha=(0,0)$, we adopt the convention that $D^{\alpha}v$ coincides
with $v$.
Also, for the sake of exposition, we may use the shortcut notation
$\partial_{x}v$, $\partial_{y}v$, $\partial_{xx}v$, $\partial_{xy}v$,
$\partial_{yy}v$, to denote the first- and second-order partial
derivatives along the coordinate directions $x$ and $y$;
$\partial_{n}v$, $\partial_{t}v$, $\partial_{nn}v$, $\partial_{nt}v$,
$\partial_{tt}v$ to denote the first- and second-order normal and
tangential derivatives along a given mesh edge;
and $\partial^{m}_{n}v$ and $\partial^{m}_{t}v$ to denote the normal
and tangential derivative of $v$ of order $m$ along a given mesh edge.
Finally, let $\mathbf{n}=(\nx,\ny)$ and $\btau=(\tx,\ty)$ be the unit
normal and tangential vectors to a given edge $e$ of an arbitrary
polygon $\P$, respectively.
We recall the following relations between the first derivatives of
$v$:
\begin{align}
  \partial_nv \nabla v^T  \textbf{n} =  \nx\partial_xv + \ny\partial_yv,\quad
  \partial_{\tau}v = \nabla v^T \btau = \tx\partial_xv + \ty\partial_yv,
  \label{eq:edge:derivatives:first}
\end{align}
and the second derivatives of $v$:
\begin{align}
  \partial_{nn}v       = \textbf{n}^T\mathsf{H}(v)\textbf{n},\quad
  \partial_{n\tau}v    = \textbf{n}^T\mathsf{H}(v)\btau,\quad
  \partial_{\tau\tau}v = \btau^T     \mathsf{H}(v)\btau,
  \label{eq:edge:derivatives:second}
\end{align}
respectively, where the matrix $\mathsf{H}(v)$ is the Hessian of $v$, i.e.,
$\mathsf{H}_{11}(v)=\partial_{xx}v$,
$\mathsf{H}_{12}(v)=\mathsf{H}_{21}(v)=\partial_{xy}v$,
$\mathsf{H}_{22}(v)=\partial_{yy}v$.

%%%%%%%%

We use the standard definitions and notation of Sobolev spaces, norms
and seminorms~\cite{Adams-Fournier:2003}.
Let $k$ be a nonnegative integer number.
The Sobolev space $\HS{k}(\omega)$ consists of all square integrable
functions with all square integrable weak derivatives up to order $k$
that are defined on the open bounded connected subset $\omega$ of
$\REAL^{2}$.
As usual, if $k=0$, we prefer the notation $\LTWO(\omega)$.
Norm and seminorm in $\HS{k}(\omega)$ are denoted by
$\norm{\cdot}{k,\omega}$ and $\snorm{\cdot}{k,\omega}$, respectively,
and $(\cdot,\cdot)_{\omega}$ denote the $\LTWO$-inner product.
We omit the subscript $\omega$ when $\omega$ is the whole
computational domain $\Omega$. 
%%
%% The bold symbol will be employed to indicate the corresponding vector
%% valued functional space, e.g., $\mathbf{L}^2(\Omega)=[L^2(\Omega)]^2$.
%%
% If $\Gamma$ is the boundary of $\Omega$, we denote the space of the
% traces of the functions in $\HONE(\Omega)$, which are compactly
% supported in $\LTWO(\Gamma)$, by the notation
% $H^{\frac{1}{2}}_{c}(\Gamma)$, and a similar notation holds for the
% traces of the element boundary.

Given the mesh partitioning $\Th=\{\P\}$ of the domain $\Omega$ into
elements $\P$, we define the broken (scalar) Sobolev space for any
integer $k>0$
\begin{align*}
  \HS{k}(\Th) 
  = \prod_{\P\in\Th}\HS{k}(\P) 
  = \big\{\,\vs\in\LTWO(\Omega)\,:\,\restrict{\vs}{\P}\in\HS{k}(\P)\,\big\}, 
\end{align*}
which we endow with the broken $\HS{k}$-norm
\begin{align}
  \label{eq:Hs:norm-broken}
  \norm{\vs}{k,\hh}^2 = \sum_{\P\in\Th}\norm{\vs}{k,\P}^{2}
  \qquad\forall\,\vs\in\HS{k}(\Th),
\end{align}
and, for $k=1$, with the broken $\HONE$-seminorm
\begin{align}
  \label{eq:norm-broken}
  \snorm{\vs}{1,h}^2 = \sum_{\P\in\Th}\norm{\nabla\vs}{0,\P}^{2} 
  \qquad\forall\,\vs\in\HONE(\Th).
\end{align}

We denote the linear space of polynomials of degree up to $\ell$
defined on $\omega$ by $\PS{\ell}(\omega)$, with the useful
conventional notation that $\PS{-1}(\omega)=\{0\}$.
We denote the space of two-dimensional vector polynomials of degree up
to $\ell$ on $\omega$ by $\big[\PS{\ell}(\omega)\big]^2$; the space of
symmetric $2\times2$-sized tensor polynomials of degree up to $\ell$
on $\omega$ by $\PS{\ell,\textrm{sym}}^{2\times2}(\omega)$.
Space $\PS{\ell}(\omega)$ is the span of the finite set of
\emph{scaled monomials of degree up to $\ell$}, that are given by
\begin{align*}
  \calM_{\ell}(\omega) =
  \bigg\{\,
    \left( \frac{\xv-\xv_{\omega}}{\hh_{\omega}} \right)^{\alpha}
    \textrm{~with~}\abs{\alpha}\leq\ell
    \,\bigg\},
\end{align*}
where 
\begin{itemize}
\item $\xv_{\omega}$ denotes the center of gravity of $\omega$ and
  $\hh_{\omega}$ its characteristic length, as, for instance, the edge
  length or the cell diameter for $\DIM=1,2$;
\item $\alpha=(\alpha_1,\alpha_2)$ is the two-dimensional multi-index
  of nonnegative integers $\alpha_i$ with degree
  $\abs{\alpha}=\alpha_1+\alpha_{2}\leq\ell$ and such that
  $\xv^{\alpha}=\xs_1^{\alpha_1}\xs_{2}^{\alpha_{2}}$ for any
  $\xv\in\REAL^{2}$.
\end{itemize}
We will also use the set of \emph{scaled monomials of degree exactly
  equal to $\ell$}, denoted by $\calM_{\ell}^{*}(\omega)$ and obtained
by setting $\abs{\alpha}=\ell$ in the definition above.

\medskip
Finally, we use the letter $C$ in the estimates below to denote a
strictly positive constant whose value can change at any instance but that is independent of the discretization parameters such as the mesh
size $\hh$.
Note that $C$ may depend on the { the polynomial order}, on the constants of the model equations or
the variational problem, like the coercivity and continuity constants,
or even constants that are uniformly defined for the family of meshes
of the approximation while $\hh\to0$, such as the mesh regularity
constant, the stability constants of the discrete bilinear forms, etc.
Whenever it is convenient, we will simplify the notation by using
expressions like $x\lesssim y$ and $x\gtrsim y$ to mean that $x\leq
Cy$ and $x\geq Cy$, respectively, $C$ being the generic constant in
the sense defined above.

\subsection{Mesh assumptions}
\label{subsec:mesh:assumptions}

Throughout the paper we assume that $\calT=\big\{\Th\big\}_{\hh}$ is a
family of decompositions of the computational domain $\Omega$, where
each mesh $\Th$ is a collection of nonoverlapping polygonal elements
$\P$ with boundary $\partial\P$,  such that $\bar\Omega= \Cup_{\P\in\Th} \bar{\P}$.
Each mesh is labeled by the \emph{mesh size} $\hh$, the diameter of
the mesh, defined as usual by $\hh=\max_{\P\in\Th}\hP$, where
$\hP=\sup_{\xv,\yv\in\P}\ABS{\xv-\yv}$.
We assume the mesh sizes of family $\calT$ form a countable subset of $\calH=(0,\infty)$ having zero as its unique accumulation point.
We denote the set of mesh vertices $\V$ by $\calV_{\hh}$ and the set
of mesh edges $\E$ by $\calE_{\hh}$
Moreover, the symbol $\hV$ is a characteristic length associated with
each vertex; more precisely, $\hV$ is the average of the diameters of
the polygons sharing vertex $\V$.
We consider the following mesh regularity assumptions:
\begin{description}
\item\TERM{M} There exists a positive constant $\gamma$, \emph{mesh
regularity constant}, which is independent of $h$ (and $\P$) and such
  that for $K\in\Omega_h$ there hold:
  \begin{itemize}
  \item\TERM{M1} $\P$ is star-shaped with respect to every point of a
    ball of radius $\gamma\hP$, where $\hP$ is the diameter of $\P$;
  \item\TERM{M2} for every edge $\E$ of the cell boundary $\partial\P$
    of every cell $\P$ of $\Th$, it holds that $\hE\geq\gamma\hP$,
    where $\hE$ denotes the length of $\E$.
  \end{itemize}
\end{description}
All the results contained in the rest of the paper are obtained under
assumptions \TERM{M1}-\TERM{M2}.

%% section-2: The virtual element method for the polyharmonic problem
%%\label{sec:polyharmonic}
%% \input{Sec2_VEM_Polyharmonic.tex}
\section{The virtual element method for the polyharmonic problem}
\label{sec:polyharmonic}

\subsection{The continuous  problem}
\label{subsec:polyharmonic:continuous:problem}

Let $\Omega\subset\REAL^2$ be a open, bounded, convex domain with
polygonal boundary $\Gamma$.
For any integer $p\geq 1$, we introduce the conforming virtual element
method for the approximation of the following problem:
\begin{subequations}\label{eq:poly:pblm:continuous}
  \begin{align}
    (-\Delta)^p  \us & = \fs\phantom{0} \qquad\text{in~}\Omega,\label{eq:poly:pblm:1}\\
    \partial^j_n \us & = 0\phantom{\fs} \qquad\text{for~}j=0,\ldots,p-1\text{~on~}\Gamma,\label{eq:poly:pblm:2}
  \end{align}
\end{subequations}
(recall the conventional notation $\partial^0_n\us=\us$).
Let
\begin{align*}
  \Vs\equiv\HS{p}_{0}(\Omega) = \big\{
  \vs\in\HS{p}(\Omega):\partial^j_n\vs=0\text{~on~}\Gamma,\,j=0,\ldots,p-1
  \big\}.
\end{align*}
Denoting the duality pairing between $\Vs$ and its dual $\Vsp$ by
$\bil{\cdot}{\cdot}$, the variational formulation of the polyharmonic
problem \eqref{eq:poly:pblm:continuous} reads as: \emph{Find
$\us\in\Vs$ such that}
\begin{equation}\label{eq:poly:pblm:wp}
  \as(\us,\vs) = \bil{\fs}{\vs} \qquad\forall\vs\in\Vs,
\end{equation}
where, for any nonnegative integer $\ell$, the bilinear form is given
by:
\begin{align}
  \as(\us,\vs) = 
  \begin{cases}
    \,\int_{\Omega} \nabla\Delta^\ell\us\cdot\nabla\Delta^\ell\vs\,\dx  & \mbox{for~$p=2\ell+1$},\\[1em]
    \,\int_{\Omega} \Delta^\ell\us\,\Delta^\ell\vs\,\dx                 & \mbox{for~$p=2\ell$}.
  \end{cases}
\end{align}
Whenever $\fs\in\LTWO(\Omega)$ we have
\begin{align}
  \bil{\fs}{\vs} = (\fs,\vs) = \int_{\Omega}\fs\vs\dV\,\dx.
  \label{eq:poly:pblm:p-rhs}
\end{align}
where $(\cdot,\cdot)$ denotes the $\LTWO$-inner product.
The existence and uniqueness of the solution to
\eqref{eq:poly:pblm:wp} follows from the Lax-Milgram Theorem because
of the continuity and coercivity of the bilinear form $\as(\cdot,\cdot)$
with respect to $\|\cdot\|_V=\vert \cdot \vert_{p,\Omega}$ which is a
norm on $\HS{p}_{0}(\Omega)$.
Moreover, since $\Omega$ is a convex polygon,
from~\cite{Gazzola-Grunau-Sweers:1991} we know that
$\us\in\HS{2p-m}(\Omega)\cap\HS{p}_{0}(\Omega)$ if
$\fs\in\HS{-m}(\Omega)$, $m\leq p$ and it holds that
$\norm{\us}{2p-m}\leq\Cs\norm{\fs}{-m}$.
In the following, we denote the coercivity and continuity constants of
$\as(\cdot,\cdot)$ by $\alpha$ and $\Ms$, respectively.

Let $\P$ be a polygonal element and set
\begin{align*}
  \asP(\us,\vs) = 
  \begin{cases}
    \,\int_{\P} \nabla\Delta^\ell\us\cdot\nabla\Delta^\ell\vs\,\dx & \mbox{for~$p=2\ell+1$},\\[1em]
    \,\int_{\P} \Delta^\ell\us\,\Delta^\ell\vs\,\dx                & \mbox{for~$p=2\ell$}.
  \end{cases}
\end{align*}
For an odd $p$, i.e., $p=2\ell+1$, a repeated application of the
integration by parts formula yields
\begin{align} 
  \asP(\us,\vs) 
  =
  & -\int_{\P} \Delta^p\us\,\vs\,\dx + \int_{\partial\P}\partial_n(\Delta^\ell\us)\,\Delta^\ell\vs\dS\nonumber \\[0.5em]
  &  +\sum_{i=1}^\ell
  \left( 
    \int_{\partial\P}\partial_n(\Delta^{p-i}\us)\,\Delta^{i-1}\vs\dS
    -\int_{\partial\P}\Delta^{p-i}\us\,\partial_n(\Delta^{i-1}\vs)\dS
  \right),
  \label{eq:poly:intbyparts:odd:p}
\end{align}
while, for an even $p$, i.e., $p=2\ell$, we have
\begin{align} 
  \asP(\us,\vs) 
  &= \int_{\P}\Delta^p\us\,\vs\,\dx
  \nonumber\\[0.5em]
  &\phantom{=}
  + \sum_{i=1}^\ell
  \left(
  \int_{\partial\P}  \partial_n(\Delta^{p-i}\us)\,\Delta^{i-1}\vs\,\dS
  -\int_{\partial\P} \Delta^{p-i}\us\,\partial_n(\Delta^{i-1}\vs)\,\dS
  \right).
  \label{eq:poly:intbyparts:even:p}
\end{align} 

\subsection{The conforming virtual element approximation}
\label{sec:VEM}

The conforming virtual element discretization of problem
\eqref{eq:poly:pblm:wp} hinges upon three mathematical objects: (1)
the finite dimensional conforming virtual element space
$\Vhrp{r}\subset V$; (2) the continuous and coercive discrete bilinear
form $\ash(\cdot,\cdot)$; (3) the linear functional
$\bil{\fsh}{\cdot}$.

Using such objects, we formulate the virtual element method as:
\emph{Find $\ush\in\Vhrp{r}$ such that}
\begin{align}
  \ash(\ush,\vsh) = \bil{\fsh}{\vsh} \quad\forall\vsh\in\Vhrp{r}.
  \label{eq:poly:VEM}
\end{align}
The existence and uniqueness of the solution $\ush$ is again a consequence
of the Lax-Milgram theorem.\cite[Theorem~2.7.7, page~62]{Brenner-Scott:2008}.

\subsubsection{Virtual element spaces}\label{subsec:VE_space}
For $p\geq 1$ and $r\geq 2p-1$, the local Virtual Element space on
element $\P$ is defined by
%%
%% \begin{align}
%%   \VhPrp{r} = \Big
%%   \{\vsh\in\HS{p}(\P):\,
%%   & \Delta^p\vsh\in\PS{r-2p}(\P),\,
%%   \vsh\in\PS{r}(\E),\,\partial^i_n\vsh\in\PS{r-i}(\E),\,i=1,\ldots,p-1~\forall
%%   \E\in\partial\P
%%   \Big\},
%% \end{align}
%%
\begin{multline*}
  \VhPrp{r} = \Big
  \{\vsh\in\HS{p}(\P):\,
  \Delta^p\vsh\in\PS{r-2p}(\P),\,
  \vsh\in\PS{r}(\E),\,\partial^i_n\vsh\in\PS{r-i}(\E),\,\\
  i=1,\ldots,p-1~\forall
  \E\in\partial\P
  \Big\},
\end{multline*}
with the conventional notation that $\PS{-1}(\P)=\{0\}$.
The virtual element space $\VhPrp{r}$ contains the space of
polynomials $\PS{r}(\P)$, for $r\geq 2p-1$.
Moreover, for $p=1$, it coincides with the conforming virtual element
space for the Poisson equation
\cite{BeiraodaVeiga-Brezzi-Cangiani-Manzini-Marini-Russo:2013}, and
for $p=2$, it coincides with the conforming virtual element space for
the biharmonic equation \cite{Brezzi-Marini:2013}.
The requirement $\vsh\in\HS{p}(\P)$ implies that suitable
compatibility conditions for $\vsh$ and its derivatives up to order
$p-1$ must hold at the vertices of the polygon (see, e.g.,
\cite[Theorems~1.5.2.4~and~1.5.7.8]{Grisvard:1985} and
\cite[Section~5]{Bernardi-Dauge-Maday:2007}).

\begin{figure}
  \begin{center}
    \begin{tabular}{cc}
      \includegraphics[scale=0.5]{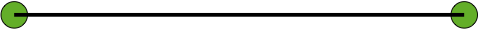} & \hspace{2cm} %% {./figures-PNG/dofs_a0_m1} & \hspace{2cm}
      \includegraphics[scale=0.5]{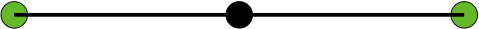} \\             %% {./figures-PNG/dofs_a0_m2} \\
      $p=1,\,r=1$ & \hspace{2cm} $p=1,\,r=2$ \\[0.75em]
      %---------------------------------------------------------------------------------------------------------------------------------
      \includegraphics[scale=0.5]{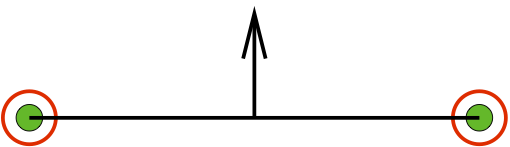} & \hspace{2cm} %% {./figures-PNG/dofs_a1_m3} & \hspace{2cm}
      \includegraphics[scale=0.5]{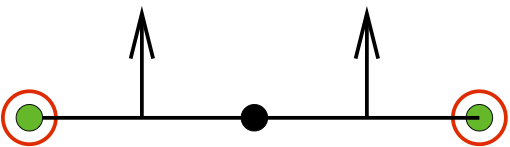} \\             %% {./figures-PNG/dofs_a1_m4} \\
      $p=2,\,r=3$ & \hspace{2cm} $p=2,\,r=4$ \\[0.75em]
      %---------------------------------------------------------------------------------------------------------------------------------
      \includegraphics[scale=0.5]{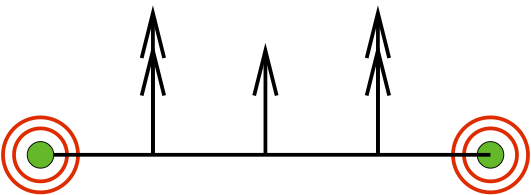} & \hspace{2cm} %% {./figures-PNG/dofs_a2_m5} & \hspace{2cm}
      \includegraphics[scale=0.5]{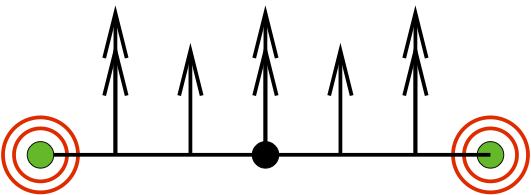} \\             %% {./figures-PNG/dofs_a2_m6} \\
      $p=3,\,r=5$ & \hspace{2cm} $p=3,\,r=6$
    \end{tabular}
    \caption{[Taken from \cite{Antonietti-Manzini-Verani:2019}]. Edge
      degrees of freedom of the virtual element space
      $V^p_{h,r}(\P)$ for the polyharmonic problem with $p=1$ (top
      panels, Laplace operator), $p=2$ (middle panels, bi-harmonic
      operator), $p=3$ (bottom panels, triharmonic operator).
      Here, $p$ is the order of the partial differential operator;
      $r=1,2,\ldots,6$ are the integer parameters that specify the
      degree of the polynomial subspace $\PS{r}(\P)$ of the VEM space
      $V_{h,r}^3(\P)$.
      The (green) dots at the vertices represent the vertex values and
      each (red) vertex circle represents an order of derivation.
      The (black) dot on the edge represents the moment of
      $\restrict{v_h}{e}$; the arrows represent the moments of
      $\restrict{\partial_nv_h}{e}$; the double arrows represent the
      moments of $\restrict{\partial_{nn}v_h}{e}$.
      The corresponding internal degrees of freedom (D4) are absent in
      the case $r=2p-1$, while reduce to a single one in the case
      $r=2p$.
    }
    \label{fig:trih:dofs}
  \end{center}
\end{figure}

We characterize the functions in $\VhPrp{r}$ through the following set of
\emph{degrees of freedom}:
\medskip
\begin{description}
\item\TERM{D1} $\hV^{|\nu|}\Ds^{\nu}\vsh(\V)$, $\ABS{\nu}\leq p-1$ for
  any vertex $\V$ of the polygonal boundary $\partial\P$;
\medskip
\item\TERM{D2} $\displaystyle\hE^{-1}\int_{\E}\qs\vsh\dS$ for any
  $\qs\in\PS{r-2p}(\E)$ and any edge $\E$ of the polygonal boundary
  $\partial\P$;
  \medskip
\item\TERM{D3}
  $\displaystyle\hE^{-1+j}\int_{\E}\qs\partial^j_n\vsh\dS$ for any
  $\qs\in\PS{r-2p+j}(\E)$, $j=1,\ldots,p-1$ and any edge $\E$ of
  $\partial\P$;
  \medskip
\item\TERM{D4} $\displaystyle\hP^{-2}\int_{\P}\qsh\vsh\,\dx$ for
  any $\qs\in\PS{r-2p}(\P)$.
\end{description}
Here, as usual, we assume that $\PS{-n}(\cdot)=\{0\}$ for $n\geq1$.
Figure~\ref{fig:trih:dofs} illustrates the degrees of freedom on a
given edge $\E$ for $p=1,2,3$ (Laplace, biharmonic, and triharmonic
case) and $r=2p-1,2p$; the corresponding internal degrees of freedom
(D4) are absent in the case $r=2p-1$, while reduce to a single one in
the case $r=2p$.
Finally, we note that in general the internal degrees of freedom
\TERM{D4} make it possible to define the $\LTWO$-orthogonal polynomial
projection of $\vsh$ onto the space of polynomial of degree $r-2p$.

\medskip
The dimension of $\VhPrp{r}$ is 
\begin{align*}
  \DIM(\VhPrp{r}) 
  &= \frac{p(p+1)}{2}\NP + \NP\sum_{j=0}^{p-1}( r-2p+j+1)    + \frac{(r-2p+1)(r-2p+2)}{2},
\end{align*}
where $\NP$ is the number of vertices, which equals the number of
edges, of $\P$.

\medskip 
In \cite{Antonietti-Manzini-Verani:2019}, it is proved that the above
choice of degrees of freedom is unisolvent in $\VhPrp{r}$.

%-----------Global space 
\medskip
Building upon the local spaces $\VhPrp{r}$ for all $\P\in\Th$, the
{\em global} conforming virtual element space $\Vhrp{r}$ is defined on
$\Omega$ as
\begin{align}
  \Vhrp{r} = \Big\{
  \vsh\in\HS{p}_{{0}}(\Omega)\,:\,\restrict{\vsh}{\P}\in\VhPrp{r}\,\,\forall\P\in\Th
  \Big\}.
  \label{eq:poly:global:space}
\end{align}
We remark that the associated global space is made of $\HS{p}(\Omega)$
functions.
Indeed, the restriction of a virtual element function $\vsh$ to each
element $\P$ belongs to $\HS{p}(\P)$ and glues with
$C^{p-1}$-regularity across the internal mesh faces.
The set of global degrees of freedom inherited by the local degrees of
freedom are:
\medskip
\begin{itemize}
\item $\hV^{|\nu|}\Ds^{\nu}\vsh(\V)$, $\ABS{\nu}\leq p-1$ for every
  interior vertex $\V$ of $\Th$;
  
  \medskip
\item $\displaystyle\hE^{-1}\int_{\E}\qs\vsh\dS$ for any
  $\qs\in\PS{r-2p}(\E)$ and every interior edge $\E\in\Eh$;

  \medskip
\item $\displaystyle\hE^{-1+j}\int_{\E}\qs\partial^j_n\vsh\,ds$ for
  any $\qs\in\PS{r-2p+j}(e)$ $j=1,\ldots,p-1$ and every interior edge
  $\E\in\Eh$;
\medskip
\item $\displaystyle\hP^{-2}\int_{\P}\qs\vsh\dx$ for any
  $\qs\in\PS{r-2p}(\P)$ and every $\P\in\Th$.
\end{itemize}%{description}

\subsubsection{Modified lowest order virtual element spaces}
\label{subsec:modified_VE_space}
In this section, we briefly discuss the possibility of introducing
modified lowest order virtual element spaces with a reduced number of
degrees of freedom with respect to the corresponding lowest order ones
that were introduced previously.
The price we pay is a reduced order of accuracy since the polynomial
functions included in such modified spaces has a lower degree.

For the sake of presentation we start from the case $p=3$, while we
refer the reader to \cite{Brezzi-Marini:2013} for the case of $p=2$ and
Section \ref{sec:vem1} where the reduced virtual space is employed in
the context of the approximation of the Cahn-Hilliard problem.
Consider the modified local virtual element space:
%%
%% \begin{align*}
%%   \tilde{V}_{h,5}^3(\P) = \Big\{ 
%%   \vsh\in\HS{3}(\P)
%%   \,:\,
%%   &\Delta^3\vsh = 0, 
%%   \vsh\in\PS{5}(\E),\,
%%   \partial_n\vsh\in\PS{3}(\E),\,
%%   \partial_{nn}\vsh\in\PS{2}(\E)~\forall\E\in\partial\P
%%   \Big\}
%% \end{align*}
%%
\begin{multline*}
  \tilde{V}_{h,5}^3(\P) = \Big\{ 
  \vsh\in\HS{3}(\P)
  \,:\,
  \Delta^3\vsh = 0, 
  \vsh\in\PS{5}(\E),\,
  \partial_n\vsh\in\PS{3}(\E),\,\\[0.5em]
  \partial_{nn}\vsh\in\PS{2}(\E)~\forall\E\in\partial\P
  \Big\}
\end{multline*}
with associated degrees of freedom: 
\begin{description}
\item\TERM{D1'} $\hV^{|\nu|}D^{\nu}\vsh(\V)$, $\ABS{\nu}\leq2$ for any
  vertex $\V$ of $\partial\P$;
  
  \medskip
\item\TERM{D2'} $\displaystyle\hE\int_{\E}\partial_{nn}\vsh\dS$ for
  any edge $\E$ of $\partial\P$.
\end{description}

In Ref.~\cite{Antonietti-Manzini-Verani:2019}, we proved that the
degrees of freedom \TERM{D1'} and \TERM{D2'} are unisolvent in
$\tilde{V}_{h,5}^3(\P)$ and this space contains the linear subspace of
polynomials of degree up to $4$.
Moreover, the associated global space obtained by gluing together all
the elemental spaces $ \tilde{V}_{h,5}^3(\P)$ reads as:
\begin{align}
  \tilde{V}_{h,5}^3 = \Big\{ 
  \vsh\in\HS{3}_0(\Omega)\,:\,\vsh\vert_{\P}\in\tilde{V}_{h,5}^3(\P)
  \,\,\forall\P\in\Th
  \Big\},
  \label{eq:trih:global:space:modified:VEM}
\end{align}
is made of $\HS{3}(\Omega)$ functions.
  
Analogously, in the general case we can build the modified lowest
order spaces containing the space of polynomials of degree up to
$2p-2$:
\begin{multline*}
  \tVhPrp{2p-1} = \Big\{\vsh\in\HS{p}(\P):\,
  \Delta^p\vsh = 0,
  \vsh\in\PS{2p-1}(\E),\,\partial^i_n\vsh\in\PS{2p-2-i}(\E),\\[0.5em]
  \,i=1,\ldots,p-1~\forall\E\in\partial\P
  \Big\},
\end{multline*}
with associated degrees of freedom:
\begin{description}
\item\TERM{D1'} $\hV^{|\nu|}D^{\nu}\vsh(\V)$, $\ABS{\nu}\leq p-1$ for
  any vertex $\V$ of $\partial\P$;
  
  \medskip
\item\TERM{D2'}
  $\displaystyle\hE^{-1+j}\int_{\E}\qs\partial_{n}^i\vsh\dS$ for any
  $\qs\in\PS{j-2}(\E)$ and edge $\E$ of $\partial\P$,
  $j=1,\ldots,p-1$.
\end{description}

\subsubsection{Discrete bilinear form}
\label{subseq:VEM:bilinear:form}
%% contruction of the elliptic projection

To define the elliptic projection $\PiPr{r}:\VhPrp{r}\to\PS{r}(\P)$,
we first need to introduce the \emph{vertex average projector}
$\widehat{\Pi}^{\P}:\VhPrp{r}\to\PS{0}(\P)$, which projects any smooth
enough function defined on $\P$ onto the space of constant
polynomials.
To this end, consider the continuous function $\psi$ defined on $\P$.
The \emph{vertex average projection} of $\psi$ onto the constant
polynomial space is given by:
\begin{align}
  \widehat{\Pi}^{\P}\psi = \frac{1}{\NP}\sum_{\vrtx\in\partial\P}\psi(\vrtx).
  \label{eq:trih:vertex:average:projection}
\end{align}

Finally, we define the elliptic projection
$\PiPr{r}:\VhPrp{r}\to\PS{r}(\P)$ as the solution of the following finite
dimensional variational problem
\begin{align}
  \asP(\PiPr{r}\vsh,\qs)                &= \asP(\vsh,\qs)\phantom{ \widehat{\Pi}^{\P}\Ds^{\nu}\vsh }\forall\qs\in\PS{r}(\P),\label{eq:poly:Pi:A}\\[0.5em]
  \widehat{\Pi}^{\P}\Ds^{\nu}\PiPr{r}\vsh &= \widehat{\Pi}^{\P}\Ds^{\nu}\vsh\phantom{\asP(\vsh,\qs)} \ABS{\nu}\leq{p-1}.      \label{eq:poly:Pi:B}
\end{align}
According to Reference~\cite{Antonietti-Manzini-Verani:2019}, such
operator has two important properties:
\begin{itemize}
\item[$(i)$] it is a polynomial-preserving operator in the sense that
  $\PiPr{r}\qs=\qs$ for every $\qs\in\PS{r}(\P)$;
\item[$(ii)$] $\PiPr{r}\vsh$ is \emph{computable} using only the
  degrees of freedom of $\vsh$.
\end{itemize}

%% VEM bilinear form
We write the symmetric bilinear form
$\ash:\Vhrp{r}\times\Vhrp{r}\to\REAL$ as the sum of local terms
\begin{align}
  \ash(\ush,\vsh) = \sum_{\P\in\Th}\ashP(\ush,\vsh),
\end{align}
where each local term $\ashP:\VhPrp{r}\times\VhPrp{r}\to\REAL$ is a
symmetric bilinear form.
We set
\begin{align}
  \ashP(\ush,\vsh) 
  = \asP(\PiPr{r}\ush,\PiPr{r}\vsh) 
  + \SP(\ush-\PiPr{r}\ush,\vsh-\PiPr{r}\vsh),
  \label{eq:poly:ah:def}
\end{align}
where $\SP:\VhPrp{r}\times\VhPrp{r}\to\mathbbm{R}$ is a symmetric
positive definite bilinear form such that
\begin{align}
  \sigma_*\asP(\vsh,\vsh)\leq\SP(\vsh,\vsh)\leq\sigma^*\asP(\vsh,\vsh)
  \qquad\forall\vsh\in\VhPrp{r}\textrm{~with~}\PiPr{r}\vsh=0,
  \label{eq:poly:S:stability:property}
\end{align}
for two some positive constants $\sigma_*$, $\sigma^*$ independent of
$\hh$ and $\P$.
The bilinear form $\ashP(\cdot,\cdot)$ has the two fundamental
properties of $r$-\emph{consistency} and
\emph{stability}~\cite{Antonietti-Manzini-Verani:2019}:
%:
\begin{description}
\item[$(i)$] $r$-\textbf{Consistency}: for every polynomial
  $\qs\in\PS{r}(\P)$ and function $\Vhrp{r}(\P)$ we have:
  \begin{align}
    \ashP(\vsh,\qs) = \asP(\vsh,\qs);
    \label{eq:poly:r-consistency}
  \end{align}
\item[$(ii)$] \textbf{Stability}: there exist two positive constants
  $\alpha_*$, $\alpha^*$ independent of $h$ and $\P$ such that for
  every $\vsh\in\VhPrp{r}$ it holds:
  \begin{align}
    \alpha_*\asP(\vsh,\vsh)\leq\ashP(\vsh,\vsh)\leq\alpha^*\asP(\vsh,\vsh).
    \label{eq:poly:stability}
  \end{align}
\end{description}

\subsubsection{Discrete load term}
\label{subseq:VEM:load:term}
We denote by $\fsh$ the piecewise polynomial approximation of $f$ on
$\Th$ given by
\begin{equation}\label{eq:rhs}
  \restrict{\fsh}{\P}=\Pizr{r-p}f,
\end{equation}
for $r\geq 2p-1$ and $\P\in\Th$. 
Then, we set
\begin{equation}\label{vem:rhs}
  \bil{\fsh}{\vsh} = \sum_{\P\in\Th}\int_{\P}\fsh\vsh\,\dx
\end{equation}
which implies, using the $\LTWO$-orthogonal projection, that
\begin{equation}\label{aux:1.1}
  \bil{\fsh}{\vsh}
  = \sum_{\P\in\Th}\int_{\P}\Pizr{r-p}\,\fs\Pizr{r-p}\vsh\,\dx
  = \sum_{\P\in\Th}\int_{\P}\fs\,\Pizr{r-p}\vsh\,\dx. 
\end{equation}
The right-hand side of \eqref{aux:1.1} is computable by a combined use
of the degrees of freedom \TERM{D1}-\TERM{D4} and the enhanced
approach of Reference~\cite{Ahmad-Alsaedi-Brezzi-Marini-Russo:2013}.

\subsubsection{VEM spaces with arbitrary degree of continuity}
\label{sec:polyharmonic-2}

In this section we briefly sketch the construction of global virtual
element spaces with arbitrary high order of continuity.
More precisely, we consider the local virtual element space defined as
before, for $r \geq 2p-1$:
\begin{multline*}
  \VhPrp{r} = \Big\{ \vsh\in\HS{p}(\P):\, 
  \Delta^p \vsh\in\PS{r-2p}(\P),\,
  \vsh\in\PS{r}(\E),\,\partial^j_n\vsh\in\PS{r-j}(\E),\,
  \\[0.5em]
  j=1,\ldots,p-1~\forall\E\in\partial\P\big\}.
\end{multline*}

Differently from the previous section, we make the degrees of freedom
depend on a given parameter $t$ with $0\leq t \leq p-1$.
For a given value of $t$ we choose the \emph{degrees of freedom} as
follows

\smallskip
\begin{description}
\item\TERM{D1} $\hV^{|\nu|}\Ds^{\nu}\vsh(\V)$, $\ABS{\nu}\leq p-1$ for
  any vertex $\V$ of $\P$;
  
  \medskip
\item\TERM{D2} $\displaystyle\hE^{-1}\int_{\E}\vsh\qs\dS$ for any
  $\qs\in\PS{r-2p}(\E)$,\, for any edge $\E$ of $\partial\P$;
  
  \medskip
\item\TERM{D3}
  $\displaystyle\hE^{-1+j}\int_{\E}\partial^j_n\vsh\qs\dS$ for any
  $\qs\in\PS{r-2p+j}(\E)$ and edge $\E\in\partial\P$,
  $j=1,\ldots,p-1$;
  
  \medskip
\item\TERM{D4'} $\displaystyle\hP^{-2}\int_{\P}\qs\vsh\,\dx$ for
  any $\qs\in\PS{r-2(p-t)}(\P)$;
\end{description}
where as usual we assume $\PS{-n}(\cdot)=\{0\}$ for $n=1,
2,3,\ldots$.

\medskip
This set of degrees is still unisolvent, cf.~\cite{Antonietti-Manzini-Verani:2019}.
Moreover, for $r\geq 2p-1$ it holds that
$\PS{r}(\P)\subset\Vs^{p}_{h,r}(\P)$.
Finally, it is worth noting that the choice \TERM{D4'}, if compared
with \TERM{D4}, still guarantees that the associated global space is
made of $\Cs{p-1}$ functions.

However, in this latter case we can use the degrees of freedom
\TERM{D1}-\TERM{D4'} to solve a differential problem involving the
$\Delta^{p-t}$ operator and $\CS{p-1}(\Omega)$ basis functions.
For the sake of exposition, let us consider the following two
examples, in the context of the Laplacian and the Bilaplacian
problem.

\medskip
\begin{enumerate}
\item Choosing $p$ and $t$ such that $p-t=1$ we obtain a
  $\CS{p-1}$-conforming virtual element method for the solution of
  the Laplacian problem.
  For example, for $p=3$, $t=2$ and $r=5$, the local space
  $V^3_{h,5}(\P)$ endowed with the corresponding degrees of freedom
  \TERM{D1}-\TERM{D4'} can be employed to build a global space made of
  $\CS{2}$ functions.
  It is also worth mentioning that the new choice \TERM{D4'},
  differently from the original choice \TERM{D4}, is essential for the
  computability of the elliptic projection, see
  \eqref{eq:poly:Pi:A}-\eqref{eq:poly:Pi:B}, with respect to the
  bilinear form
  $\asP(\cdot,\cdot)=\int_\P\nabla(\cdot)\nabla(\cdot)\,\dx$.
  
  \medskip
\item Choosing $p$ and $t$ such that $p-t=2$ we have a
  $\CS{p-1}$-conforming virtual element method for the solution of the
  Bilaplacian problem.
  For example, for $p=3$, $t=1$ and $r=5$, similarly to the previous
  case, the space $V^3_{h,5}(\P)$ together with \TERM{D1}-\TERM{D4'}
  provides a global space of $\CS{2}$ functions that can be employed
  for the solution of the biharmonic problem.
\end{enumerate}

It is worth remembering that $\CS{1}$-regular virtual element basis
function has been employed, e.g., in \cite{BeiraodaVeiga-Manzini:2015}
to study residual based a posteriori error estimators for the virtual
element approximation of second order elliptic problems.
Moreover, the solution of coupled elliptic problems of different order
can take advantage from this flexibility of the degree of continuity
of the basis functions.
Indeed, for the sake of clarity consider the conforming virtual
element approximation of the following simplified situation:

\begin{align*}
  -\Delta\us_1    &= \fs_1 \phantom{\fs_2\us_1\partial_n\us_20} \hspace{-5mm}\text{in~}\Omega_1,\\
  \Delta^2\us_2   &= \fs_2 \phantom{\fs_1\us_1\partial_n\us_20} \hspace{-5mm}\text{in~}\Omega_2,\\
  \us_1           &= \us_2 \phantom{\fs_1\fs_2\partial_n\us_20} \hspace{-5mm}\text{on~}\Gamma=\Omega_1\cap\Omega_2,\\
  \partial_n\us_1 &= \partial_n\us_2 \phantom{\fs_1\fs_1\us_10} \hspace{-5mm}\text{on~}\Gamma,\\
  \us_1           &= 0 \phantom{\fs_1\fs_2\us_2\partial_n\us_2} \hspace{-5mm}\text{on~}\partial\Omega_1\setminus\Gamma,\\
  \us_2           &= 0 \phantom{\fs_1\fs_2\us_2\partial_n\us_2} \hspace{-5mm}\text{on~}\partial\Omega_2\setminus\Gamma,\\
  \partial_n\us_2 &= 0 \phantom{\fs_1\fs_2\us_2\partial_n\us_2} \hspace{-5mm}\text{on~}\partial\Omega_2\setminus\Gamma.
\end{align*}
Handling the coupling conditions on $\Gamma$ asks for the use of
$\CS{1}$-regular virtual basis functions not only in $\Omega_2$ where
the bilaplacian problem is defined, but also in $\Omega_1$, where the
second order elliptic problem is defined.
Indeed, a simple use of $\CS{0}$-basis functions in $\Omega_1$, which
would be natural given the second order of the problem, would not
allow the imposition (or at least a simple imposition) of the gluing
condition on the normal derivatives.

\subsubsection{Convergence results}\label{subsec:error_energy}

The following convergence result in the energy norm holds (see
\cite{Antonietti-Manzini-Verani:2019} for the proof).

\begin{theorem}
  \label{theorem:poly:energy:convg:rate}
  Let $\fs\in\HS{r-p+1}(\Omega)$ be the forcing term at the right-hand
  side, $\us$ the solution of the variational
  problem~\eqref{eq:poly:pblm:wp} and $\ush\in\Vhrp{r}$ the solution
  of the virtual element method~\eqref{eq:poly:VEM}.
  Then, it holds that
  \begin{align}
    \normV{\us-\ush}
    \leq\Cs\hh^{r-(p-1)}\big( \snorm{\us}{r+1} + \snorm{\fs}{r-p+1} \big).
    \label{eq:poly:energy:convg:rate}
  \end{align}
\end{theorem}

Moreover, the following convergence results in lower order norms can
established~\cite{Antonietti-Manzini-Verani:2019}.
\begin{theorem}[Even $p$, even norms]
  Let $\fs\in\HS{r-p+1}(\Omega)$, $\us$ the solution of the
  variational problem~\eqref{eq:poly:pblm:wp} with $p=2\ell$ and
  $\vsh\in\Vhrp{r}$ the solution of the virtual element
  method~\eqref{eq:poly:VEM}.
  Then, there exists a positive constant $\Cs$ independent of $\hh$
  such that
  \begin{align}
    \snorm{\us-\ush}{2i}
    \leq\Cs\hh^{r+1-2i}\Big(\snorm{\us}{r+1}+\snorm{\fs}{r-(p-1)}\Big),
  \end{align}
  for every integer $i=0,\ldots,\ell-1$.
\end{theorem}
\begin{theorem}[Even $p$, odd norms]
  Let $\fs\in\HS{r-p+1}(\Omega)$, and $\us$ the solution of the
  variational problem~\eqref{eq:poly:pblm:wp} with $p=2\ell$ and
  $\ush\in\Vhrp{r}$ the solution of the virtual element
  method~\eqref{eq:poly:VEM}.
  Then, there exists a positive constant $\Cs$ independent of $\hh$
  such that
  \begin{align}
    \snorm{\us-\ush}{2i+1}
    \leq\Cs\hh^{(r+1)-(2i+1)}\Big(\snorm{\us}{r+1}+\snorm{\fs}{r-(p-1)}\Big),
  \end{align}
  for every integer $i=0,\ldots,\ell-1$.
\end{theorem}
\begin{theorem}[Odd $p$, even norms]
  Let $\us$ be the solution of the variational
  problem~\eqref{eq:poly:pblm:wp} and $\ush\in\Vhrp{r}$
  the solution of the virtual element
  method~\eqref{eq:poly:VEM}.
  Then, there exists a positive constant $\Cs$ independent of $\hh$
  such that
  \begin{align}
    \snorm{\us-\ush}{2i}
    \leq\Cs\hh^{ (r+1)-2i }\Big(\snorm{\us}{r+1} + \snorm{\fs}{r-(p-1)}\Big),
  \end{align}
    for every integer $i=0,\ldots,\ell-1$.
\end{theorem}
\begin{theorem}[Odd $p$, odd norms]
  Let $\us$ be the solution of the variational
  problem~\eqref{eq:poly:pblm:wp} and $\ush\in\Vhrp{r}$ the solution
  of the virtual element method~\eqref{eq:poly:VEM}.
  Then, there exists a positive constant $\Cs$ independent of $\hh$
  such that
  \begin{align}
    \snorm{\us-\ush}{2i+1}
    \leq\Cs\hh^{(r+1)-(2i+1)}\Big(\snorm{\us}{r+1}+\snorm{\fs}{r-(p-1)}\Big),
  \end{align}
  for every integer $i=0,\ldots,\ell-1$.
\end{theorem}

%% section-3: The virtual element method for the Cahn-Hilliard problem
%%\section{The virtual element method for the Cahn-Hilliard problem}
%%\label{sec:Cahn-Hilliard}
%% \input{Sec3_VEM_Cahn_Hilliard.tex}
%%
%% operators to be defined
%%
%% \Pi^{\nabla}
%% \Pi^{\Delta}
%% \Pi^{0}
%%
%% symbols
%% use \V for the vertices instead of \nu
%% use $\Nu_\P$ for the set of vertices of eleemnt\P

\section{The virtual element method for the Cahn-Hilliard problem}
\label{sec:Cahn-Hilliard}

\subsection{The continuous problem}
\label{sec:cont}
%-------------------------------------------------------------------
Let $\Omega\subset\REAL^2$ be an open, bounded domain with polygonal
boundary $\Gamma$, $\psi:\REAL\to\REAL$ with $\psi(x)=(1-x^2)^2/4$ and
$\phi(x)=\psi^\prime(x)$.
We consider the Cahn-Hilliard problem:
\textit{Find $\us(\xs,t):\Omega\times[0,\Ts]\rightarrow\REAL$ such
  that:}
\begin{subequations}
  \begin{align}
    \dot\us - \Delta\big(\phi(\us) - \gamma^2\Delta\us\big)              = 0\phantom{\us_0(\cdot)} &\quad\textrm{in~}\Omega\times[0,T],          \label{eq:CHforte:A}\\[0.3em]
    \us(\cdot,0)                                                         = \us_0(\cdot)\phantom{0} &\quad\textrm{in~}\Omega,                    \label{eq:CHforte:B}\\[0.1em]
    \partial_{n}\us = \partial_{ n}\big(\phi(\us) - \gamma^2\Delta\us\big) = 0\phantom{\us_0(\cdot)} &\quad\textrm{on~}\partial\Omega\times[0,\Ts],\label{eq:CHforte:C}
  \end{align}
\end{subequations}
where $\partial_{ n}$ denotes the (outward) normal derivative and
$\gamma\in\REAL^+$, $0<\gamma\ll1$, represents the interface
parameter.
On the domain boundary we impose a no flux-type condition on $\us$ and
the chemical potential $\phi(\us)-\gamma^2\Delta\us$.

\medskip
To define the variational formulation of
problem~\eqref{eq:CHforte:A}-\eqref{eq:CHforte:C} we introduce the
three bilinear forms:
\begin{equation*}
  \begin{aligned}
    \as^\Delta(\vs,\ws) &= \int_\Omega (\nabla^2\vs) : (\nabla^2\ws)\,\dxv  & \quad\forall\vs,\ws\in\HS{2}(\Omega),\\[0.5em]
    \as^\nabla(\vs,\ws) &= \int_\Omega \nabla\vs\cdot\nabla\ws\,\dxv       & \quad\forall\vs,\ws\in\HONE(\Omega),\\[0.5em]
    \as^0(\vs,\ws)     &= \int_\Omega \vs\,\ws\,\dxv                      & \quad\forall\vs,\ws\in\LTWO(\Omega),
  \end{aligned}
\end{equation*}
($\nabla^2$ being the Hessian operator) and the semi-linear form
\begin{align*}
  \rs(\zs;\vs,\ws) = \int_\Omega \phi'(\zs)\nabla\vs\cdot\nabla\ws\,\dxv\quad\forall\zs,\vs,\ws\in\HS{2}(\Omega).
\end{align*}
Finally, introducing the functional space
\begin{equation}\label{eq:V}
  \Vs=\big\{\vs\in\HONE(\Omega)\,:\,\partial_{ n}{\vs}=0\textrm{~on~}\Gamma\big\},
\end{equation}
which is a subspace of $\HONE(\Omega)$.

\medskip
The weak formulation of problem
\eqref{eq:CHforte:A}-\eqref{eq:CHforte:C} reads as:
\textit{Find $u(\cdot,t)\in\Vs$ such that}
\begin{subequations}
  \begin{align}
    & \as^0(\dot\us,\vs) + \gamma^2\as^\Delta(\us,\vs) + \rs(\us;\us,\vs) = 0 \quad\forall\vs\in\Vs,\label{eq:contpbl:A}\\[0.2em]
    & \us(\cdot,0) =\us_0.\label{eq:contpbl:B}
  \end{align}
\end{subequations}

\subsection{The conforming Virtual Element approximation}

In this section, we introduce the main building blocks for the
conforming virtual discretization of the Cahn-Hilliard equation,
report a convergence result and collect some numerical results
assessing the theoretical properties of the proposed scheme.

\subsubsection{A $C^1$ Virtual Element space}
\label{sec:vem1}
%-------------------------------------------------------------------
We briefly recall the construction of the virtual element space
$\Vhcahn\subset\HS{2}(\Omega)$ that we use to discretize
\eqref{eq:contpbl:A}-\eqref{eq:contpbl:B}; see
\cite{Antonietti-BeiraodaVeiga-Scacchi-Verani:2016} for more details.

\medskip
Given an element $\P\in\Th$, the \emph{augmented} local space $\VEt$
is defined by
\begin{align}
  \VEt
  = \Big\{\vs\in\HS{2}(\P)\,:\,
  &\Delta^2\vs\in\PS{2}(\P),
  \restrict{\vs}{\partial\P}\in\CS{0}(\partial\P),\restrict{\vs}{\E}\in\PS{3}(\E)\quad\forall\E\in\partial\P,
  \nonumber\\
  &\restrict{\nabla\vs}{\partial\P}\in\big[\CS{0}(\partial\P)\big]^2,\partial_{ n}\restrict{\vs}{\E}\in\PS{1}(\E)\quad\forall\E\in\partial\P
  \Big\},
  \label{eq:Vtilde}
\end{align}
with $\partial_{n}$ denoting the (outward) normal derivative.

\begin{remark} The space $\VEt$ corresponds to the space 
  $\tVhPrp{2p-1}$ with $p=2$ introduced in Section
  \ref{subsec:modified_VE_space}.
\end{remark}

We consider the two sets of linear operators from $\VEt$ into
$\mathbb{R}$ denoted by \TERM{D1} and \TERM{D2} and defined as
follows:
\begin{description}
\item[]\TERM{D1} contains linear operators evaluating $\vsh$ at the
  $n=n(\P)$ vertices of $\P$;
\item[]\TERM{D2} contains linear operators evaluating $\nabla\vsh$ at
  the $n=n(\P)$ vertices of $\P$.
\end{description}

The output values of the two sets of operators \TERM{D1} and \TERM{D2}
are sufficient to uniquely determine $\vsh$ and $\nabla\vsh$ on the
boundary of $\P$ (cf. Section \ref{subsec:modified_VE_space}).

We use of the following local bilinear forms for all $\P\in\Omega_h$
\begin{align}
  \asP^\Delta(\vs,\ws) &= \int_\P(\nabla^2\vs):(\nabla^2\ws)\,\dxv
  \phantom{\int_\P \vs\,\ws\,\dxv}\hspace{-1cm}
  \forall\vs,\ws\in\HS{2}(\P),\label{eq:loc-forms-cont} \\
  \asP^\nabla(\vs,\ws) &= \int_\P \nabla\vs\cdot\nabla\ws\,\dxv       \quad\forall\vs,\ws\in\HONE(\P), \\
  \asP^0(\vs,\ws)     &= \int_\P \vs\,\ws\,\dxv
  \phantom{\int_\P(\nabla^2\vs):(\nabla^2\ws)\,\dxv}\hspace{-1cm}
  \forall\vs,\ws\in\LTWO(\P).
\end{align}

Now, we introduce the elliptic projection operator
$\PiDr{2}\colon\VEt\rightarrow\PS{2}(\P)$ defined by
\begin{align}
  \asP^\Delta(\PiDr{2}\vsh,\qs) &= \asP^\Delta(\vsh,\qs)
  \phantom{(\!(\vsh,\qs)\!)_\P}
  \forall\qs\in\PS{2}(\P),\label{wq:pi_delta}\\
  (\!( \PiDr{2}\vsh,\qs)\!)_\P &= (\!(\vsh,\qs)\!)_\P
  \phantom{\asP^\Delta(\vsh,\qs)}
  \forall\qs\in\PS{1}(\P),%\label{pi_delta}
\end{align}
for all $\vsh\in\VEt$ where $(\!(\cdot,\cdot)\!)_\P$ is the Euclidean
scalar product acting on the vectors that collect the vertex function
values, i.e.
\begin{align*}
  (\!(\vsh,\wsh)\!)_\P = \sum_{\V\in\VRTX{\P}} 
  \!\!\vsh(\V)\:\wsh(\V) \quad\forall\vsh,\wsh\in\CS{0}(\P) .
\end{align*}
As shown in \cite{Antonietti-BeiraodaVeiga-Scacchi-Verani:2016}, the
operator $\PiDr{2}\colon\VEt\rightarrow\PS{2}(\P)$ is well
defined and uniquely determined on the basis of the information
carried by the linear operators in \TERM{D1} and \TERM{D2}.

Hinging upon the augmented space $\VEt$ and employing the projector
$\PiDr{2}$ we define our virtual local space
\begin{align}
  \VE = \big\{
  \vs\in\VEt\,:\,\asP^0(\PiDr{2}(\vs),\qs)  = \asP^0 (\vs,\qs)\quad\forall\qs\in\PS{2}(\P)
  \big\}.
  \label{eq:VE}
\end{align}
Since $\VE\subset\VEt$, operator $\PiDr{2}$ is well defined
on $\VE$ and computable by using the values provided by \TERM{D1} and
\TERM{D2}.
Moreover, the set of operators \TERM{D1} and \TERM{D2} constitutes a
set of degrees of freedom for the space $\VE$.
Finally, there holds $\PS{2}(\P)\subseteq\VE$.

\medskip
We now introduce two further projectors on the local space $\VE$,
namely $\Pizr{2}$ and $\PiPr{2}$, that will be employed together
with the above projector $\PiDr{2}$ to build the discrete
counterparts of the bilinear forms
in~\eqref{eq:loc-forms-cont}.
Operator $\Pizr{2}:\VE\rightarrow\PS{2}(\P)$ is the standard $\LTWO$
projector on the space of quadratic polynomials in $\P$.
This is computable by means of the values of the degrees of freedom
\TERM{D1} and \TERM{D2}
(cf. \cite{Antonietti-BeiraodaVeiga-Scacchi-Verani:2016}).
To define $\PiPr{2}:\VE\rightarrow\PS{2}(\P)$ we need the
additional bilinear form
$\as^{\nabla}(\cdot,\cdot):\VE\times\VE\to\REAL$ that is given by
\begin{align*}
  \as^\nabla(\vs,\ws) = \int_\Omega \nabla\vs\cdot\nabla\ws\,\dxv \quad\forall\vs,\ws\in\HONE(\Omega).
\end{align*}
Operator $\PiPr{2}$ is the elliptic projection defined with
respect to $\as^{\nabla}(\cdot,\cdot)$:
\begin{subequations}
  \begin{align}
    \asP^\nabla(\PiPr{2}\vsh,\qs) &= \asP^\nabla(\vsh,\qs)\quad\forall\qs\in\PS{2}(\P),\label{eq:pi_nabla:0}\\
    \int_\P \PiPr{2}\vsh\dxv     &= \int_\P \vsh\,\dxv.\label{eq:pi_nabla:1}
  \end{align}
\end{subequations}
Such operator is well defined and uniquely determined by the values of
\TERM{D1} and
\TERM{D2}~\cite{Antonietti-BeiraodaVeiga-Scacchi-Verani:2016}.

We are now ready to introduce the global virtual element space, which
defined as follows
\begin{align*}
  \Vhcahn =
  \big\{\vs\in\Vs\,:\,\restrict{\vs}{\P}\in\VE\quad\forall\P\in\Th\big\}.
\end{align*}
The virtual element functions in $\Vhcahn$ and their gradients are
continuous fields on $\Omega$, so this functional space is a
conforming subspace of $\HS{2}(\Omega)$.
The \emph{global degrees of freedom} of $\Vhcahn$ are obtained by
collecting the elemental degrees of freedom, so the dimension of
$\Vhcahn$ is three times the number of the mesh vertices, and every
virtual element function $\vsh$ defined on $\Omega$ is uniquely
determined by
\begin{itemize}
\item[]$(i)$~its values at the mesh vertices;
\item[]$(ii)$~its gradient values at the mesh vertices.
\end{itemize}
Finally, we recommended to scale the degrees of freedom \TERM{D2} by
some local characteristic mesh size $\hV$ in order to obtain a
better condition number of the final system.

% -------------------------------------------------------------------
\subsubsection{Virtual element bilinear forms}
\label{sec:vem2}
% -------------------------------------------------------------------
We start by introducing the discrete versions of the elemental
bilinear form forms in \eqref{eq:loc-forms-cont}.
Let $\P\in\Th$ be a generic mesh element and
$\ssP(\cdot,\cdot):\VE\times\VE\to\REAL$ the positive definite
bilinear form given by:
\begin{align*}
  \ssP(\vsh,\wsh) = \sum_{\V\in\VRTX{\P}} 
  \!\! \Big(\vsh(\V)\:\wsh(\nu) + \hV^2 \: \nabla\vsh(\V)\cdot\nabla\wsh(\V)\Big)
  \quad\forall\vsh,\wsh\in\VE,
\end{align*}
where $\hV$ is a characteristic mesh size length associated with node
$\V$, e.g., the maximum diameter among the elements having $\V$ as a
vertex.

Recalling~\eqref{eq:loc-forms-cont}, we consider the virtual element
bilinear forms:
\begin{align}
  \ashP^\Delta(\vsh,\wsh) &= \asP^\Delta(\PiDr{2}\vsh, \PiDr{2}\wsh) + \hP^{-2} \ssP\big(\vsh - \PiDr{2}\vsh,\wsh - \PiDr{2}\wsh\big),\label{eq:loc-discr-forms:Delta}\\
  \ashP^\nabla(\vsh,\wsh) &= \asP^\nabla(\PiPr{2}\vsh, \PiPr{2}\wsh) +         \ssP\big(\vsh - \PiPr{2}\vsh,\wsh - \PiPr{2}\wsh),    \label{eq:loc-discr-forms:nabla}\\
  \ashP^0    (\vsh,\wsh) &= \asP^0(\Pizr{2}\vsh, \Pizr{2}\wsh)            + \hP^2    \ssP\big(\vsh - \Pizr{2}\vsh,    \wsh - \Pizr{2}\wsh\big)    \label{eq:loc-discr-forms:0}
\end{align}
for all $\vsh$, $\wsh\in\VE$.
Under the mesh regularity conditions of
Section~\ref{subsec:mesh:assumptions}, we can prove the consistency
and stability of the discrete bilinear forms.
Let the symbol $\dagger$ stands for ``$\Delta$'', ``$\nabla$'' or
``$0$''.
We have:
\begin{itemize}
\item[]\TERM{A}~\emph{(polynomial consistency)}
  $\ashP^\dagger(\ps,\vsh)=\asP^\dagger(\ps,\vsh)\quad\forall\ps\in\PS{2}(\P),\,\vsh\in\VE$;

  \smallskip
\item[]\TERM{B}~\emph{(stability)} there exist two positive constants
  $\cs_*$ and $\cs^*$ independent of $\hh$ and the element $\P\in\Th$
  such that
  \begin{align*}
    \cs_*\asP^\dagger(\vsh,\vsh)\leq\ashP^\dagger(\vsh,\vsh)\leq\cs^*\asP^\dagger(\vsh,\vsh)
    \quad\forall\vsh\in\VE.
  \end{align*}
\end{itemize}

A consequence of the above properties is that the bilinear form
$\ashP^\dagger(\cdot,\cdot)$ is continuous with respect to the
relevant norm, which is $\HS{2}$ for~\eqref{eq:loc-discr-forms:Delta},
$\HONE$ for~\eqref{eq:loc-discr-forms:nabla}, and $\LTWO$ for
\eqref{eq:loc-discr-forms:0}.
For every choice of $\dagger$, the corresponding global bilinear form
is
\begin{align*}
  \ash^\dagger(\vsh,\wsh)=\sum_{\P\in\Th}\ashP^\dagger(\vsh,\wsh)\quad\forall\vsh,\wsh\in\Vhcahn.
\end{align*}\eqref{eq:CHforte:A}-\eqref{eq:CHforte:C}

We now turn our attention to the semilinear form
$\rs(\cdot;\cdot,\cdot)$, which we can also write as the sum of
elemental contributions:\\
\setlength{\belowdisplayskip}{0.3em}%
\setlength{\abovedisplayskip}{0.3em}%
\begin{align*}
  & \rs (\zs;\vs,\ws) = \sum_{\P\in\Th}\rsP(\zs;\vs,\ws) \quad \forall\zs,\vs\ws\in\HS{2}(\Omega)
  \intertext{where}
  & \rsP(\zs;\vs,\ws) = \int_\P (3\zs^2-1) \nabla\vs\cdot\nabla\ws\dxv\quad\forall\P\in\Th.
\end{align*}

% $\phi''(x) = 3x^2-1 $
%%
\medskip
On each element $\P$, we approximate the term $\zs(\xs)^2$ by means of
its cell average, which we compute using the $\LTWO(\P)$ bilinear form
$\ashP^0(\cdot,\cdot)$:
\begin{align*}
  \restrict{\zsh^2}{\P}\approx\mP^{-1}\ashP^0(\zsh,\zsh) ,
\end{align*}
where we recall that $\mP$ is the area of element $\P$.
This approach has the correct approximation properties and
preserves the positivity of $\zs^2$.

We therefore propose the following approximation of the local
nonlinear forms
\begin{align*}
  \rshP(\zsh;\vsh,\wsh) = \restrict{\widehat{\phi^\prime(\zsh)}}{\P}\:\ashP^\nabla(\vsh,\wsh)
  \quad\forall\zsh,\vsh,\wsh\in\VE, 
\end{align*}
where 
$\restrict{\widehat{\phi^\prime(\zsh)}}{\P} = 3\mP^{-1}\ashP^0(\zsh,\zsh) - 1$.
The global form is then assembled as 
\begin{align*}
  \rsh(\zsh;\vsh,\wsh) = \sum_{\P\in\Th}\rshP(\zsh;\vsh,\wsh)
  \quad\forall\zsh,\vsh\wsh\in\Vhcahn.
\end{align*}

%% -------------------------------------------------------------------
\subsubsection{The discrete problem}
\label{sec:vem3}
%% -------------------------------------------------------------------
%%
The virtual element discretization of problem
\eqref{eq:contpbl:A}\eqref{eq:contpbl:B} follows a Galerkin approach
in space combined with a backward Euler time-stepping scheme.
Consider the functional space
\begin{align*}
  \Vhbc
  = \Vhcahn\cap\Vs
  =  \big\{ \vs\in\Vhcahn\,:\,\partial_{ n}\vs=0\textrm{~on~}\partial\Omega\big\},
\end{align*}
which includes the boundary conditions.
Then, we introduce the the semi-discrete approximation: \textit{Find
  $\ush(\cdot,t)$ in $\Vhbc$ such that} \\
\begin{align}
  \ash^0(\dot{\ush},\vsh) + \gamma^2 \ash^\Delta(\ush,\vsh) + \rsh(\ush;\ush,\vsh) &= 0  \quad \forall\vsh\in\Vhbc,\label{eq:discrpbl:A}\\[0.25em]
  \ush(0,\cdot)=\us_{0,\hh}(\cdot),\label{eq:discrpbl:B}
\end{align}
where $\us_{0,\hh}$ is a suitable approximation of $\us_0$ in $\Vhbc$
and $\ash^0(\cdot,\cdot)$, $\ash^\Delta$ and $\rsh$ are the virtual
element bilinear forms defined in the previous section.

\medskip
To formulate the fully discrete scheme, we subdivide the time interval
$[0,T]$ into $\Ns$ uniform sub-intervals of length $k=T/N$ by means of
the time nodes $0=\ts_0<\ts_1<\ldots<\ts_{\Ns-1}<\ts_{\Ns}=\Ts$, and
denote the virtual element approximation of the solution $\us(\cdot,\ts)$ at $\us(\cdot,\ts^{i})$ in
$\Vhbc$ by $\us_{\hh,k}^{i}$.
The fully discrete problem reads as: \textit{Given
  $\us_{hk}^0=\us_{0,h}\in\Vhbc$, find $\us_{hk}^i\in\Vhbc$,
  $i=1,\ldots,\Ns$ such that}
\begin{align}
  k^{-1}\ash^0(\us_{hk}^{i}-\us_{hk}^{i-1},\vsh)
  + \gamma^2\ash^\Delta(\us_{hk}^{i},\vsh) + \rsh(\us_{hk}^{i},\us_{hk}^{i};\vsh) = 0
  \quad\forall\vsh\in\Vhbc.
  \label{eq:fullydisc}
\end{align}

The semidiscrete Virtual Element formulation given in
\eqref{eq:discrpbl:A}-\eqref{eq:discrpbl:B} converges to the exact
solution of problem \eqref{eq:contpbl:A}-\eqref{eq:contpbl:B}
according to the result stated in this theorem and proved
in~\cite{Antonietti-BeiraodaVeiga-Scacchi-Verani:2016}.
\begin{theorem}\label{conv_theo}
  Let $\us$ be the solution of problem
  \eqref{eq:contpbl:A}-\eqref{eq:contpbl:B}.
  Let $\ush$ be the virtual element approximation provided
  by~\eqref{eq:discrpbl:A}-\eqref{eq:discrpbl:B} and assume that
  \begin{align*}
    \norm{\ush}{\LINF(\Omega)}\leq\Cs
  \end{align*}
  for all $t\in(0,T]$ and some positive constant $\Cs$ independent of $h$.
  Then, it holds that
  \begin{align*}
    \NORM{\us - \ush}{\LTWO(\Omega)} \lesssim \hh^2
  \end{align*}
  for every $\ts\in[0,T]$.
\end{theorem}

% -------------------------------------------------------------------------------
\subsection{Numerical results}
\label{sec:num}
% -------------------------------------------------------------------------------
In this test, taken
from \cite{Antonietti-BeiraodaVeiga-Scacchi-Verani:2016} we study the
convergence of our VEM discretization applied to the Cahn-Hilliard
problem with a load term $\fs$ obtained by enforcing as exact solution
$\us(x,y,t)=t\cos(2\pi\xs)\cos(2\pi\ys)$.
The parameter $\gamma$ is set to $1/10$ and the time step size
$\Delta\ts$ is $1e-7$.
The $\HTWO$, $\HONE$ and $\LTWO$ errors are computed at $t=0.1$ on
four quadrilateral meshes discretizing the unit square.
The time discretization is performed by the Backward Euler method.
The resulting non-linear system \eqref{eq:fullydisc} at each time step
is solved by the Newton method, using the $l^2$ norm of the relative
residual as a stopping criterion.
The tolerance for convergence is $1e-6$.

\begin{table}[htbp]
  \caption{ $\HTWO$, $\HONE$ and $\LTWO$ errors and convergence rates
    $\alpha$ computed on four quadrilateral meshes discretizing the
    unit square~\cite{Antonietti-BeiraodaVeiga-Scacchi-Verani:2016}.}
  \begin{center}
    \begin{tabular}{c|cc|cc|cc}
      \hline
      $h$    &$|\us-\ush|_{H^2(\Omega)}$  &$\alpha$  &$|\us-\ush|_{H^1(\Omega)}$  &$\alpha$  &$\NORM{\us-\ush}{L^2(\Omega)}$  &$\alpha$\\
      \hline
      1/16   &1.35e-1          &--        &8.57e-2          &--        &8.65e-2            &--\\
      1/32   &5.86e-2          &1.20      &2.20e-2          &1.96      &2.20e-2            &1.97\\
      1/64   &2.79e-2          &1.07      &5.53e-3          &1.99      &5.52e-3            &1.99\\
      1/128  &1.38e-2          &1.02      &1.37e-3          &2.01      &1.37e-3            &2.01\\
      \hline
    \end{tabular}
    \label{tab_test1}
  \end{center}
\end{table}

\begin{figure}[t]
  \begin{center}
    \includegraphics[width=3.45cm]{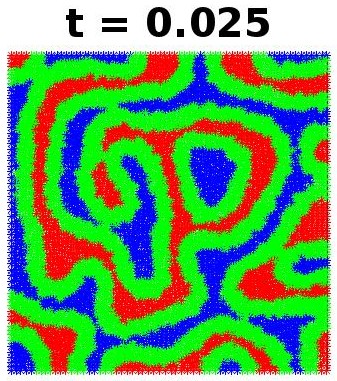} % {figures-JPG/fig_CH_t0025.jpg}
    \includegraphics[width=3.5cm ]{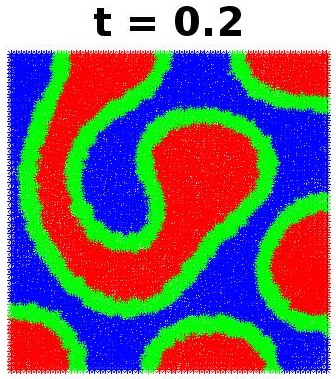} % {figures-JPG/fig_CH_t02.jpg}
    \includegraphics[width=3.5cm ]{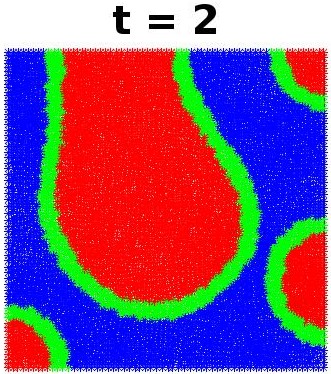} % {figures-JPG/fig_CH_t2.jpg}
    \caption{Spinoidal decomposition on the unit square at three temporal frames
      for a Voronoi polygonal mesh of 4096 elements.
%      meshes)~\cite{Antonietti-BeiraodaVeiga-Scacchi-Verani:2016}.
    }\label{spin_paul}
  \end{center}
\end{figure}

The results reported in Table~\ref{tab_test1} show that the VEM method
converges is convergent with a convergence rate close to $2$ in the
$\LTWO$ norm as expected from Theorem~\ref{conv_theo}.
In the $\HTWO$ and $\HONE$ seminorms, the method converges with order
$1$ and $2$ respectively, as we can expect from the FEM theory and the
approximation properties of the virtual element space. Finally, in Figure \ref{spin_paul} we report the results of a
spinoidal decomposition.
For completeness, we recall that spinoidal decomposition is a physical
phenomenon consisting of the separation of a mixture of two or more
components to bulk regions of each, which occurs when a
high-temperature mixture of different components is rapidly cooled. We employ an initial datum $\us_0$ chosen to be a uniformly
distributed random perturbation between $-1$ and $1$.
Results are consistent with the literature,
cf. \cite{Antonietti-BeiraodaVeiga-Scacchi-Verani:2016}.
%%for a more detailed discussion and further numerical results).
%%

%% section-4: The virtual element method for the elastodynamics problem
%% \input{Sec4_VEM_Elastodynamics.tex}
\section{The virtual element method for the elastodynamics problem}
\label{sec:elastodynamics}

\subsection{The continuous  problem}
\label{sec2:model}

We consider an elastic body occupying the open, bounded polygonal
domain $\Omega\subset\REAL^2$ with Lipschitz boundary $\Gamma$.
We assume that boundary $\Gamma$ can be split into the two disjoint
subsets $\GamD$ and $\GamN$, so that $\Gamma=\GamD\cup\GamN$ and with
the one-dimensional Lebesgue measure (length)
$\ABS{\GamD\cap\GamN}=0$.
For the well-posedness of the mathematical model, we further require
length of $\GamD$ is nonzero, i.e., $\ABS{\GamD}>0$.
Let $T>0$ denote the final time.
We consider
the external load $\fv\in\LTWO\big(0,T;[\LTWO(\Omega)]^2\big)$,
the boundary function $\gvN\in\CS{1}\big(0,T;[\HS{1/2}_{0,\GamN}]^2\big)$,
and the initial functions $\uv_{0}\in[H^{1}_{0,\GamD}(\Omega)]^2$,
$\uv_{1}\in[L^2(\Omega)]^2$.
For such time-dependent vector fields, we may indicate the dependence
on time explicitly, e.g., $\fv(t):=\fv(\cdot,t)\in[\LTWO(\Omega)]^2$,
or drop it out to ease the notation when it is obvious from the
context.

The equations governing the two-dimensional initial/boundary-value
problem of linear elastodynamics for the displacement vector
$\uv:\Omega\times[0,T]\to\REAL^2$ are:
\begin{align}
  \rho\ddot{\uv} - \nabla\cdot\bsig(\uv)
  &= \fv\phantom{\zero\gvN\uv_{0}\uv_{1}}\textrm{in}~\Omega\times(0,T],\label{eq:pblm:strong:A}\\[0.5em]
  \uv
  &= \zero\phantom{\fv\gvN\uv_{0}\uv_{1}}\textrm{on}~\GamD \times(0,T],\label{eq:pblm:strong:B}\\[0.5em]
  \bsig(\uv)\nv
  &= \gvN\phantom{\fv\zero\uv_{0}\uv_{1}}\textrm{on}~\GamN \times(0,T],\label{eq:pblm:strong:C}\\[0.5em]
  \uv
  &= \uv_{0}\phantom{\fv\zero\gvN\uv_{1}}\textrm{in}~\Omega\times\{0\},\label{eq:pblm:strong:D}\\[0.5em]
  \dot{\uv} &=
  \uv_{1}\phantom{\fv\zero\gvN\uv_{0}}\textrm{in}~\Omega\times\{0\}.\label{eq:pblm:strong:E}
\end{align}
Here, $\rho$ is the mass density, which we suppose to be a strictly
positive and uniformly bounded function and $\bsig(\uv)$ is the stress
tensor.
In~\eqref{eq:pblm:strong:B} we assume homogeneous Dirichlet boundary
conditions on $\GamD$.
This assumption is made only to ease the exposition and the analysis,
as our numerical method is easily extendable to nonhomogeneous
Dirichlet boundary conditions.

\medskip
\noindent
We denote the space of the symmetric, $2\times2$-sized, real-valued
tensors by $\symmSpace$ and assume that the stress tensor
$\bsig:\Omega\times[0,T]\to\symmSpace$ is expressed, according to
Hooke's law, by $\bsig(\uv) = \calD\beps(\uv)$,
where, $\beps(\uv)$ denotes the symmetric gradient of $\uv$ , i.e.,
$\beps(\uv)=\big(\nabla\uv+(\nabla\uv)^{T}\big)\slash{2}$, and
$\calD=\calD(\xv)\,\,:\symmSpace\longrightarrow\symmSpace$
is the \emph{stiffness} tensor
\begin{align}
  \calD\btau = 2\mu\btau + \lambda\TRACE(\btau)\matI
  \label{eq:bsig:def}
\end{align}
for all $\btau\in\mathbb{R}^{2\times2}_{\textrm{sym}}$.
In this definition, $\matI$ and $\TRACE(\cdot)$ are the identity
matrix and the trace operator; $\lambda$ and $\mu$ are the first and
second Lam\'e coefficients, which we assume to be in $\LINF(\Omega)$
and nonnegative.
The compressional (P) and shear (S) wave velocities of the medium are
respectively obtained through the relations $c_P = \sqrt{(\lambda +
  2\mu)/\rho}$ and $c_S = \sqrt{\mu/\rho}$.

\medskip
\noindent
Let $\Vv=\big[\HONEgm(\Omega)\big]^2$ be the space of $H^1$
vector-valued functions with null trace on $\Gamma_D$.
We consider the two bilinear forms
$\ms(\cdot,\cdot),\,\as(\cdot,\cdot)\,:\,\Vv\times\Vv\to\mathbb{R}$ defined
as
\begin{eqnarray}
  \ms(\wv,\vv) &=& \int_{\Omega}\rho\wv\cdot\vv\dx      \quad\forall\wv,\,\vv\in\Vv,\label{eq:ms:def}\\[0.5em]
  \as(\wv,\vv) &=& \int_{\Omega}\bsig(\wv):\beps(\vv)\dx \quad\forall\wv,\,\vv\in\Vv,\label{eq:as:def}
  %% ---
\end{eqnarray}  
  and the linear functional $\Fs(\cdot)\,:\,\Vv\to\mathbb{R}$ defined
    as
  %% ---
\begin{equation}  
  \Fs(\vv) = \int_{\Omega}\fv\cdot\vv\dx + \int_{\GamN}\gvN\cdot\vv \dS\quad\forall\vv\in\Vv .\label{eq:Fv:def}
\end{equation}

The variational formulation of the linear elastodynamics equations
reads as: \emph{For all $t\in(0,T]$ find $\uv(t)\in\Vv$ such that for
  $t=0$ it holds that $\uv(0)=\uv_0$ and $\dot{\uv}(0)=\uv_1$ and}
\begin{align}
  \ms(\ddot{\uv},\vv) + \as(\uv,\vv) &= \Fs(\vv)\qquad\forall\vv\in\Vv.
  \label{eq:weak:A}
\end{align}
As shown, for example, by Raviart and Thomas (see
Theorem~8-3.1 \cite{Raviart-Thomas:1983}) the variational
problem~\eqref{eq:weak:A} is well posed and its unique solution
satisfies
$\uv\in\CS{0}\big(0,T;\Vv\big)\cap\CS{1}\big(0,T;[\LTWO(\Omega)]^2\big)$.

\subsection{The conforming Virtual Element approximation}
    
In this section we introduce the main building blocks for the
conforming virtual element discretization of the elastodynamics
equation, report stability and convergence results and collect some
numerical results assessing the theoretical properties of the proposed
scheme.

\subsubsection{Virtual element spaces}
Let $k\geq1$ be an integer number.
The global virtual element space is defined as
\begin{align}
  \Vvhk := \Big\{\vv\in\Vv\,:\,\vv_{|_\P}\in\Vvhk(\P)\,\textrm{for~every~}\P\in\Th\Big\}
\end{align}
where $\Vvhk(\P)=\big[\Vhk(\P)\big]^2$, with 
\begin{align}
  \Vhk(\P) := \Big\{
  \vsh\in\HONE(\P)\,:\,
  &
  \vs_{\hh|_{\partial\P}}\in\CS{}(\partial\P),\,
  \vs_{\hh|_\E}\in\PS{k}(\E)\,\,\forall\E\in\partial\P,\,\,
  \Delta\vsh\in\PS{k}(\P),\nonumber\\[0.25em]
  &\big(\vsh-\Pin{k}\vsh,\mu_{\hh}\big)_{\P}=0
  \,\,\forall\mu_{\hh}\in\PS{k}(\P)\backslash{\PS{k-2}(\P)}
  \Big\},
\end{align}
where $\Pin{k}:\HONE(\P)\cap\CS{0}(\overline{\P})\to\PS{k}(\P)$ is the
usual elliptic projection of a function $\vsh$ on the space of
polynomials of degree $k$, cf. \eqref{eq:poly:Pi:A}-\eqref{eq:poly:Pi:B}.

Each virtual element function $\vsh\in\Vhk(\P)$ is uniquely
characterized by
\begin{description}
\item[]\TERM{C1} the values of $\vsh$ at the vertices of $\P$;
\item[]\TERM{C2} the moments of $\vsh$ of order up to $k-2$ on each
  one-dimensional edge $\E\in\partial\P$:
  \begin{equation}\label{eq:CF:dofs:01}
    \frac{1}{\mE}\int_{\E}\vsh\,\ms\dS,
    \,\,\forall\ms\in\calM_{k-2}(\E),\,
    \forall\E\in\partial\P;
  \end{equation}
\item[]\TERM{C3} the moments of $\vsh$ of order up to $k-2$ on $\P$:
  \begin{equation}\label{eq:CF:dofs:02}
    \frac{1}{\mP}\int_{\P}\vsh\,\ms\dx,
    \,\,\forall\ms\in\calM_{k-2}(\P).
  \end{equation}
\end{description}

As usual, the degrees of freedom of the global space $\Vvhk$ are
provided by collecting all the local degrees of freedom {(which allow
  the computation of the elliptic projection $\Pin{k}$)}, and their
unisolvence is an immediate consequence of the unisolvence of the
local degrees of freedom for the elemental spaces $\Vhk(\P)$.

\medskip
\noindent

\subsubsection{Discrete bilinear forms}
%%
%% VEM discretization, VEM bilinear form
In the virtual element setting, we define the bilinear forms
$\msh(\cdot,\cdot)$ and $\ash(\cdot,\cdot)$ as the sum of elemental
contributions, which are respectively denoted by $\mshP(\cdot,\cdot)$
and $\ashP(\cdot,\cdot)$:
\begin{align*}
  \msh(\cdot,\cdot)&\,:\,\Vvhk\times\Vvhk\to\REAL,
  \quad\textrm{with}\quad\msh(\vvh,\wvh) = \sum_{\P\in\Th}\mshP(\vvh,\wvh),\\[0.25em]
  %% ----------------
  \ash(\cdot,\cdot)&\,:\,\Vvhk\times\Vvhk\to\REAL,
  \quad\textrm{with}\quad\ash(\vvh,\wvh) = \sum_{\P\in\Th}\ashP(\vvh,\wvh).
\end{align*}

The local bilinear form $\mshP(\cdot,\cdot)$ is given by
\begin{align}
  \mshP(\vvh,\wvh) 
  = \int_{\P}\rho\Piz{k}\vvh\cdot\Piz{k}\wvh\dV 
  + \SP_{m}(\vvh,\wvh),
  \label{eq:mshP:def}
\end{align}
where $\SP_{m}(\cdot,\cdot)$ is the local stabilization term.
The bilinear form $\mshP$ depends on the orthogonal projections
$\Piz{k}\vvh$ and $\Piz{k}\wvh$, which are computable from the degrees
of freedom of $\vvh$ and $\wvh$.
The local form $\SP_{m}(\cdot,\cdot)\,:\,\Vvhk\times\Vvhk\to\REAL$ can
be \emph{any} symmetric and coercive bilinear form that is computable
from the degrees of freedom and for which there exist two strictly
positive real constants $\sigma_*$ and $\sigma^*$ such that
\begin{align}
  \sigma_*\msP(\vvh,\vvh)\leq\SP_{m}(\vvh,\vvh)\leq\sigma^*\msP(\vvh,\vvh)
  \quad\vvh\in\textrm{ker}\big(\Piz{k}\big)\cap\Vvhk(\P).
\end{align}
Computable stabilizations $\SP_{m}(\cdot,\cdot)$ are provided by
resorting to the two-dimensional stabilizations of the effective
choices for the scalar case proposed in the
literature\cite{Mascotto:2018,Dassi-Mascotto:2018}.
\medskip
The local bilinear form $\ashP$ is given by
\begin{align}
  \ashP(\vvh,\wvh) 
  = \int_{\P}\matD\Piz{k-1}(\beps(\vvh)):\Piz{k-1}(\beps(\wvh))\dV
  + \SP_a(\vvh,\wvh),
  \label{eq:ashP:def}
\end{align}
where $\SP_{a}(\cdot,\cdot)$ is the local stabilization term.
The bilinear form $\ashP$ depends on the orthogonal projections
$\Piz{k-1}\nabla\vvh$ and $\Piz{k-1}\nabla\wvh$, which are computable
from the degrees of freedom of $\vvh$ and $\wvh$.
On its turn, $\SP_{a}(\cdot,\cdot)\,:\,\Vvhk\times\Vvhk\to\REAL$ can
be \emph{any} symmetric and coercive bilinear form that is computable
from the degrees of freedom and for which there exist two strictly
positive real constants $\overline{\sigma}_*$ and
$\overline{\sigma}^*$ such that
\begin{align}
  \overline{\sigma}_*\asP(\vvh,\vvh)
  \leq\SP_{m}(\vvh,\vvh)\leq
  \overline{\sigma}^*\asP(\vvh,\vvh)
  \quad\vvh\in\textrm{ker}\big(\Piz{k}\big)\cap\Vvhk(\P).
\end{align}
Moreover, the bilinear form $\SP_a(\cdot,\cdot)$ must scale with
respect to $\hh$ like $\asP(\cdot,\cdot)$, i.e., as $\mathcal{O}(1)$.
As before, we can define computable stabilizations
$\SP_{a}(\cdot,\cdot)$ by resorting to the two-dimensional
stabilizations for the scalar case proposed in the
literature~\cite{Mascotto:2018,Dassi-Mascotto:2018}.
\medskip
As usual, the discrete bilinear forms $\ashP(\cdot,\cdot)$ and
$\mshP(\cdot,\cdot)$ satisfy the \emph{$k$-consistency} and
\emph{stability} properties.
The stability constants may depend on physical parameters and the
polynomial degree
$k$~\cite{BeiraodaVeiga-Chernov-Mascotto-Russo:2016,Antonietti-Manzini-Mazzieri-Mourad-Verani:2020}.

\subsubsection{Discrete load term}
We approximate the right-hand side~\eqref{eq:VEM:semi-discrete} of the
variational formulation by means of the linear functional
$\Fsh(\cdot):\Vvhk\to\REAL^2$ given by
\begin{align}
  \Fsh(\vvh) 
  = \int_{\Omega}\fv\cdot\Piz{k-2}(\vvh)\dV 
  + \sum_{\E\in\GamN}\int_{{\E} }\gvN\cdot{ {\vvh} }\dS
  \quad\forall\vvh\in\Vvhk.
  \label{eq:Fsh:def}
\end{align}
%%
%where $\Pi^{0,\E}_{k}(\vvh)$ is the $\LTWO$-orthogonal projection of
%the components of the virtual element vector field $\vvh$ on the space
%of scalar polynomials defined on edge $\E$.
%% 
The linear functional $\Fsh(\cdot)$ is clearly computable since the
edge trace $\restrict{\vvh}{\E}$ is a known polynomial and
$\Piz{k}(\vvh)$ is computable from the degrees of freedom of $\vvh$.
Moreover, $\Fsh(\cdot)$ is a bounded functional.
In fact, when $\gvN=0$ using the stability of the projection operator
and the Cauchy-Schwarz inequality, we note that
\begin{align}
  \ABS{\Fsh(\vvh)}
  &\leq\left\vert
  \int_{\Omega}\fv(t)\cdot\Piz{k-2}(\vvh)\dV
  \right\vert
  \leq \NORM{\fv(t)}{0}\NORM{\Piz{k-2}(\vvh)}{0}
  \nonumber\\[0.5em]
  &\leq \NORM{\fv(t)}{0}\NORM{\vvh}{0}
  \quad\forall t\in[0,T].
  \label{eq:Fsh:stab}
\end{align}
This estimate is used in the proof of the stability of the
semi-discrete virtual element approximation { (see Theorem \ref{theorem:semi-discrete:stability})}.

\subsubsection{The discrete problem}
The semi-discrete virtual element approximation of~\eqref{eq:weak:A}
reads as: \emph{For all $t\in(0,T]$ find $\uvh(t)\in\Vvhk$ such that
for $t=0$ it holds that $\uvh(0)=(\uv_{0})_{\INTP}$ and
$\dot{\uv}_{\hh}(0)= (\uv_{1})_{\INTP}$ and}
\begin{align}
  \msh(\ddot{\uv}_h,\vvh) 
  + \ash(\uvh,\vvh) = \Fsh(\vvh)
  \quad\forall\vvh\in\Vvhk.
  \label{eq:VEM:semi-discrete}
\end{align}
Here, 
$\uvh(t)$ is the virtual element approximation of $\uv$ and $\vvh$ is
the generic test function in $\Vvhk$, while $(\uv_{0})_{\INTP}$ and
$(\uv_{1})_{\INTP}$ are the virtual element interpolants of the
initial solution functions $\uv(0)$ and $\dot{\uv}(0)$.

%% fully discrete scheme
We carry out the time integration by applying the leap-frog time
marching scheme~\cite{Quarteroni-Sacco-Saleri:2007} to the second
derivative in time $\ddot{\uv}_{\hh}$.
To this end, we subdivide the interval $(0,T]$ into $N_T$ subintervals
of amplitude $\Delta t=T\slash{N_T}$ and at every time level
$t^n=n\Delta t$ we consider the variational problem for $n\geq1$:
\begin{multline}
  \msh(\uvh^{n+1},\vvh) 
  - 2\msh(\uvh^{n},\vvh) 
  +  \msh(\uvh^{n-1},\vvh)
  + \Delta t^2\ash(\uvh^{n},\vvh)
  \\[0.5em]
  = \Delta t^2\Fsh^{n}(\vvh)
  \quad\forall\vvh\in\Vvhk,
  \label{eq:VEM:fully-discrete}
\end{multline}
and initial step
\begin{multline*}
  \msh(\uvh^{1},\vvh) 
  - \msh(\uv_{0},\vvh) 
  - \Delta t  \msh(\uv_{1},\vvh) 
  + \frac{\Delta t^2}{2}\ash(\uv_{0},\vvh)
  \\[0.5em]
  = \frac{\Delta t^2}{2}\Fsh^{0}(\vvh)
  \quad\forall\vvh\in\Vhk.
\end{multline*}
The leap-frog scheme is second-order accurate, explicit and
conditionally stable.~\cite{Quarteroni-Sacco-Saleri:2007}
It is straightforward to show that these properties are inherited by
the fully-discrete scheme~\eqref{eq:VEM:fully-discrete}.

\subsubsection{Stability and convergence analysis for the semi-discrete problem}
\label{sec4:convergence}
We employ the \emph{energy} norm
\begin{align}
  \TNORM{\vvh(t)}{}^2 
  = \NORM{ \rho^{\frac{1}{2}}\dot{\vv}_{\hh}(t) }{0}^2 
  + \abs{ \vvh(t) }_{1}^2,
  \qquad t\in[0,T],
  \label{eq:three-bar-norm}
\end{align}
which is defined for all $\vvh\in\Vvhk$.
The local stability property of the bilinear forms $\msh(\cdot,\cdot)$
and $\ash(\cdot,\cdot)$ implies the equivalence relation
\begin{align}
  \msh(\dot{\vv}_{\hh},\dot{\vv}_{\hh}) + \ash(\vvh,\vvh)
  \lesssim\TNORM{\vvh(t)}{}^2\lesssim \msh(\dot{\vv}_{\hh},\dot{\vv}_{\hh}) +
  \ash(\vvh,\vvh)
  \label{eq:energy:norm:equivalence}
\end{align}
for all time-dependent virtual element functions $\vvh(t)$ with square
integrable derivative $\dot{\vv}_{\hh}(t)$.

The hidden constants in~\eqref{eq:energy:norm:equivalence} are
independent of the mesh size parameter
$\hh$~\cite{Antonietti-Manzini-Mazzieri-Mourad-Verani:2020}.
However, they may depend on the stability parameters, the physical
parameters and the polynomial degree $k$
\cite{BeiraodaVeiga-Chernov-Mascotto-Russo:2018}.
It is worth noting that the dependence on $k$ does not seem to have a
relevant impact on the optimality of the convergence rates in the
numerical experiments of Section~\ref{sec5:numerical}.
The following stability result has been proved in
\cite{Antonietti-Manzini-Mazzieri-Mourad-Verani:2020}.

\begin{theorem}\label{theorem:semi-discrete:stability}
  Let $\fv\in\LTWO\big(0,T;[\LTWO(\Omega)]^2\big)$ and let
  $\uvh\in\CS{2}\big(0,T;\Vvhk\big)$ be the solution of
  \eqref{eq:VEM:semi-discrete}.
  Then, it holds
  \begin{equation}
    \TNORM{\uvh(t)}{} \lesssim 
    \TNORM{ (\uv_{0})_I }{} + \int_0^t \NORM{ \fv(\tau) }{0,\Omega} d\tau. 
  \end{equation}
  The hidden constant in $\lesssim$ is independent of $\hh$, but may
  depend on the model parameters and approximation constants and the
  polynomial degree $k$.
\end{theorem}
We point out that in the case of $\fv$ null external force,
i.e. $\fv=\mathbf{0}$, the above bound reduces to $$\TNORM{\uvh(t)}{}
\lesssim \TNORM{ (\uv_{0})_I }{} $$ that is the virtual element
approximation is dissipative.
    
Now, we recall \cite{Antonietti-Manzini-Mazzieri-Mourad-Verani:2020}
the convergence of the semi-discrete virtual element approximation in
the energy norm \eqref{eq:three-bar-norm}.
\begin{theorem}
  Let $\uv\in\CS{2}\big(0,T;[H^{m+1}(\Omega)]^2\big)$,
  $m\in\INTG$, be the exact solution of problem~\eqref{eq:weak:A}.
  Let $\uvh\in\Vvhk$ be the solution of the semi-discrete
  problem~\eqref{eq:VEM:semi-discrete}. For
  $\fv\in\LTWO\big((0,T);\big[\HS{m-1}(\Omega)\big]^2\big)$ we have
  that
  \begin{align}
    &\sup_{0<t\leq T}\TNORM{\uv(t)-\uvh(t)}{}
    \lesssim
    \frac{\hh^\mu}{k^m}\sup_{0<t\leq T}\Big(\NORM{\dot{\uv}(t)}{m+1} + \NORM{\uv(t)}{m+1}\Big)\nonumber\\[0.5em]
    &\qquad+\int_{0}^{T}\left(
      \frac{\hh^{\mu+1}}{k^{m}}\left( \NORM{ \ddot{\uv}(\tau) }{m+1} + \NORM{ \dot{\uv}(\tau) }{m+1}\right) 
      + \frac{\hh^{\mu}}{k^{m}}\big ( \NORM{ \ddot{\uv}(\tau) }{m+1} + \NORM{ \dot{\uv}(\tau) }{m+1} \big) 
    \right)\,d\tau\nonumber\\[0.5em]
    &\qquad+\int_{0}^{T}\hh\Norm{ \big(I-\Pi^0_{k-2}\big)\fv(\tau) }{0}\,d\tau,
    \label{eq:corollary:estimate}
  \end{align}
  where $\mu=\min(k,m)$.
  The hidden constant in ``$\lesssim$`` is independent of $\hh$, but
  may depend on the model parameters and approximation constants, the
  polynomial degree $k$, and the final observation time $T$.
\end{theorem}

Finally, we state the convergence result in the $\LTWO$ norm, whose
proof is again found
in~\cite{Antonietti-Manzini-Mazzieri-Mourad-Verani:2020}.

\begin{theorem}
  \label{theorem:L2:convergence}
  Let $\uv$ be the exact solution of problem~\eqref{eq:weak:A} under
  the assumption that domain $\Omega$ is $\HTWO$-regular and
  $\uvh\in\Vvhk$ the solution of the virtual element method stated
  in~\eqref{eq:VEM:semi-discrete}.
  If
  $\uv,\dot{\uv},\ddot{\uv}\in\LTWO\big(0,T;\big[\HS{m+1}(\Omega)\cap\HONEzr(\Omega)\big]^2\big)$,
  with integer $m\geq0$, then the following estimate holds for almost
  every $t\in[0,T]$ by setting $\mu=\min(m,k)$:
  \begin{align}
    \NORM{\uv(t)-\uvh(t)}{0} 
    \lesssim& \NORM{\uvh(0)-\uv_0}{0} + \NORM{ \dot{\uv}_{\hh}(0)-\uv_{1} }{0} 
    + \frac{\hh^{\mu+1}}{ {k^{m+1}} }\Big(
    \NORM{\ddot{\uv}}{\LTWO( 0,T; [\HS{m+1}(\Omega)]^2 )} 
    \nonumber\\ &
    +
    \NORM{\dot {\uv}}{\LTWO( 0,T; [\HS{m+1}(\Omega)]^2 )} + 
    \NORM{      \uv }{\LTWO( 0,T; [\HS{m+1}(\Omega)]^2 )} 
    \Big)
    \nonumber\\ &
    + 
    \int_{0}^{T}\Norm{\big(1-\Piz{k-2}\big)\fv(\tau)}{0}^2d\tau.
    %%\hh^{\mu} \NORM{\fv}{\LTWO( 0,T; (\HS{m+1}(\Omega))^2 )},
  \end{align}
  The hidden constant in ``$\lesssim$`` is independent of $\hh$, but
  may depend on the model parameters and approximation constants
  $\varrho$, $\mu^*$, and {the polynomial degree $k$}, and the final
  observation time $T$.
\end{theorem}

\subsection{Numerical Results}
\label{sec5:numerical}

In this section, we report from
\cite{Antonietti-Manzini-Mazzieri-Mourad-Verani:2020} a set of
numerical results assessing the convergence properties of the virtual
element discretization by using a manufactured solution on three
different mesh families, each one possessing some special feature.

%% MESHES
\begin{figure}
  \centering
  \begin{tabular}{ccc}
    \begin{overpic}[scale=0.2]{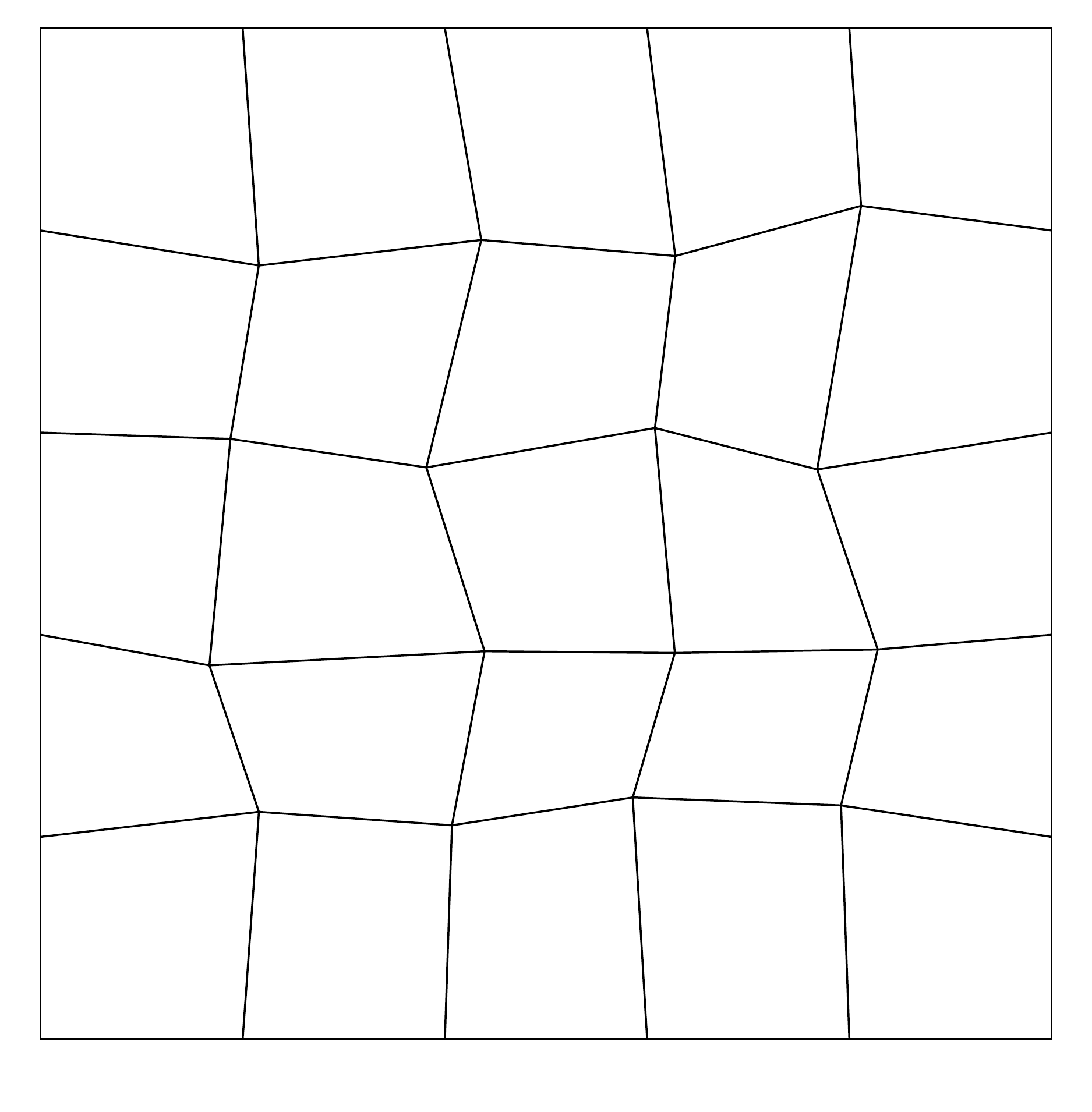}%{./figures-EPS-PDF/Convergence/Mesh-Dir/PDF-cropped/mesh2D_quads_0}
    \end{overpic} 
    &
    \begin{overpic}[scale=0.2]{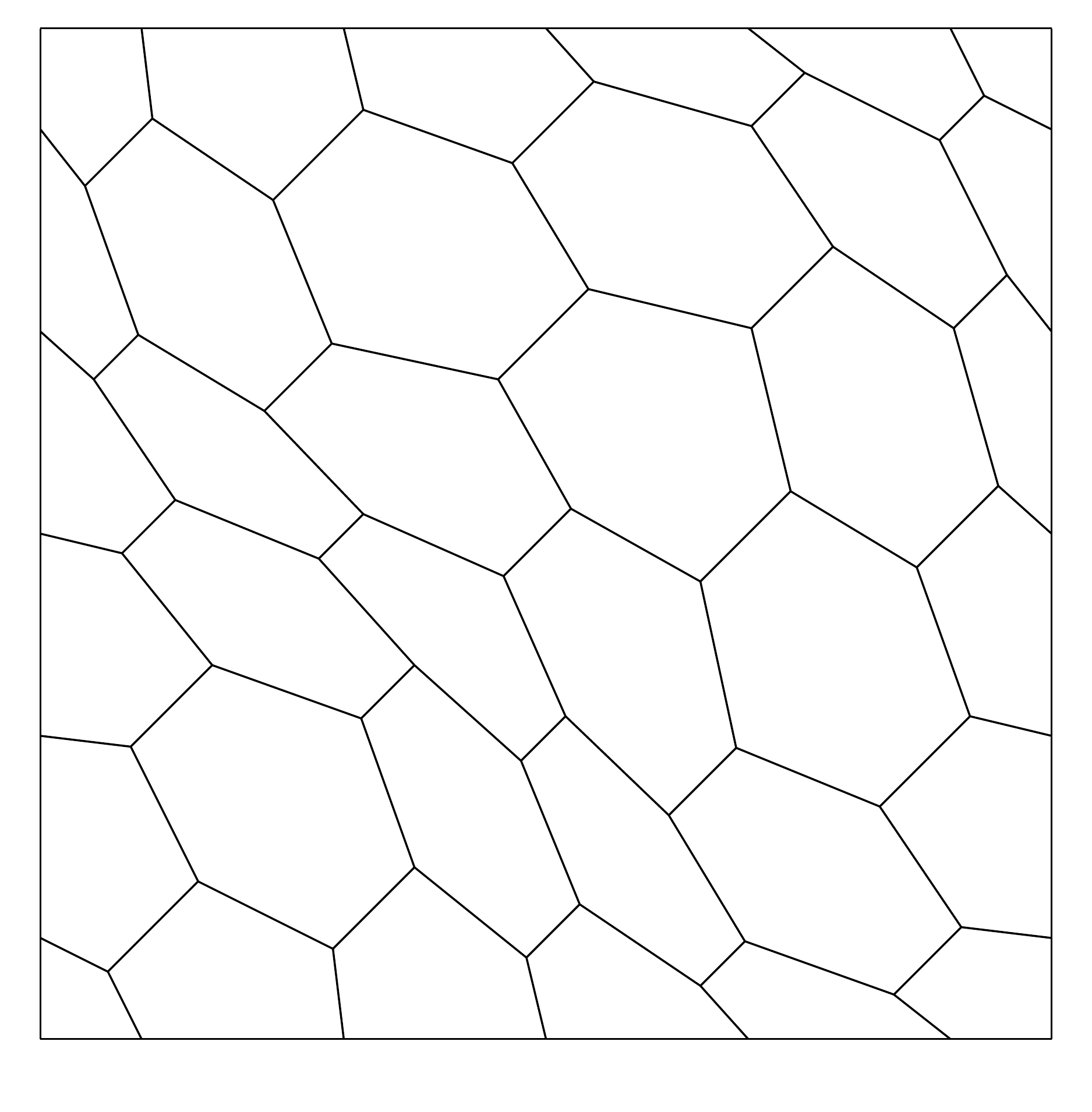}%{./figures-EPS-PDF/Convergence/Mesh-Dir/PDF-cropped/mesh2D_hexa_0}
    \end{overpic} 
    &
    \begin{overpic}[scale=0.2]{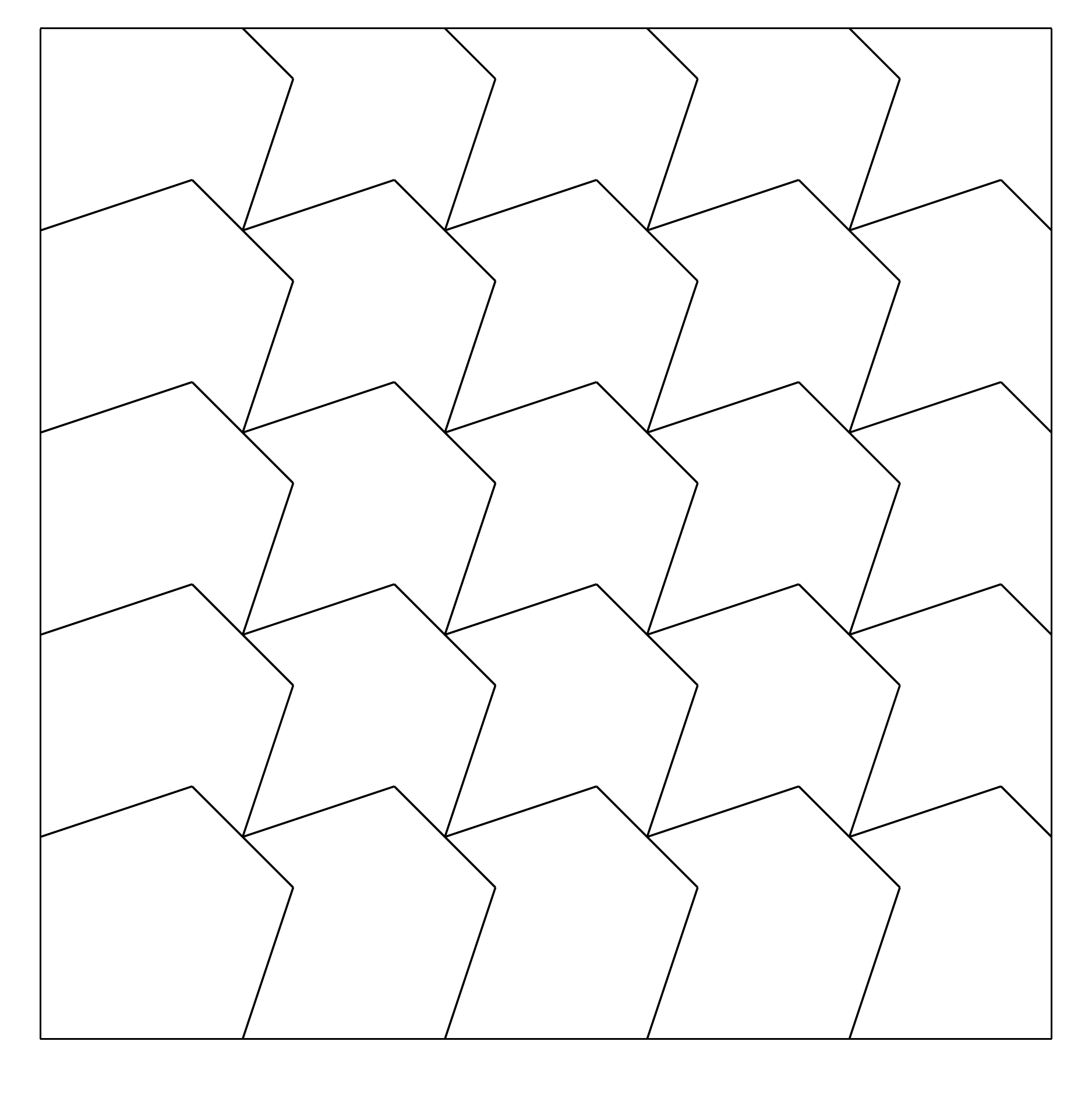}%{./figures-EPS-PDF/Convergence/Mesh-Dir/PDF-cropped/mesh2D_octa_0}
    \end{overpic}
    \\[0.25em]
    \begin{overpic}[scale=0.2]{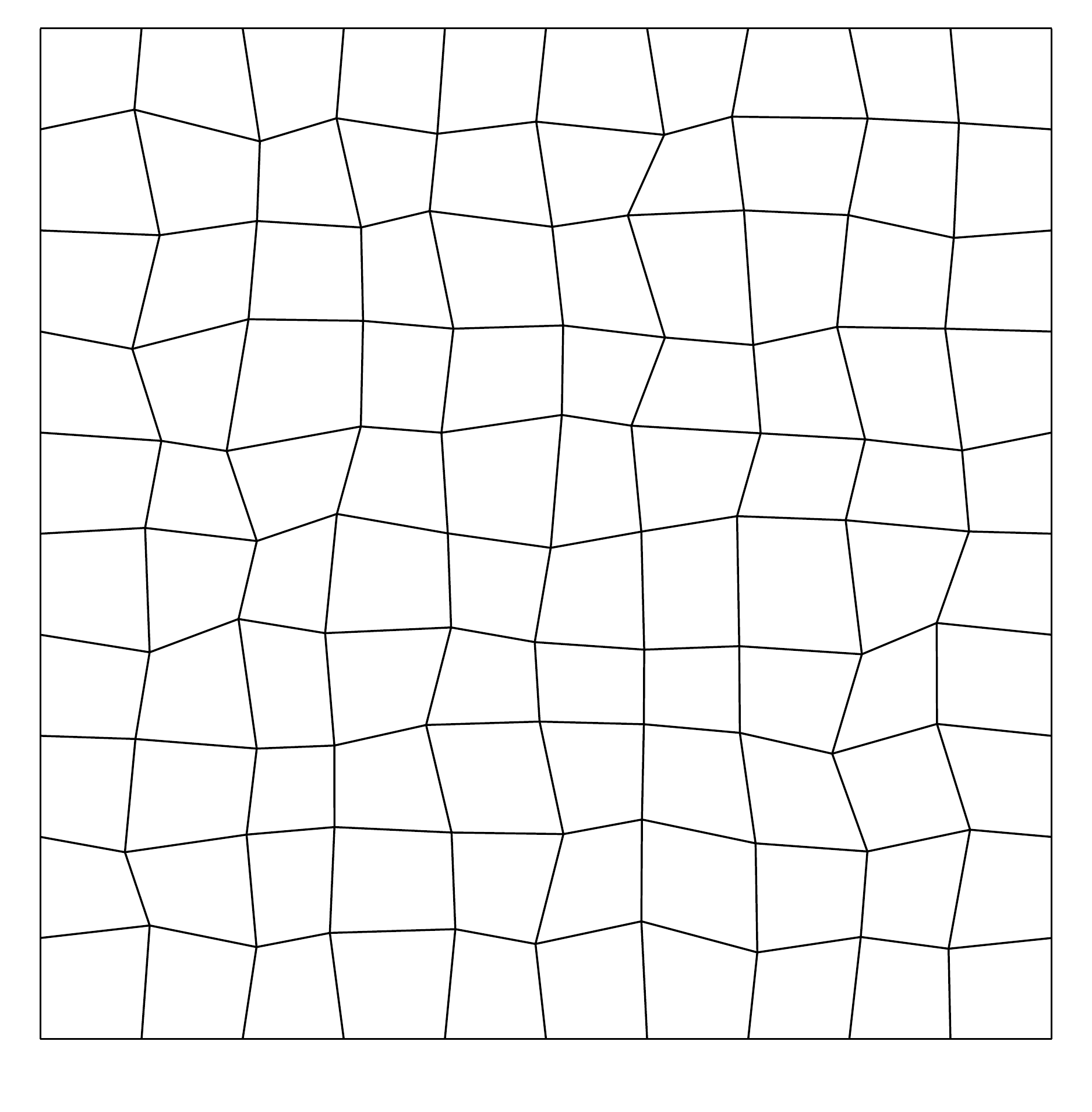}%{./figures-EPS-PDF/Convergence/Mesh-Dir/PDF-cropped/mesh2D_quads_1}
    \end{overpic} 
    &
    \begin{overpic}[scale=0.2]{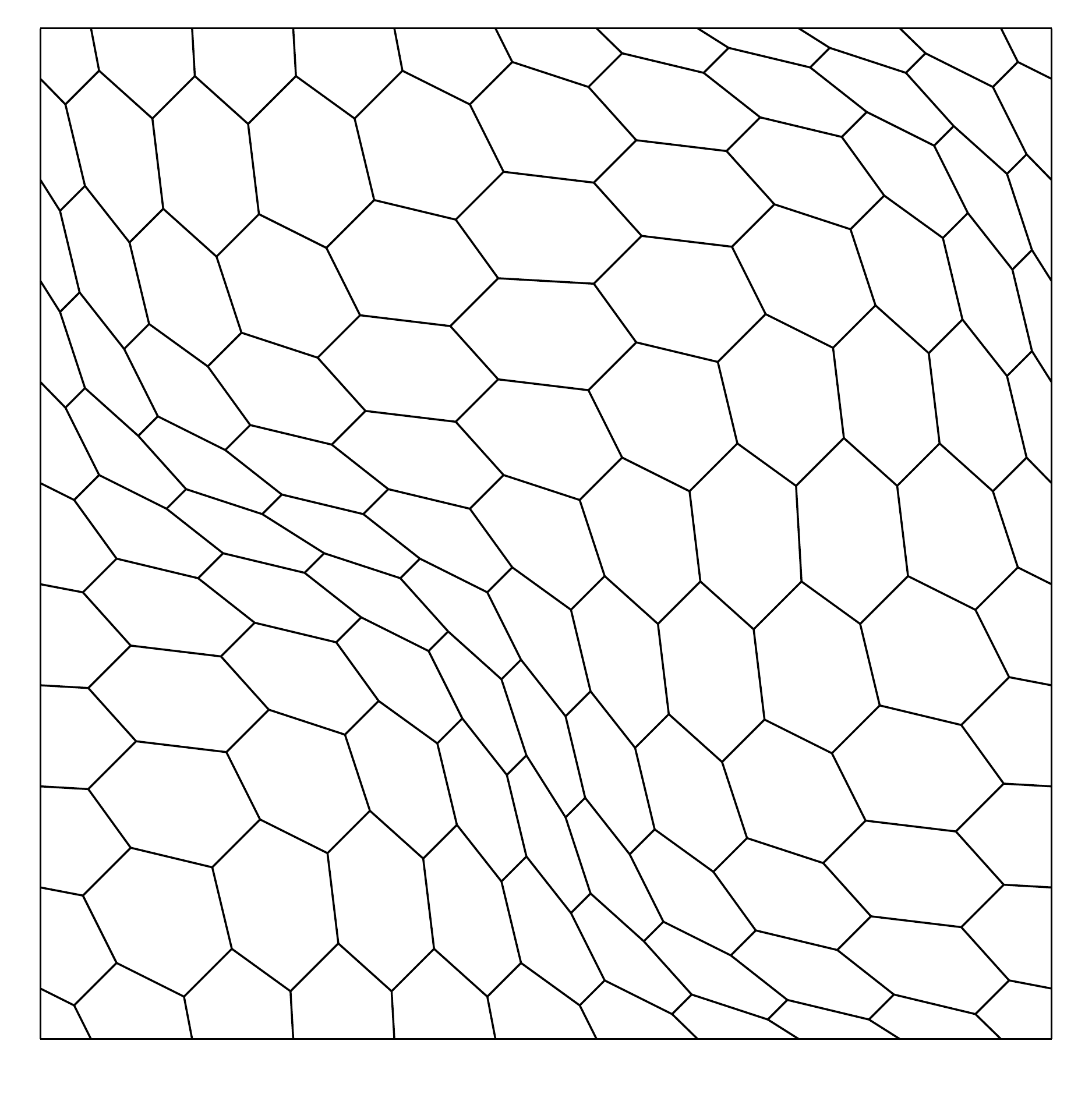}%{./figures-EPS-PDF/Convergence/Mesh-Dir/PDF-cropped/mesh2D_hexa_1}
    \end{overpic} 
    &
    \begin{overpic}[scale=0.2]{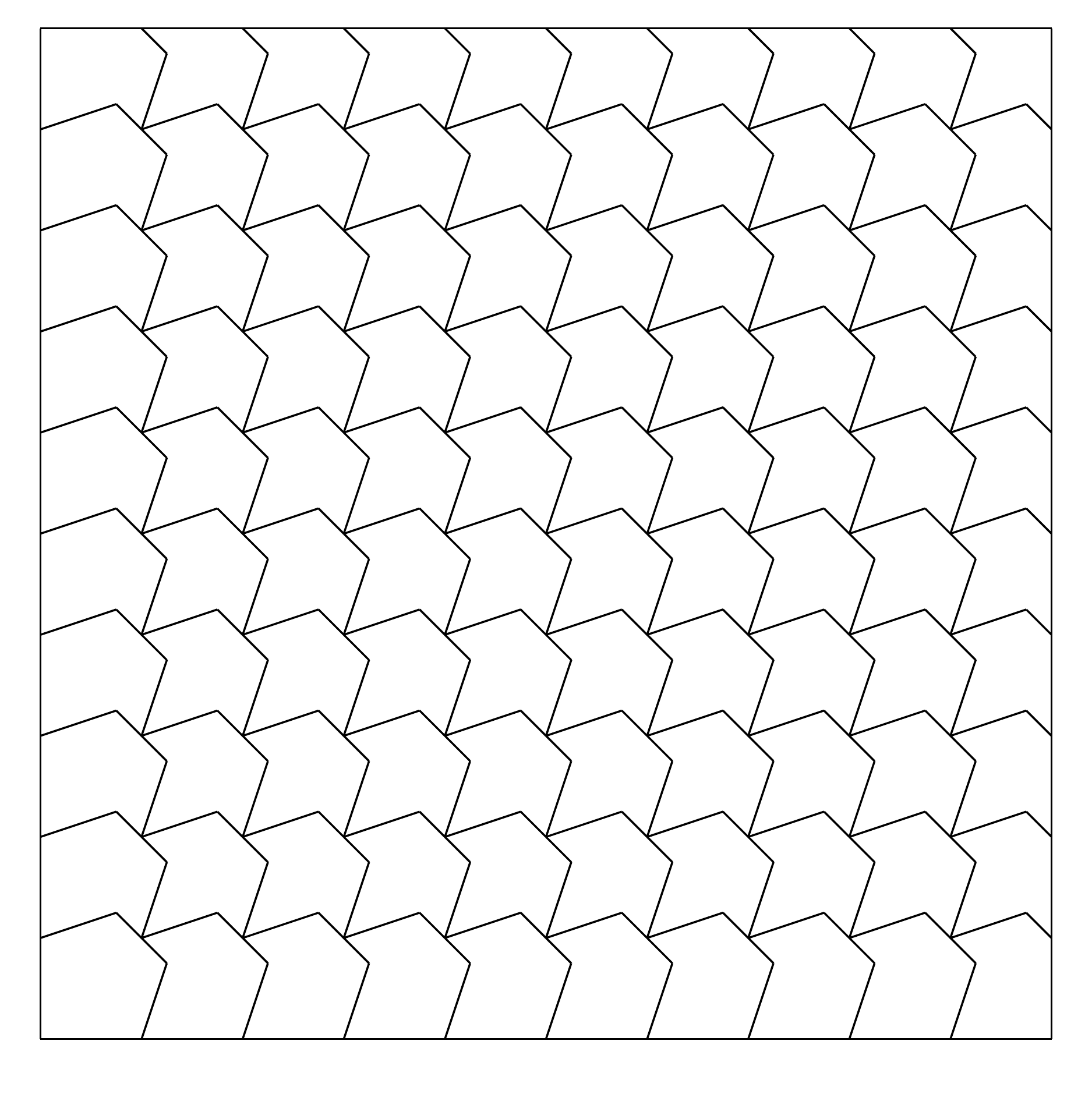}%{./figures-EPS-PDF/Convergence/Mesh-Dir/PDF-cropped/mesh2D_octa_1}
    \end{overpic}
    \\[-0.25em]
    \textit{Mesh~1} & \textit{Mesh~2} & \textit{Mesh~3}
  \end{tabular}
  \caption{ Base meshes (top row) and first refined meshes (bottom
    row) of the following mesh families from left to right: randomized
    quadrilateral mesh; mainly hexagonal mesh; nonconvex octagonal
    mesh~\cite{Antonietti-Manzini-Mazzieri-Mourad-Verani:2020}.  }
  \label{fig:Meshes}
\end{figure}

%% Regular quadrilateral mesh 
\begin{figure}
  \centering
  \begin{tabular}{cc}
    \begin{overpic}[scale=0.325]{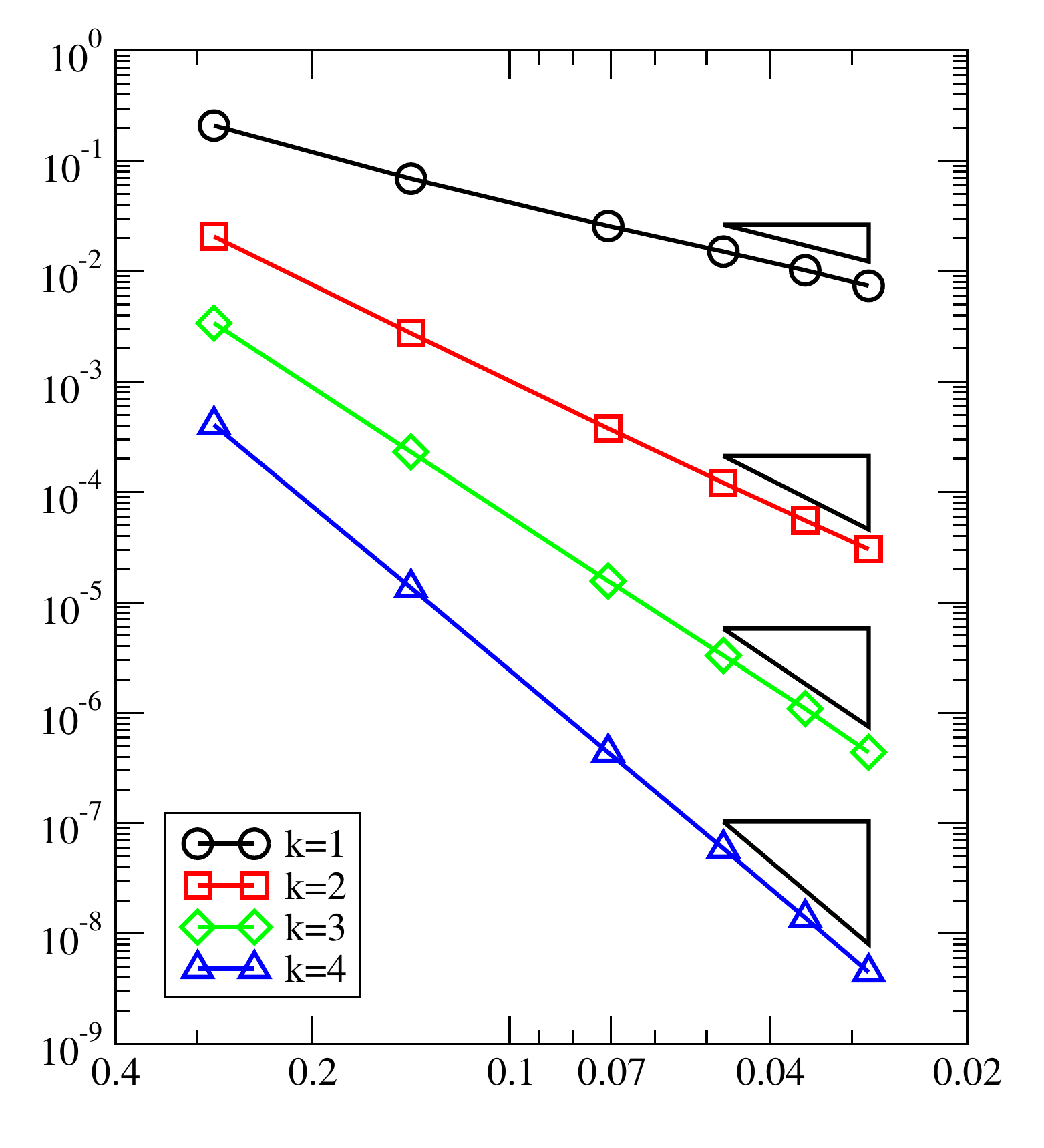}% {./figures-EPS-PDF/Convergence/PDF-cropped/quads_regular_err_L2}
      \put(-5,9){\begin{sideways}\textbf{$\mathbf{L^2}$ relative approximation error}\end{sideways}}
      \put(32,-2) {\textbf{Mesh size $\mathbf{h}$}}
      \put(68,82){\textbf{2}}
      \put(68,62){\textbf{3}}
      \put(68,47){\textbf{4}}
      \put(68,30){\textbf{5}}
    \end{overpic} 
    & \qquad
    \begin{overpic}[scale=0.325]{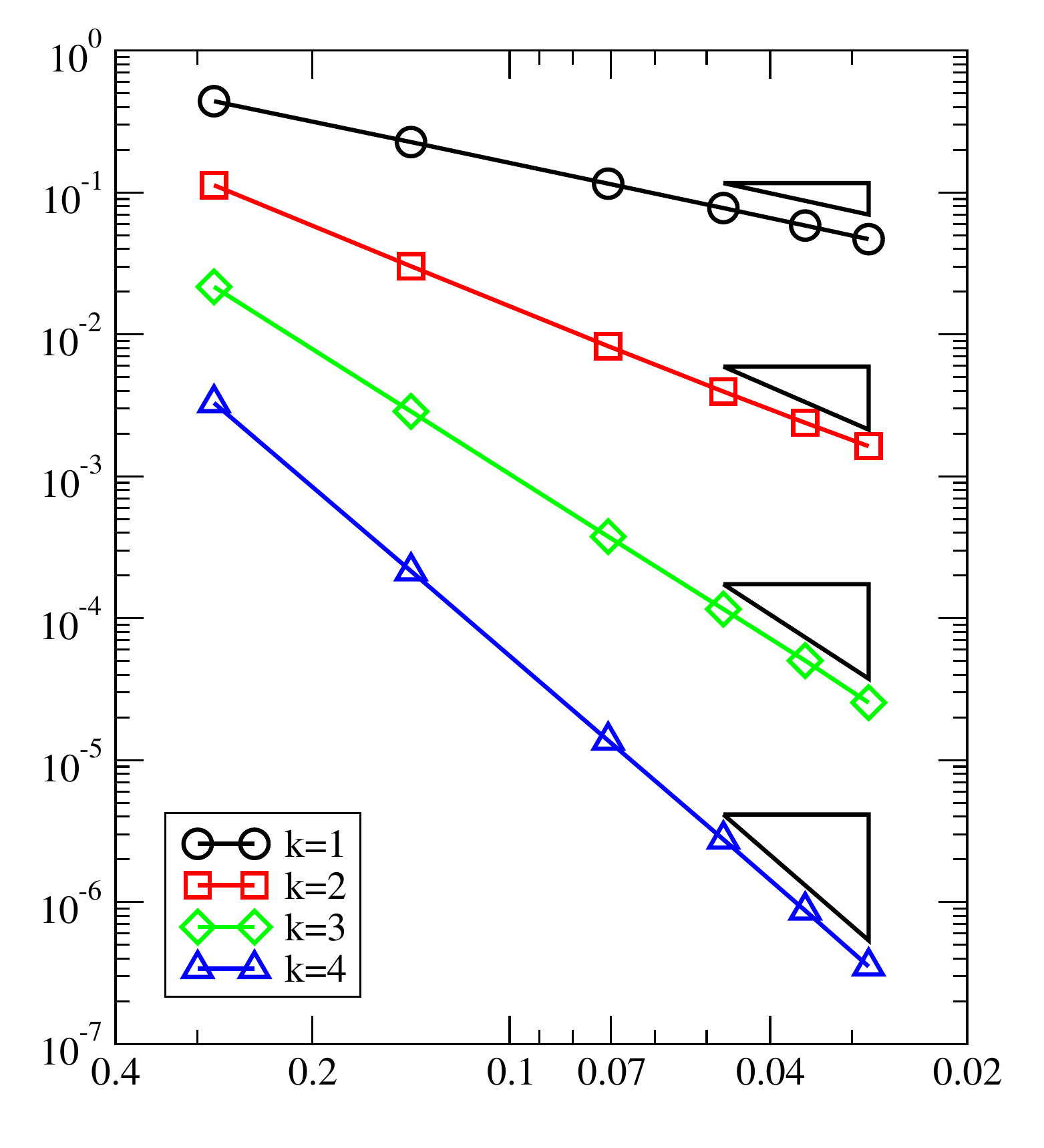}% {./figures-EPS-PDF/Convergence/PDF-cropped/quads_regular_err_H1}
      \put(-5,9){\begin{sideways}\textbf{$\mathbf{H^1}$ relative approximation error}\end{sideways}}
      \put(32,-2) {\textbf{Mesh size $\mathbf{h}$}}
      \put(68,86){\textbf{1}}
      \put(68,69){\textbf{2}}
      \put(68,50){\textbf{3}}
      \put(68,30){\textbf{4}}
    \end{overpic}
    \\[0.5em]
  \end{tabular}
    \caption{ Convergence plots for the virtual element approximation of
    Problem~\eqref{eq:pblm:strong:A}-\eqref{eq:pblm:strong:E} with
    exact solution~\eqref{eq:benchmark:solution} using
    family~\textit{Mesh~1} of randomized quadrilateral meshes.
    Error curves are computed using the $\LTWO$ norm (left panels) and
    $\HONE$ norm (right panels) and are plot versus the number of
    degrees of
    freedom~\cite{Antonietti-Manzini-Mazzieri-Mourad-Verani:2020}. }
  \label{fig:quads:rates}
\end{figure}

%% Mainly hexagonal mesh (hexa remapped)
\begin{figure}
  \centering
  \begin{tabular}{cc}
    \begin{overpic}[scale=0.325]{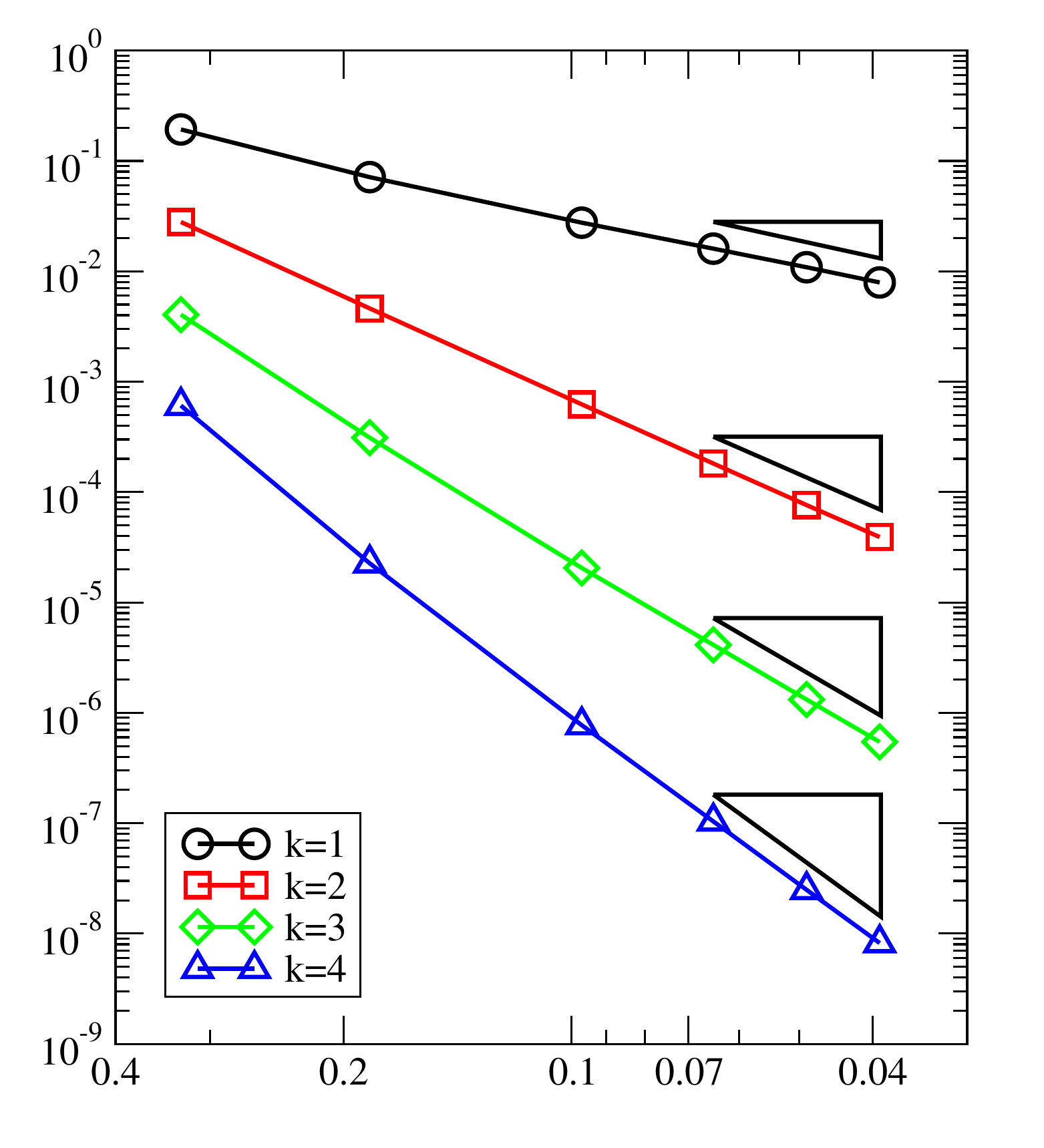}% {./figures-EPS-PDF/Convergence/PDF-cropped/hexa_remapped_err_L2}
      \put(-5,9){\begin{sideways}\textbf{$\mathbf{L^2}$ relative approximation error}\end{sideways}}
      \put(32,-2) {\textbf{Mesh size $\mathbf{h}$}}
      \put(69,82){\textbf{2}}
      \put(69,63){\textbf{3}}
      \put(69,47){\textbf{4}}
      \put(69,32){\textbf{5}}
    \end{overpic} 
    & \qquad
    \begin{overpic}[scale=0.325]{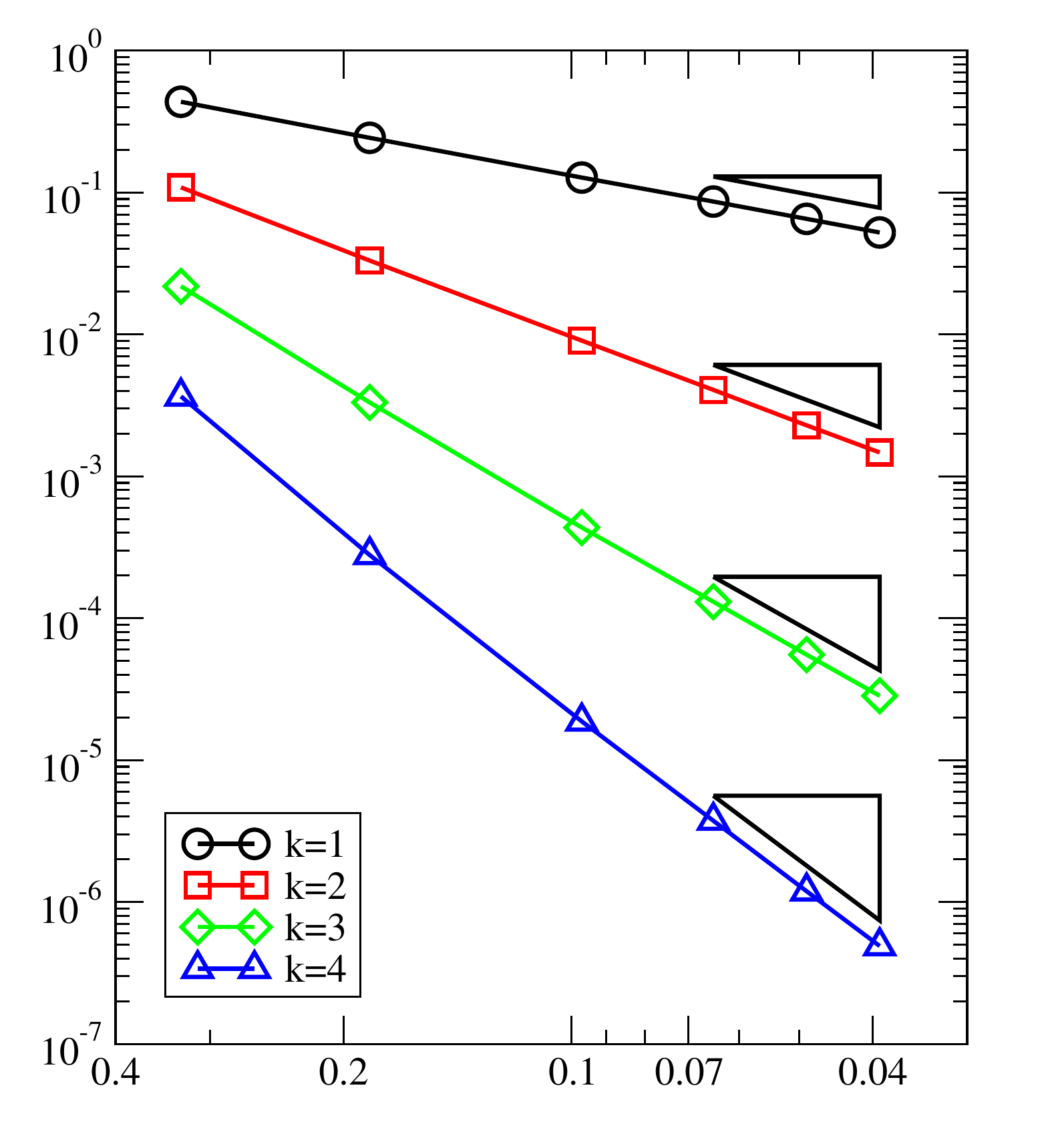}% {./figures-EPS-PDF/Convergence/PDF-cropped/hexa_remapped_err_H1}
      \put(-5,9){\begin{sideways}\textbf{$\mathbf{H^1}$ relative approximation error}\end{sideways}}
      \put(32,-2) {\textbf{Mesh size $\mathbf{h}$}}
      \put(68,86){\textbf{1}}
      \put(68,69){\textbf{2}}
      \put(68,51){\textbf{3}}
      \put(68,32){\textbf{4}}
    \end{overpic}
    \\[0.5em]
  \end{tabular}
    \caption{ Convergence plots for the virtual element approximation of
    Problem~\eqref{eq:pblm:strong:A}-\eqref{eq:pblm:strong:E} with
    exact solution~\eqref{eq:benchmark:solution} using
    family~\textit{Mesh~2} of mainly hexagonal meshes.
    Error curves are computed using the $\LTWO$ norm (left panels) and
    $\HONE$ norm (right panels) and are plot versus the number of
    degrees of
    freedom~\cite{Antonietti-Manzini-Mazzieri-Mourad-Verani:2020}. }
  \label{fig:hexa:rates}
\end{figure}

%% Nonconvex octagon mesh
\begin{figure}
  \centering
  \begin{tabular}{cc}
    \begin{overpic}[scale=0.325]{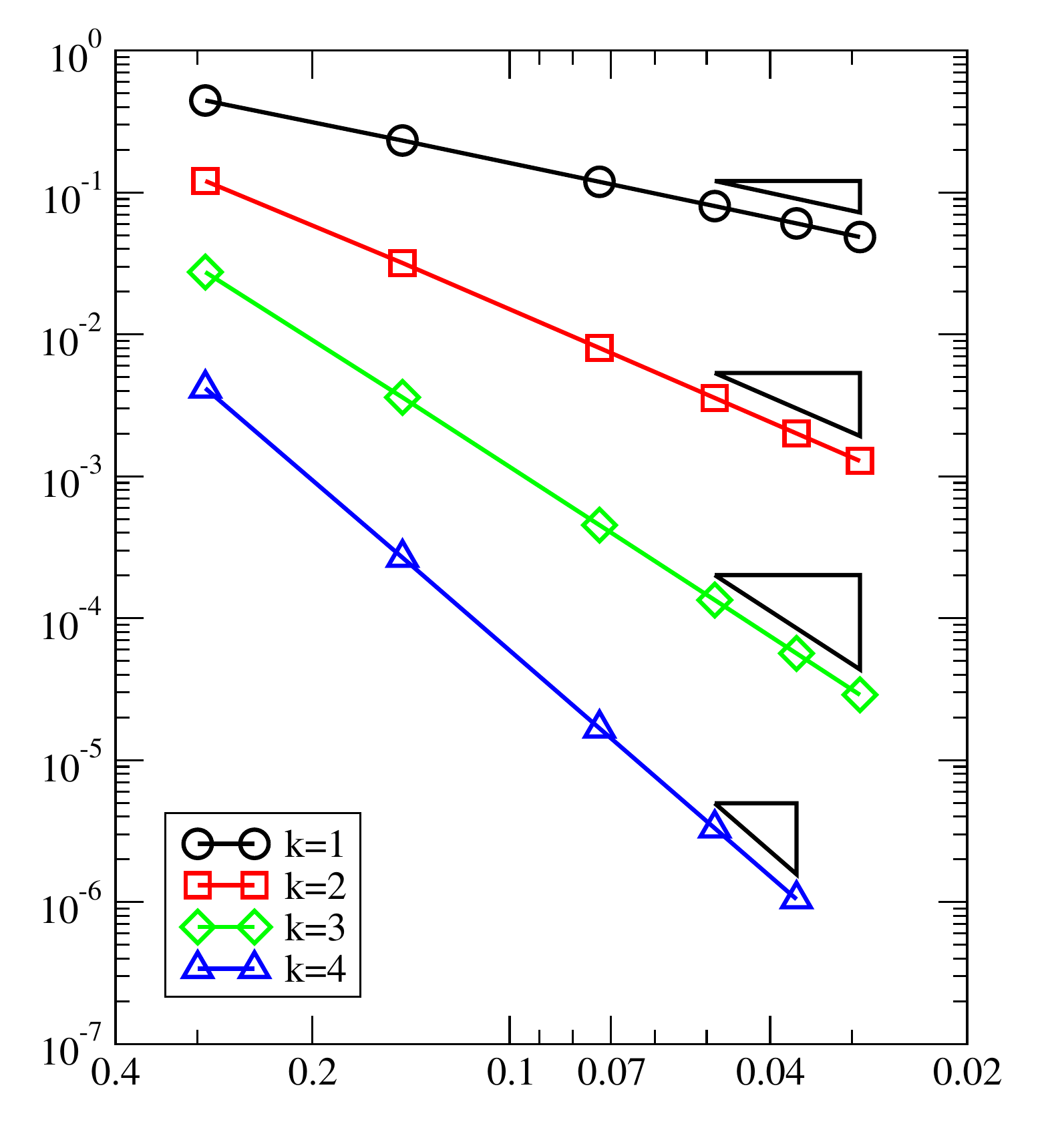}% {./figures-EPS-PDF/Convergence/PDF-cropped/octa_nonconvex_err_H1}
      \put(-5,9){\begin{sideways}\textbf{$\mathbf{L^2}$ relative approximation error}\end{sideways}}
      \put(32,-2) {\textbf{Mesh size $\mathbf{h}$}}
      \put(69,86){\textbf{2}}
      \put(69,69){\textbf{3}}
      \put(68,51){\textbf{4}}
      \put(65,31){\textbf{5}}
    \end{overpic} 
    & \qquad
    \begin{overpic}[scale=0.325]{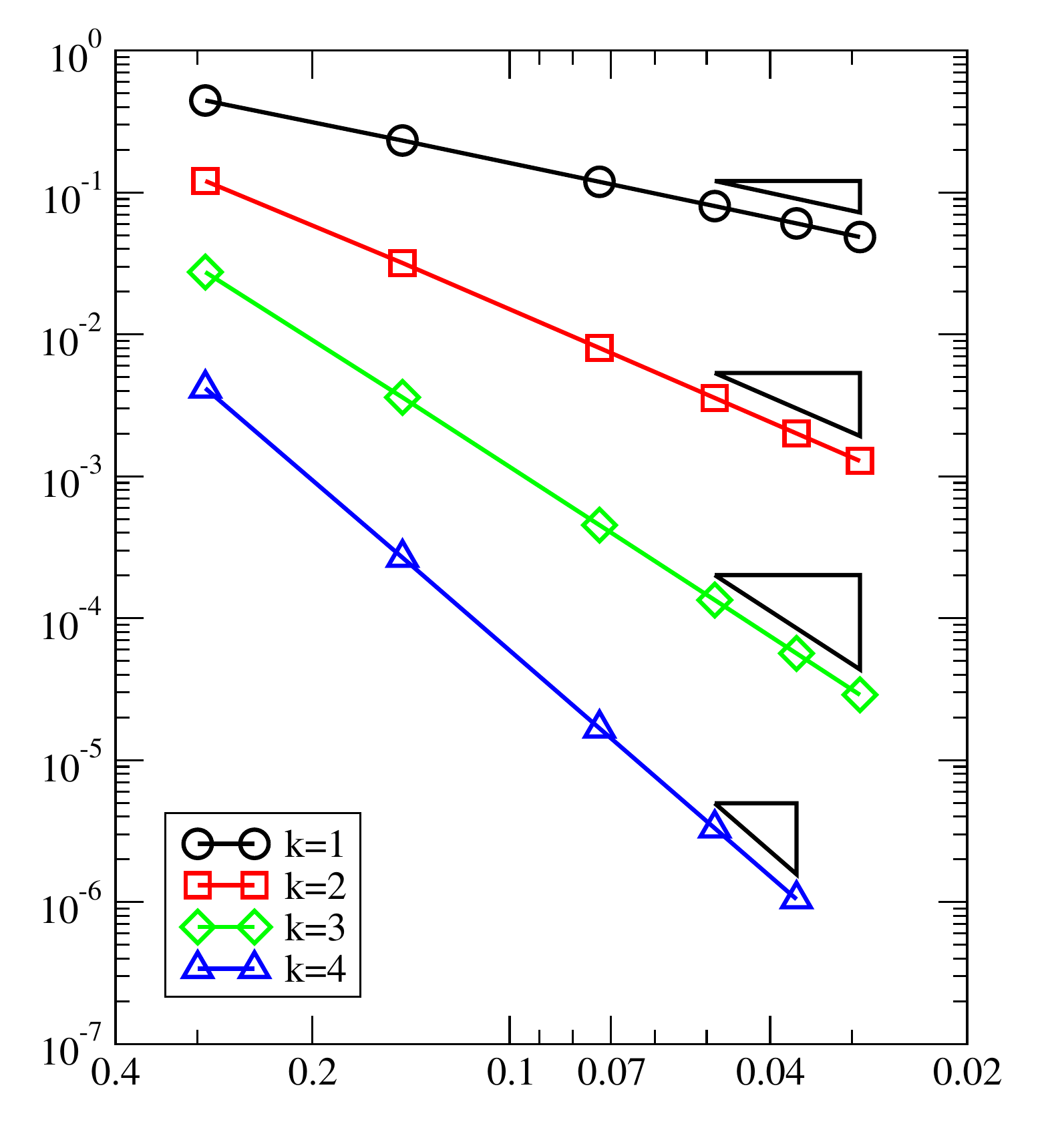}% {./figures-EPS-PDF/Convergence/PDF-cropped/octa_nonconvex_err_H1}
      \put(-5,9){\begin{sideways}\textbf{$\mathbf{H^1}$ relative approximation error}\end{sideways}}
      \put(32,-2) {\textbf{Mesh size $\mathbf{h}$}}
      \put(68,86){\textbf{1}}
      \put(68,69){\textbf{2}}
      \put(68,51){\textbf{3}}
      \put(64,31){\textbf{4}}
    \end{overpic}
    \\[0.5em]
  \end{tabular}
    \caption{ Convergence plots for the virtual element approximation of
    Problem~\eqref{eq:pblm:strong:A}-\eqref{eq:pblm:strong:E} with
    exact solution~\eqref{eq:benchmark:solution} using
    family~\textit{Mesh~3} of nonconvex octagonal meshes.
    Error curves are computed using the $\LTWO$ norm (left panels) and
    $\HONE$ norm (right panels) and are plot versus the number of
    degrees of
    freedom~\cite{Antonietti-Manzini-Mazzieri-Mourad-Verani:2020}. }
  \label{fig:octa:rates}
\end{figure}

%% k-refinement: Randomized quadrilateral meshes
\begin{figure}
  \centering
  \begin{tabular}{cc}
    \begin{overpic}[scale=0.325]{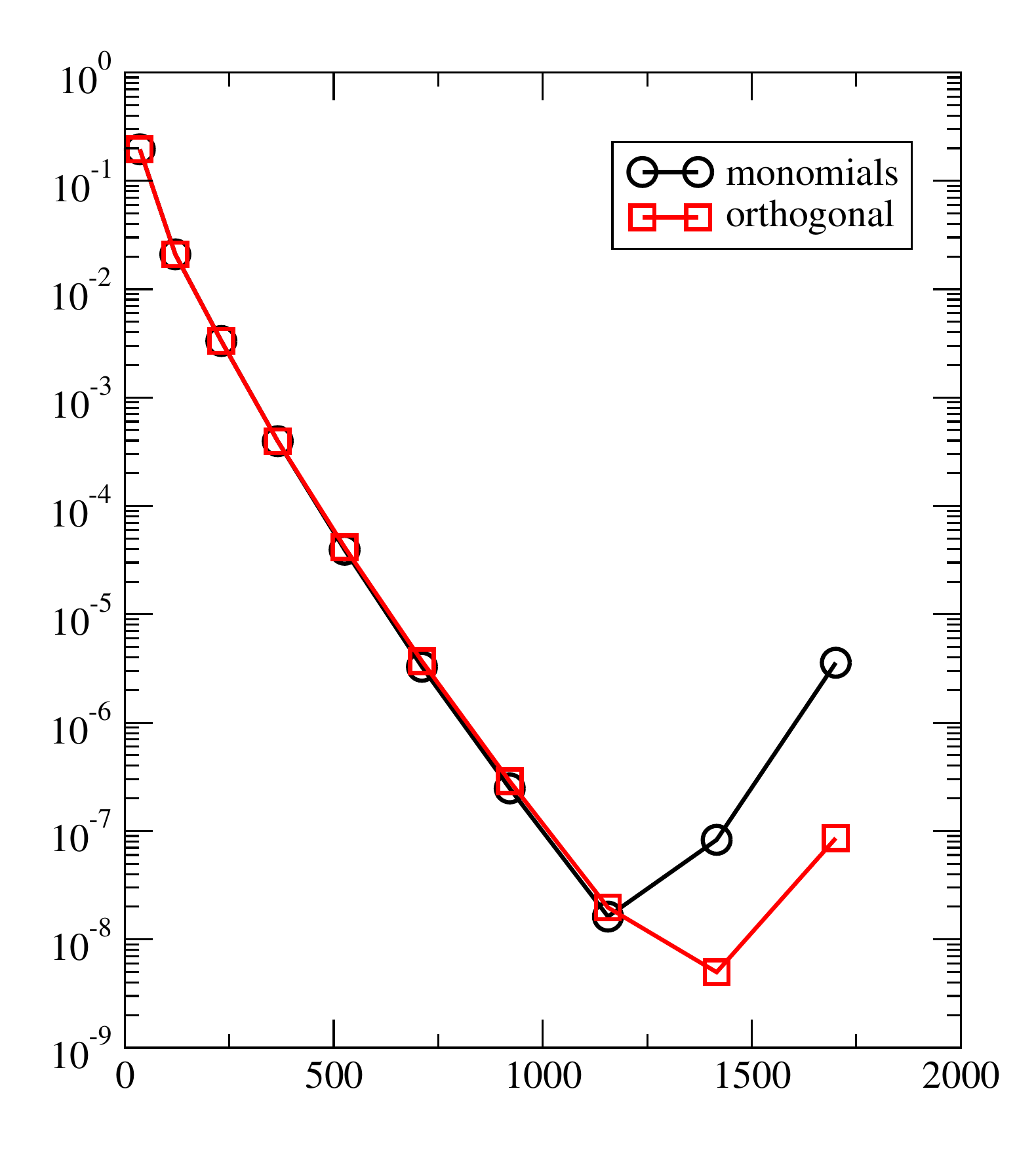}% {./figures-EPS-PDF/Convergence/PDF-cropped/pref_L2}
      \put(-5,11){\begin{sideways}\textbf{$\mathbf{L^2}$ relative approximation error}\end{sideways}}
      \put(20,-2) {\textbf{\#degrees of freedom}}
    \end{overpic} 
    & \qquad
    \begin{overpic}[scale=0.325]{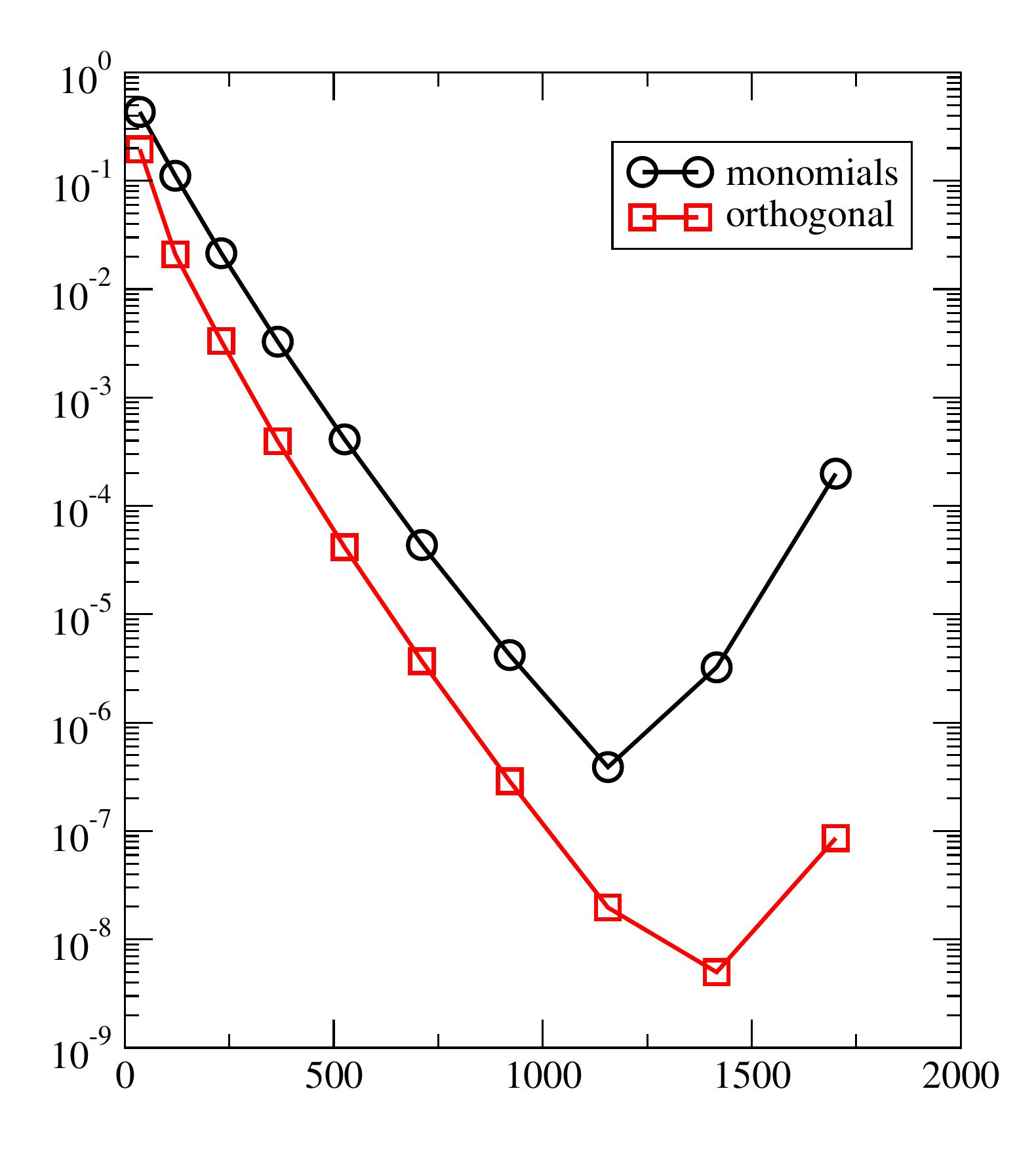}% {./figures-EPS-PDF/Convergence/PDF-cropped/pref_H1}
      \put(-5,11){\begin{sideways}\textbf{$\mathbf{H^1}$ relative approximation error}\end{sideways}}
      \put(20,-2) {\textbf{\#degrees of freedom}}
    \end{overpic}
    \\[0.5em]
  \end{tabular}
  \caption{ Convergence plots for the virtual element approximation of
    Problem~\eqref{eq:pblm:strong:A}-\eqref{eq:pblm:strong:E} with
    exact solution~\eqref{eq:benchmark:solution} using
    family~\textit{Mesh~1} of randomized quadrilateral meshes.
    Error curves are computed using k-refinement the $\LTWO$ norm
    (left panel) and $\HONE$ norm (right panel) and are plot versus
    the number of degrees of freedom by performing a refinement of
    type ``p'' on a $5\times 5$ mesh.
    Each plot shows the two convergence curves that are obtained using
    monomials (circles) and orthogonalized polynomials
    (squares)~\cite{Antonietti-Manzini-Mazzieri-Mourad-Verani:2020}. }
  \label{fig:expo:rate}
\end{figure}

In particular, we let $\Omega=(0,1)^2$ for $t\in[0,T]$, $T=1$, and
consider initial condition $\uv_{0}$, boundary condition $\gv$ and
forcing term $\fv$ determined from the exact solution:
\begin{align}
  \uv(x,y,t) = 
  \cos\left(\frac{2\pi\,t}{T}\right)
  \left(
    \begin{array}{c}
      \sin^2(\pi\xs)\sin(2\pi\ys)\\[0.5em]
      \sin(2\pi\xs)\sin^2(\pi\ys)
    \end{array}
  \right).
  \label{eq:benchmark:solution}
\end{align}
To this end, we consider three different mesh partitionings, denoted
by:
\begin{itemize}
\item \textit{Mesh~1}, randomized quadrilateral mesh;
\item \textit{Mesh~2}, mainly hexagonal mesh with continuously
  distorted cells;
\item \textit{Mesh~3}, nonconvex octagonal mesh.
\end{itemize}
The base mesh and the first refined mesh of each mesh sequence are
shown in Figure~\ref{fig:Meshes}.

The discretization in time is given by applying the leap-frog method
with $\Delta t=10^{-4}$ and carried out for $10^{4}$ time cycles in
order to reach time $T=1$.

For these calculations, we used the VEM approximation based on the
conforming space $\Vhk$ with {$k=1,2,3,4$} and the convergence curves
for the three mesh sequences above are reported in
Figures~\ref{fig:quads:rates},~\ref{fig:hexa:rates}
and~\ref{fig:octa:rates}.
The expected rate of convergence is shown in each panel by the
triangle closed to the error curve and indicated by an explicit
label.
The results are in agreement with the theoretical estimates.
To conclude, Figure~\ref{fig:expo:rate} shows the semilog error curves
obtained through a``p''-type refinement calculation for the previous
benchmark, i.e.\ for a fixed $5\times 5$ mesh of type $I$ the order of
the virtual element space is increased from $k=1$ to $k=10$.
We employ two different implementations, namely in the first case the
space of polynomials of degree $k$ is generated by the standard scaled
monomials, while in the second one we consider an orthogonal
polynomial basis.
The behavior of the VEM when using nonorthogonal and orthogonal
polynomial basis shown in Figure~\ref{fig:expo:rate} is in accordance
with the literature, see, e.g.,
\cite{Berrone-Borio:2017,Mascotto:2018}.

%% section-5
%% final remarks

%% ---------------
%% Acknowledgments
%% ---------------

\section*{Acknowledgement}
PFA and MV acknowledge the financial support of PRIN research grant
number 201744KLJL ``\emph{Virtual Element Methods: Analysis and
Applications}'' funded by MIUR.  PFA, IM, and MV, and SS acknowledges
the financial support of INdAM-GNCS.
GM acknowledges the financial support of the ERC Project CHANGE, which
has received funding from the European Research Council under the
European Union’s Horizon 2020 research and innovation program (grant
agreement no.~694515).

%% ------------
%% BIBLIOGRAPHY
%% ------------

%% \bibliographystyle{plain}
%% \bibliography{chapter-MOX}

\begin{thebibliography}{10}

\bibitem{Adams-Fournier:2003}
R.~A. Adams and J.~J.~F. Fournier.
\newblock {\em Sobolev spaces}.
\newblock Pure and Applied Mathematics. Academic Press, 2 edition, 2003.

\bibitem{Ahmad-Alsaedi-Brezzi-Marini-Russo:2013}
B.~Ahmad, A.~Alsaedi, F.~Brezzi, L.~D. Marini, and A.~Russo.
\newblock Equivalent projectors for virtual element methods.
\newblock {\em Comput. Math. Appl.}, 66(3):376--391, 2013.

\bibitem{Aldakheel-Hudobivnik-Hussein-Wriggers:2018}
F.~Aldakheel, B.~Hudobivnik, A.~Hussein, and P.~Wriggers.
\newblock Phase-field modeling of brittle fracture using an efficient virtual
  element scheme.
\newblock {\em Comput. Methods Appl. Mech. Engrg.}, 341:443--466, 2018.

\bibitem{Antonietti-AyusodeDios-Mazzieri-Quarteroni:2016}
P.~F. Antonietti, B.~{Ayuso~de~Dios}, I.~Mazzieri, and A.~Quarteroni.
\newblock Stability analysis of discontinuous {G}alerkin approximations to the
  elastodynamics problem.
\newblock {\em J. Sci. Comput.}, 68(1):143--170, 2016.

\bibitem{Antonietti-BeiraodaVeiga-Scacchi-Verani:2016}
P.~F. Antonietti, L.~{Beir\~{a}o~da~Veiga}, S.~Scacchi, and M.~Verani.
\newblock A {$C^1$} virtual element method for the {C}ahn-{H}illiard equation
  with polygonal meshes.
\newblock {\em SIAM J. Numer. Anal.}, 54(1):34--56, 2016.

\bibitem{Antonietti-Bonaldi-Mazzieri:2018}
P.~F. Antonietti, F.~Bonaldi, and I.~Mazzieri.
\newblock A high-order discontinuous {G}alerkin approach to the elasto-acoustic
  problem.
\newblock {\em Comput. Methods Appl. Mech. Engrg.}, 358:112634, 29, 2020.

\bibitem{Antonietti-Manzini-Mazzieri-Mourad-Verani:2020}
P.~F. Antonietti, G.~Manzini, I.~Mazzieri, H.~Mourad, and M~Verani.
\newblock The virtual element method for linear elastodynamics models.
  {C}onvergence, stability and dissipation-dispersion analysis.
\newblock arXiv:1912.07122, 2020.
\newblock accepted on {\em International Journal for Numerical Methods in
  Engineering}.

\bibitem{Antonietti-Manzini-Verani:2018}
P.~F. Antonietti, G.~Manzini, and M.~Verani.
\newblock The fully nonconforming virtual element method for biharmonic
  problems.
\newblock {\em Math. Models Methods Appl. Sci.}, 28(2):387--407, 2018.

\bibitem{Antonietti-Manzini-Verani:2019}
P.~F. Antonietti, G.~Manzini, and M.~Verani.
\newblock The conforming virtual element method for polyharmonic problems.
\newblock {\em Comput. Math. Appl.}, 79(7):2021--2034, 2020.

\bibitem{Antonietti-Mazzieri:2018}
P.~F. Antonietti and I.~Mazzieri.
\newblock High-order discontinuous {G}alerkin methods for the elastodynamics
  equation on polygonal and polyhedral meshes.
\newblock {\em Comput. Methods Appl. Mech. Engrg.}, 342:414--437, 2018.

\bibitem{Antonietti-Mazzieri-Quarteroni-Rapetti:2012}
P.~F. Antonietti, I.~Mazzieri, A.~Quarteroni, and F.~Rapetti.
\newblock Non-conforming high order approximations of the elastodynamics
  equation.
\newblock {\em Comput. Methods Appl. Mech. Engrg.}, 209/212:212--238, 2012.

\bibitem{AyusodeDios-Lipnikov-Manzini:2016}
B.~{Ayuso~de~Dios}, K.~Lipnikov, and G.~Manzini.
\newblock The non-conforming virtual element method.
\newblock {\em ESAIM Math. Model. Numer.}, 50(3):879--904, 2016.

\bibitem{Barrett-Langdon-Nurnberg:2004}
J.~W. Barrett, S.~Langdon, and R.~N\"{u}rnberg.
\newblock Finite element approximation of a sixth order nonlinear degenerate
  parabolic equation.
\newblock {\em Numer. Math.}, 96(3):401--434, 2004.

\bibitem{BeiraodaVeiga-Brezzi-Cangiani-Manzini-Marini-Russo:2013}
L.~{Beir\~{a}o~da~Veiga}, F.~Brezzi, A.~Cangiani, G.~Manzini, L.~D. Marini, and
  A.~Russo.
\newblock Basic principles of virtual element methods.
\newblock {\em Math. Models Methods Appl. Sci.}, 23(1):199--214, 2013.

\bibitem{BeiraodaVeiga-Brezzi-Marini:2013}
L.~{Beir\~{a}o~da~Veiga}, F.~Brezzi, and L.~D. Marini.
\newblock Virtual elements for linear elasticity problems.
\newblock {\em SIAM J. Numer. Anal.}, 51(2):794--812, 2013.

\bibitem{BeiraodaVeiga-Brezzi-Marini-Russo:2016c}
L.~{Beir\~ao~da~Veiga}, F.~Brezzi, L.~D. Marini, and A.~Russo.
\newblock Mixed virtual element methods for general second order elliptic
  problems on polygonal meshes.
\newblock {\em ESAIM: Mathematical Modelling and Numerical Analysis},
  50(3):727--747, 2016.

\bibitem{BeiraodaVeiga-Brezzi-Marini-Russo:2016b}
L.~{Beir\~ao~da~Veiga}, F.~Brezzi, L.~D. Marini, and A.~Russo.
\newblock Virtual element methods for general second order elliptic problems on
  polygonal meshes.
\newblock {\em Math. Models Methods Appl. Sci.}, 26(4):729--750, 2016.

\bibitem{BeiraodaVeiga-Chernov-Mascotto-Russo:2016}
L.~{Beir\~{a}o~da~Veiga}, A.~Chernov, L.~Mascotto, and A.~Russo.
\newblock Basic principles of {$hp$} virtual elements on quasiuniform meshes.
\newblock {\em Math. Models Methods Appl. Sci.}, 26(8):1567--1598, 2016.

\bibitem{BeiraodaVeiga-Chernov-Mascotto-Russo:2018}
L.~{Beir\~{a}o~da~Veiga}, A.~Chernov, L.~Mascotto, and A.~Russo.
\newblock Exponential convergence of the {$hp$} virtual element method in
  presence of corner singularities.
\newblock {\em Numer. Math.}, 138(3):581--613, 2018.

\bibitem{BeiraodaVeiga-Dassi-Russo:2020}
L.~{Beir\~{a}o~da~Veiga}, F.~Dassi, and A.~Russo.
\newblock A {$C^1$} virtual element method on polyhedral meshes.
\newblock {\em Comput. Math. Appl.}, 79(7):1936--1955, 2020.

\bibitem{BeiraodaVeiga-Lipnikov-Manzini:2011}
L.~{Beir\~ao~da~Veiga}, K.~Lipnikov, and G.~Manzini.
\newblock Arbitrary order nodal mimetic discretizations of elliptic problems on
  polygonal meshes.
\newblock {\em SIAM J. Numer. Anal.}, 49(5):1737--1760, 2011.

\bibitem{BeiraodaVeiga-Lipnikov-Manzini:2014}
L.~{Beir\~ao~da~Veiga}, K.~Lipnikov, and G.~Manzini.
\newblock {\em The Mimetic Finite Difference Method}, volume~11 of {\em MS\&A.
  Modeling, Simulations and Applications}.
\newblock Springer, {I} edition, 2014.

\bibitem{BeiraodaVeiga-Lovadina-Mora:2015}
L.~{Beir\~{a}o~da~Veiga}, C.~Lovadina, and D.~Mora.
\newblock A virtual element method for elastic and inelastic problems on
  polytope meshes.
\newblock {\em Comput. Methods Appl. Mech. Engrg.}, 295:327--346, 2015.

\bibitem{BeiraodaVeiga-Manzini:2014}
L.~{Beir\~ao~da~Veiga} and G.~Manzini.
\newblock A virtual element method with arbitrary regularity.
\newblock {\em IMA J. Numer. Anal.,}, 34(2):782--799, 2014.
\newblock DOI: 10.1093/imanum/drt018, (first published online 2013).

\bibitem{BeiraodaVeiga-Manzini:2015}
L.~{Beir\~{a}o~da~Veiga} and G.~Manzini.
\newblock Residual {\it a posteriori} error estimation for the virtual element
  method for elliptic problems.
\newblock {\em ESAIM Math. Model. Numer. Anal.}, 49(2):577--599, 2015.

\bibitem{Benvenuti-Chiozzi-Manzini-Sukumar:2019:CMAME:journal}
E.~Benvenuti, A.~Chiozzi, G.~Manzini, and N.~Sukumar.
\newblock Extended virtual element method for the {L}aplace problem with
  singularities and discontinuities.
\newblock {\em Comput. Methods Appl. Mech. Engrg.}, 356:571 -- 597, 2019.

\bibitem{Bernardi-Dauge-Maday:2007}
C.~Bernardi, M.~Dauge, and Y.~Maday.
\newblock Polynomials in the {S}obolev world.
\newblock Technical report, HAL, 2007.
\newblock hal-00153795,.

\bibitem{Berrone-Borio:2017}
S.~Berrone and A.~Borio.
\newblock Orthogonal polynomials in badly shaped polygonal elements for the
  virtual element method.
\newblock {\em Finite Elem. Anal. Des.}, 129:14--31, 2017.

\bibitem{Berrone-Borio-Manzini:2018:CMAME:journal}
S.~Berrone, A.~Borio, and Manzini.
\newblock {SUPG} stabilization for the nonconforming virtual element method for
  advection–diffusion–reaction equations.
\newblock {\em Comput. Methods Appl. Mech. Engrg.}, 340:500--529, 2018.

\bibitem{Berrone-Pieraccini-Scialo-Vicini:2015}
S.~Berrone, S.~Pieraccini, S.~Scial{\`o}, and F.~Vicini.
\newblock A parallel solver for large scale {DFN} flow simulations.
\newblock {\em SIAM J. Sci. Comput.}, 37(3):C285--C306, 2015.

\bibitem{Borden-Hughes-Landis-Verhoosel:2014}
M.~J. Borden, T.~J.~R. Hughes, C.~M. Landis, and C.~V. Verhoosel.
\newblock A higher-order phase-field model for brittle fracture: formulation
  and analysis within the isogeometric analysis framework.
\newblock {\em Comput. Methods Appl. Mech. Engrg.}, 273:100--118, 2014.

\bibitem{Bramble-Falk:1985}
J.~H. Bramble and R.~S. Falk.
\newblock A mixed-{L}agrange multiplier finite element method for the
  polyharmonic equation.
\newblock {\em RAIRO Mod\'{e}l. Math. Anal. Num\'{e}r.}, 19(4):519--557, 1985.

\bibitem{Brenner-Scott:2008}
S.~C. Brenner and R.~Scott.
\newblock {\em The mathematical theory of finite element methods}, volume~15.
\newblock Springer Science \& Business Media, 2008.

\bibitem{Brezzi-Buffa-Lipnikov:2009}
F.~Brezzi, A.~Buffa, and K.~Lipnikov.
\newblock Mimetic finite differences for elliptic problems.
\newblock {\em M2AN Math. Model. Numer. Anal.}, 43:277--295, 2009.

\bibitem{Brezzi-Falk-Marini:2014}
F.~Brezzi, R.~S. Falk, and L.~D. Marini.
\newblock Basic principles of mixed virtual element methods.
\newblock {\em ESAIM Math. Model. Numer. Anal.}, 48(4):1227--1240, 2014.

\bibitem{Brezzi-Marini:2013}
F.~Brezzi and L.~D. Marini.
\newblock Virtual element methods for plate bending problems.
\newblock {\em Comput. Methods Appl. Mech. Engrg.}, 253:455--462, 2013.

\bibitem{Cahn:1961}
J.~W. Cahn.
\newblock On spinodal decomposition.
\newblock {\em Acta Metall.}, 9:795--801, 1961.

\bibitem{Cahn-Hilliard:1958}
J.~W. Cahn and J.~E. Hilliard.
\newblock Free energy of a nonuniform system. {I}. {I}nterfacial free energy.
\newblock {\em The Journal of Chemical Physics}, 28:258--0, 1958.

\bibitem{Cahn-Hilliard:1959}
J.~W. Cahn and J.~E. Hilliard.
\newblock Free energy of a nonuniform system. {III}. nucleation in a
  two-component incompressible fluid.
\newblock {\em The Journal of Chemical Physics}, 31:688--, 1959.

\bibitem{Cangiani-Georgoulis-Pryer-Sutton:2016}
A.~Cangiani, E.~H. Georgoulis, T.~Pryer, and O.~J. Sutton.
\newblock A posteriori error estimates for the virtual element method.
\newblock {\em Numer. Math.}, 137:857--893, 2017.

\bibitem{Cangiani-Gyrya-Manzini:2016}
A.~Cangiani, V.~Gyrya, and G.~Manzini.
\newblock The non-conforming virtual element method for the {S}tokes equations.
\newblock {\em SIAM J. Numer. Anal.}, 54(6):3411--3435, 2016.

\bibitem{Cangiani-Gyrya-Manzini-Sutton:2017:GBC:chbook}
A.~Cangiani, V.~Gyya, G.~Manzini, and Sutton. O.
\newblock Chapter 14: Virtual element methods for elliptic problems on
  polygonal meshes.
\newblock In K.~Hormann and N.~Sukumar, editors, {\em Generalized Barycentric
  Coordinates in Computer Graphics and Computational Mechanics}, pages 1--20.
  CRC Press, Taylor \& Francis Group, 2017.

\bibitem{Cangiani-Manzini-Russo-Sukumar:2015}
A.~Cangiani, G.~Manzini, A.~Russo, and N.~Sukumar.
\newblock Hourglass stabilization of the virtual element method.
\newblock {\em Internat. J. Numer. Methods Engrg.}, 102(3-4):404--436, 2015.

\bibitem{Cangiani-Manzini-Sutton:2017}
A.~Cangiani, G.~Manzini, and O.~Sutton.
\newblock Conforming and nonconforming virtual element methods for elliptic
  problems.
\newblock {\em IMA J. Numer. Anal.,}, 37:1317--1354, 2017.
\newblock (online August 2016).

\bibitem{Certik-Gardini-Manzini-Mascotto-Vacca:2020}
O.~Certik, F.~Gardini, G.~Manzini, L.~Mascotto, and G.~Vacca.
\newblock The p- and hp-versions of the virtual element method for elliptic
  eigenvalue problems.
\newblock {\em Comput. Math. Appl.}, 79(7):2035--2056, 2020.

\bibitem{Certik-Gardini-Manzini-Vacca:2018:ApplMath:journal}
O.~Certik, F.~Gardini, G.~Manzini, and G.~Vacca.
\newblock The virtual element method for eigenvalue problems with potential
  terms on polytopic meshes.
\newblock {\em Applications of Mathematics}, 63(3):333--365, 2018.

\bibitem{Chave-DiPietro-Marche-Pigeonneau:2016}
F.~Chave, D.~A. {Di~Pietro}, F.~Marche, and F.~Pigeonneau.
\newblock A hybrid high-order method for the {C}ahn-{H}illiard problem in mixed
  form.
\newblock {\em SIAM J. Numer. Anal.}, 54(3):1873--1898, 2016.

\bibitem{Chen-Shen:2013}
F.~Chen and J.~Shen.
\newblock Efficient energy stable schemes with spectral discretization in space
  for anisotropic {C}ahn-{H}illiard systems.
\newblock {\em Communications in Computational Physics}, 13(5):1189–1208,
  2013.

\bibitem{Chen-Huang:2020}
L.~Chen and X.~Huang.
\newblock Nonconforming virtual element method for {$2m$}th order partial
  differential equations in {$\mathbb{R}^n$}.
\newblock {\em Math. Comp.}, 89(324):1711--1744, 2020.

\bibitem{Chinosi-Marini:2016}
C.~Chinosi and L.~D. Marini.
\newblock Virtual element method for fourth order problems: {$L^2$}-estimates.
\newblock {\em Comput. Math. Appl.}, 72(8):1959--1967, 2016.

\bibitem{Dassi-Mascotto:2018}
F.~Dassi and L.~Mascotto.
\newblock Exploring high-order three dimensional virtual elements: bases and
  stabilizations.
\newblock {\em Comput. Math. Appl.}, 75(9):3379--3401, 2018.

\bibitem{DiPietro-Droniou-Manzini:2018}
D.~A. Di~Pietro, J.~Droniou, and G.~Manzini.
\newblock Discontinuous skeletal gradient discretisation methods on polytopal
  meshes.
\newblock {\em J. Comput. Phys.}, 355:397--425, 2018.

\bibitem{Elliott-French:1987}
C.~M. Elliott and D.~A. French.
\newblock Numerical studies of the {C}ahn-{H}illiard equation for phase
  separation.
\newblock {\em IMA J. Appl. Math.}, 38(2):97--128, 1987.

\bibitem{Elliott-French:1989}
C.~M. Elliott and D.~A. French.
\newblock A nonconforming finite-element method for the two-dimensional
  {C}ahn-{H}illiard equation.
\newblock {\em SIAM J. Numer. Anal.}, 26(4):884--903, 1989.

\bibitem{Elliott-French-Milner:1989}
C.~M. Elliott, D.~A. French, and F.~A. Milner.
\newblock A second-order splitting method for the {C}ahn-{H}illiard equation.
\newblock {\em Numer. Math.}, 54(5):575--590, 1989.

\bibitem{Elliott-Larsson:1992}
C.~M. Elliott and S.~Larsson.
\newblock Error estimates with smooth and nonsmooth data for a finite element
  method for the {C}ahn-{H}illiard equation.
\newblock {\em Math. Comp.}, 58(198):603--630, S33--S36, 1992.

\bibitem{Elliott-Zheng:1986}
C.~M. Elliott and Z.~Songmu.
\newblock On the {C}ahn-{H}illiard equation.
\newblock {\em Arch. Rational Mech. Anal.}, 96(4):339--357, 1986.

\bibitem{Faccioli-Maggio-Quarteroni-Tagliani:1996}
E.~Faccioli, F.~Maggio, A.~Quarteroni, and A.~Taghan.
\newblock Spectral‐domain decomposition methods for the solution of acoustic
  and elastic wave equations.
\newblock {\em The Leading Edge}, 61:1160--1174, 07 1996.
\newblock Faccioli1996.

\bibitem{Gallistl:2017}
D.~Gallistl.
\newblock Stable splitting of polyharmonic operators by generalized {S}tokes
  systems.
\newblock {\em Math. Comp.}, 86(308):2555--2577, 2017.

\bibitem{Gardini-Manzini-Vacca:2019:M2AN:journal}
F.~Gardini, G.~Manzini, and G.~Vacca.
\newblock The nonconforming virtual element method for eigenvalue problems.
\newblock {\em ESAIM Math. Model. Numer.}, 53:749--774, 2019.
\newblock Accepted for publication: 29 November 2018. DOI-DUMMY:
  10.1051/m2an/2018074.

\bibitem{Gazzola-Grunau-Sweers:1991}
F.~Gazzola, H.-C. Grunau, and G.~Sweers.
\newblock {\em Polyharmonic boundary value problems}, volume 1991 of {\em
  Lecture Notes in Mathematics}.
\newblock Springer-Verlag, Berlin, 2010.
\newblock Positivity preserving and nonlinear higher order elliptic equations
  in bounded domains.

\bibitem{Gomez-Calo-Bazilevs-Hughes:2008}
H.~G\'{o}mez, V.~M. Calo, Y.~Bazilevs, and T.~J.~R. Hughes.
\newblock Isogeometric analysis of the {C}ahn-{H}illiard phase-field model.
\newblock {\em Comput. Methods Appl. Mech. Engrg.}, 197(49-50):4333--4352,
  2008.

\bibitem{Grisvard:1985}
P.~Grisvard.
\newblock {\em Elliptic problems in nonsmooth domains}, volume~24 of {\em
  Monographs and Studies in Mathematics}.
\newblock Pitman (Advanced Publishing Program), Boston, MA, 1985.

\bibitem{Gudi-Neilan:2011}
T.~Gudi and M.~Neilan.
\newblock An interior penalty method for a sixth-order elliptic equation.
\newblock {\em IMA J. Numer. Anal.}, 31(4):1734--1753, 2011.

\bibitem{Kay-Styles-Suli:2009}
D.~Kay, V.~Styles, and E.~S\"{u}li.
\newblock Discontinuous {G}alerkin finite element approximation of the
  {C}ahn-{H}illiard equation with convection.
\newblock {\em SIAM J. Numer. Anal.}, 47(4):2660--2685, 2009.

\bibitem{Komatitsch-Tromp:1999}
D.~Komatitsch and J.~Tromp.
\newblock Introduction to the spectral element method for three-dimensional
  seismic wave propagation.
\newblock {\em Geophysical Journal International}, 139(3):806--822, 1999.

\bibitem{Korteweg:1901}
D.~J. Korteweg.
\newblock Sur la forme que prenent les \'equations du mouvements des fluides si
  l'on tient compte des forces capilaires caus\'ees par des variations de
  densit\'e consid\'erables mains continues et sur la th\'eorie de la
  capillarit\'e dans l'hypoth\'ese d'une varation continue de la densit\'e.
\newblock Arch. N\'eerl Sci. Exactes Nat. Ser. II, 1901.

\bibitem{Landau-Ginzburg:1965}
L.~D. Landau.
\newblock On the theory of superconductivity.
\newblock In D.~ter Haar, editor, {\em Collected papers of L. D. Landau}, pages
  546--568. Pergamon, 1965.

\bibitem{Lipnikov-Manzini:2014}
K.~Lipnikov and G.~Manzini.
\newblock A high-order mimetic method for unstructured polyhedral meshes.
\newblock {\em J. Comput. Phys.}, 272:360--385, 2014.

\bibitem{Lipnikov-Manzini-Shashkov:2014}
K.~Lipnikov, G.~Manzini, and M.~Shashkov.
\newblock Mimetic finite difference method.
\newblock {\em J. Comput. Phys.}, 257 -- Part B:1163--1227, 2014.

\bibitem{Xin-Zhangxin:2019}
X.~Liu and Z.~Chen.
\newblock A virtual element method for the {C}ahn-{H}illiard problem in mixed
  form.
\newblock {\em Appl. Math. Lett.}, 87:115--124, 2019.

\bibitem{Lovadina-Mora-Velasquez:2019}
C.~Lovadina, D.~Mora, and I.~Vel\'{a}squez.
\newblock A virtual element method for the von {K}\'arm\'an equations.
\newblock Technical report, Preprint CI2MA:2019-36, 2019.

\bibitem{Manzini-Lipnikov-Moulton-Shashkov:2017}
G.~Manzini, K.~Lipnikov, J.~D. Moulton, and M.~Shashkov.
\newblock Convergence analysis of the mimetic finite difference method for
  elliptic problems with staggered discretizations of diffusion coefficients.
\newblock {\em SIAM J. Numer. Anal.}, 55(6):2956--2981, 2017.

\bibitem{Manzini-Russo-Sukumar:2014}
G.~Manzini, A.~Russo, and N.~Sukumar.
\newblock New perspectives on polygonal and polyhedral finite element methods.
\newblock {\em Math. Models Methods Appl. Sci}, 24(8):1621--1663, 2014.

\bibitem{Mascotto:2018}
L.~Mascotto.
\newblock Ill-conditioning in the virtual element method: stabilizations and
  bases.
\newblock {\em Numer. Methods Partial Differential Equations},
  34(4):1258--1281, 2018.

\bibitem{Mora-Rivera-Rodriguez:2015}
D.~Mora, G.~Rivera, and R.~Rodr\'{i}guez.
\newblock A virtual element method for the {S}teklov eigenvalue problem.
\newblock {\em Math. Methods Appl. Sci.}, 25(08):1421--1445, 2015.

\bibitem{Mora-Rivera-Velasquez:2018}
D.~Mora, G.~Rivera, and I.~Vel\'{a}squez.
\newblock A virtual element method for the vibration problem of {K}irchhoff
  plates.
\newblock {\em ESAIM Math. Model. Numer. Anal.}, 52(4):1437--1456, 2018.

\bibitem{Mora-Velasquez:2018}
D.~Mora and I.~Vel\'{a}squez.
\newblock A virtual element method for the transmission eigenvalue problem.
\newblock {\em Math. Models Methods Appl. Sci.}, 28(14):2803--2831, 2018.

\bibitem{Mora-Velasquez:2020}
D.~Mora and I.~Vel\'{a}squez.
\newblock Virtual element for the buckling problem of {K}irchhoff-{L}ove
  plates.
\newblock {\em Comput. Methods Appl. Mech. Engrg.}, 360:112687, 22, 2020.

\bibitem{Park-Chi-Paulino:2019-CMAME}
K.~Park, H.~Chi, and G.~H. Paulino.
\newblock On nonconvex meshes for elastodynamics using virtual element methods
  with explicit time integration.
\newblock {\em Comput. Methods Appl. Mech. Engrg.}, 356:669--684, 2019.

\bibitem{Park-Chi-Paulino:2019-IJNME}
K.~Park, H.~Chi, and G.~H. Paulino.
\newblock Numerical recipes for elastodynamic virtual element methods with
  explicit time integration.
\newblock {\em Internat. J. Numer. Methods Engrg.}, 121(1):1--31, 2020.

\bibitem{Paulino-Gain:2015}
G.~H. Paulino and A.~L. Gain.
\newblock Bridging art and engineering using {E}scher-based virtual elements.
\newblock {\em Struct. and Multidisciplinary Optim.}, 51(4):867--883, 2015.

\bibitem{Perugia-Pietra-Russo:2016}
I.~Perugia, P.~Pietra, and A.~Russo.
\newblock A plane wave virtual element method for the {H}elmholtz problem.
\newblock {\em ESAIM Math. Model. Num.}, 50(3):783--808, 2016.

\bibitem{Quarteroni-Sacco-Saleri:2007}
A.~Quarteroni, R.~Sacco, and F.~Saleri.
\newblock {\em Numerical Mathematics}, volume Vol.~37 of {\em Texts in Applied
  Mathematics}.
\newblock Springer, 2007.

\bibitem{Raviart-Thomas:1983}
P.-A. Raviart and J.-M. Thomas.
\newblock {\em Introduction \`a l'analyse num\'erique des \'equations aux
  d\'eriv\'ees partielles}.
\newblock Collection Math\'ematiques Appliqu\'ees pour la Ma\^\i trise.
  [Collection of Applied Mathematics for the Master's Degree]. Masson, Paris,
  1983.

\bibitem{Riviere-Wheeler:2003}
B.~Rivi\`ere and M.~F. Wheeler.
\newblock Discontinuous finite element methods for acoustic and elastic wave
  problems.
\newblock In {\em Current trends in scientific computing ({X}i'an, 2002)},
  volume 329 of {\em Contemp. Math.}, pages 271--282. Amer. Math. Soc.,
  Providence, RI, 2003.

\bibitem{Rowlinson:1979}
J.~S. Rowlinson.
\newblock Translation of {J}. {D}. van der {W}aals' ``{T}he thermodynamic
  theory of capillarity under the hypothesis of a continuous variation of
  density''.
\newblock {\em J. Statist. Phys.}, 20(2):197--244, 1979.

\bibitem{Schedensack:2016}
M.~Schedensack.
\newblock A new discretization for {$m$}th-{L}aplace equations with arbitrary
  polynomial degrees.
\newblock {\em SIAM J. Numer. Anal.}, 54(4):2138--2162, 2016.

\bibitem{Torabi-Lowengrub-Voigt-Wise:2009}
S.~Torabi, J.~Lowengrub, A.~Voigt, and S.~Wise.
\newblock A new phase-field model for strongly anisotropic systems.
\newblock {\em Proc. R. Soc. Lond. Ser. A Math. Phys. Eng. Sci.},
  465(2105):1337--1359, 2009.
\newblock With supplementary material available online.

\bibitem{Vacca:2016}
G.~Vacca.
\newblock Virtual element methods for hyperbolic problems on polygonal meshes.
\newblock {\em Comput. Math. Appl.}, 74(5):882--898, 2017.

\bibitem{Wang-Xu:2013}
M.~Wang and J.~Xu.
\newblock Minimal finite element spaces for {$2m$}-th-order partial
  differential equations in {$R^n$}.
\newblock {\em Math. Comp.}, 82(281):25--43, 2013.

\bibitem{Wells-Kuhl-Garikipati:2006}
G.~N. Wells, E.~Kuhl, and K.~Garikipati.
\newblock A discontinuous {G}alerkin method for the {C}ahn-{H}illiard equation.
\newblock {\em J. Comput. Phys.}, 218(2):860--877, 2006.

\bibitem{Wriggers-Rust-Reddy:2016}
P.~Wriggers, W.~T. Rust, and B.~D. Reddy.
\newblock A virtual element method for contact.
\newblock {\em Comput. Mech.}, 58(6):1039--1050, 2016.

\bibitem{Zhao-Chen-Zhang:2016}
J.~Zhao, S.~Chen, and B.~Zhang.
\newblock The nonconforming virtual element method for plate bending problems.
\newblock {\em Math. Models Methods Appl. Sci.}, 26(9):1671--1687, 2016.

\bibitem{Zhao-Zhang-Chen-Mao:2018}
J.~Zhao, B.~Zhang, S.~Chen, and S.~Mao.
\newblock The {M}orley-type virtual element for plate bending problems.
\newblock {\em J. Sci. Comput.}, 76(1):610--629, 2018.

\end{thebibliography}

\end{document}